\documentclass[graybox,envcountchap,sectrefs,envcountsame]{svmono} 
\usepackage{amsfonts}
\usepackage{amsmath}
\usepackage{amssymb,psfrag}
\usepackage{enumitem}
\usepackage{bbm}
\usepackage{tikz}
\usepackage{tkz-graph}

\setlist[itemize]{leftmargin=*}
\setlist[enumerate]{leftmargin=*,label=(\roman*),ref=(\roman*)}

\spnewtheorem{observation}[theorem]{Observation}{\bf}{\it}

\usepackage[numeric,initials,nobysame,msc-links,abbrev]{amsrefs}
\renewcommand{\eprint}[1]{\href{https://arxiv.org/abs/#1}{arXiv:#1}}
\newcommand{\pageafter}[1]{#1~pp.}
\BibSpec{article}{%
+{} {\PrintAuthors} {author}
+{,} { \textit} {title}
+{.} { } {part}
+{:} { \textit} {subtitle}
+{,} { \PrintContributions} {contribution}
+{.} { \PrintPartials} {partial}
+{,} { } {journal}
+{} { \textbf} {volume}
+{} { \PrintDatePV} {date}
+{,} { \issuetext} {number}
+{,} { \pageafter} {pages}
+{,} { } {status}
+{,} { \PrintDOI} {doi}
+{,} { available at \eprint} {eprint}
+{} { \parenthesize} {language}
+{} { \PrintTranslation} {translation}
+{;} { \PrintReprint} {reprint}
+{.} { } {note}
+{.} {} {transition}
+{} {\SentenceSpace \PrintReviews} {review}
}
\BibSpec{collection.article}{%
+{} {\PrintAuthors} {author}
+{,} { \textit} {title}
+{.} { } {part}
+{:} { \textit} {subtitle}
+{,} { \PrintContributions} {contribution}
+{,} { \PrintConference} {conference}
+{} {\PrintBook} {book}
+{,} { } {booktitle}
+{,} { \PrintDateB} {date}
+{,} { \pageafter} {pages}
+{,} { } {status}
+{,} { \PrintDOI} {doi}
+{,} { available at \eprint} {eprint}
+{} { \parenthesize} {language}
+{} { \PrintTranslation} {translation}
+{;} { \PrintReprint} {reprint}
+{.} { } {note}
+{.} {} {transition}
+{} {\SentenceSpace \PrintReviews} {review}
}

\usepackage{verbatim}
\usetikzlibrary{arrows}
\usetikzlibrary[patterns]
\usetikzlibrary{shapes.misc}
\usetikzlibrary{decorations.pathreplacing}
\usetikzlibrary{arrows.meta}
\usepackage{subcaption}
\usepackage{array}
\usepackage{multirow}

\usepackage{hyperref}

\tikzset{cross/.style={cross out, draw=black, minimum size=2*(#1-\pgflinewidth), inner sep=0pt, outer sep=0pt},
cross/.default={1pt}}

\usepackage{newtxtext}
\usepackage{newtxmath}

\newcommand{\Var}{\operatorname{Var}} 
\newcommand{\Ent}{\operatorname{Ent}} 

\newcommand{\cA}{\ensuremath{\mathcal A}}
\newcommand{\cB}{\ensuremath{\mathcal B}}
\newcommand{\cC}{\ensuremath{\mathcal C}}
\newcommand{\cD}{\ensuremath{\mathcal D}}
\newcommand{\cE}{\ensuremath{\mathcal E}}

\newcommand{\cG}{\ensuremath{\mathcal G}}
\newcommand{\cH}{\ensuremath{\mathcal H}}
\newcommand{\cI}{\ensuremath{\mathcal I}}

\newcommand{\cK}{\ensuremath{\mathcal K}}
\newcommand{\cL}{\ensuremath{\mathcal L}}

\newcommand{\cS}{\ensuremath{\mathcal S}}
\newcommand{\cT}{\ensuremath{\mathcal T}}
\newcommand{\cU}{\ensuremath{\mathcal U}}

\newcommand{\cX}{\ensuremath{\mathcal X}}

\newcommand{\bbE}{{\ensuremath{\mathbb E}} }
\newcommand{\bbF}{{\ensuremath{\mathbb F}} }

\newcommand{\bbH}{{\ensuremath{\mathbb H}} }

\newcommand{\bbN}{{\ensuremath{\mathbb N}} }

\newcommand{\bbP}{{\ensuremath{\mathbb P}} }
\newcommand{\bbQ}{{\ensuremath{\mathbb Q}} }
\newcommand{\bbR}{{\ensuremath{\mathbb R}} }

\newcommand{\bbX}{{\ensuremath{\mathbb X}} }

\newcommand{\bbZ}{{\ensuremath{\mathbb Z}} }

\newcommand{\sB}{\ensuremath{\mathscr B}}

\newcommand{\bone}{{\ensuremath{\mathbf 1}} }
\newcommand{\bzero}{{\ensuremath{\mathbf 0}} }

\newcommand{\be}{{\ensuremath{\mathbf e}} }

\newcommand{\bu}{{\ensuremath{\mathbf u}} }
\newcommand{\bv}{{\ensuremath{\mathbf v}} }

\newcommand{\bx}{{\ensuremath{\mathbf x}} }
\newcommand{\by}{{\ensuremath{\mathbf y}} }
\newcommand{\bz}{{\ensuremath{\mathbf z}} }

\newcommand{\var}{\operatorname{Var}}
\newcommand{\<}{\langle}
\renewcommand{\>}{\rangle}
\newcommand{\diam}{{\ensuremath{\mathrm{diam}}} }

\newcommand{\1}{{\ensuremath{\mathbbm{1}}} }
\newcommand{\trel}{T_{\mathrm{rel}}}

\newcommand{\qc}{{\ensuremath{q_{\mathrm{c}}}} }
\newcommand{\qct}{{\ensuremath{\tilde q_{\mathrm{c}}}} }
\newcommand{\tbp}{{\ensuremath{\tau_0^{\mathrm{BP}}}} }

\renewcommand{\leq}{\leqslant}
\renewcommand{\geq}{\geqslant}
\renewcommand{\le}{\leqslant}
\renewcommand{\ge}{\geqslant}
\DeclareDocumentCommand \to { o o } {%
  \IfNoValueTF {#1} {\IfNoValueTF{#2}{\rightarrow}{\xrightarrow[#2]}}%
{\IfNoValueTF{#2}{\xrightarrow{#1}}{\xrightarrow[#2]{#1}}}}

\begin{document}
\frontmatter
\title{Kinetically constrained models}
\author{Ivailo Hartarsky and Cristina Toninelli}
\maketitle

{\dedication{To our friend, mentor and collaborator Fabio Martinelli, \newline
without whom this book would have ended here.}}

\tableofcontents
\preface

This manuscript focuses on Kinetically Constrained Models (KCM), a topic which  lies at the intersection between probability and statistical mechanics. KCM are a class of Markov processes. They belong to the larger class of interacting particle systems with stochastic dynamics on discrete lattices. KCM were introduced in the physics literature in the 1980's to model the liquid-glass transition, a longstanding open problem in condensed matter physics. The key feature of KCM is that the update at a given lattice site can occur only  if the configuration verifies a  kinetic constraint requiring that there are no particles in a suitable neighbourhood. Extensive numerical simulations  indicate that  KCM display a remarkable behavior  typical of glassy systems. Therefore,  they have  been the subject of several investigations in  the last forty years    with the aim of providing a deeper  understanding of the liquid-glass transition and of  more general jamming transitions.

Mathematically, KCM pose very challenging and interesting problems. In fact, the presence of the constraints induces non-attractiveness,  the occurrence of several invariant measures, and  the failure of many powerful tools to analyze relaxation to equilibrium (coercive inequalities, coupling, censoring\dots). Remarkably, the degeneracy of the rates caused by the constraints is not a mere technical obstacle which prevents using the classic tools. Indeed, the behavior of KCM is qualitatively different  from that of interacting particle systems without constraints. Peculiar features include anomalously long mixing times, aging, singularities in the dynamical large deviation function, dynamical heterogeneities, and atypical ergodicity breaking transitions corresponding to the emergence of a large variety of amorphous structures.  All in all, we can definitely say that KCM open a new chapter in the well established field of interacting particle systems.
 
Major progress has been made in the last  twenty 
years towards a full and rigorous understanding of the large time behavior of  KCM at stationarity. We present these results, illustrating  both the high level ideas and some  novel technical  tools that have been devised to deal with the presence of constraints and with the lack of attractiveness. On the way, we unveil some remarkable connections of KCM with other mathematical subjects, in particular with  bootstrap percolation cellular automata. 
We also present a choice of open problems concerning particularly the out of equilibrium dynamics. Indeed, despite some achievements, robust tools to analyse KCM in this regime are still lacking and several  beautiful questions remain open, even for simple choices of constraints.

This book aims at being accessible to both mathematicians and physicists. Hopefully it will be a useful tool to reinforce the bridge between the two communities which, in our opinion, have still much to learn from each other on KCM and glassy dynamics.

\section*{Outline}
The content of the manuscript is as follows.
\begin{itemize}
    \item In Chapter \ref{chap:intro} we provide the physics background and motivation for studying KCM.
    \item In Chapter \ref{chap:models} we formally introduce KCM along with the relevant notation and key quantities of interest. It may be viewed as defining the scope of the manuscript.
    \item In Chapter \ref{chap:BP} we discuss deterministic monotone cellular automata known as bootstrap percolation and their fundamental relation to KCM.
    \item In Chapter \ref{chap:1d} we explore KCM in one dimension and introduce some basic tools, notably the bisection-constrained method. We focus particularly on the
    Fredrickson--Andersen 1-spin facilitated model (FA-1f) and on the East model, which not only serve as a warm-up for more advanced models, but also as a tool for their study.
    \item In Chapter \ref{chap:FA2f} we consider the Fredrickson--Andersen 2-spin facilitated model in 2 dimensions. We develop progressively more sophisticated tools for its study, culminating with determining its sharp asymptotic behaviour at low temperature. These tools, which are flexible enough to be generalised to treat other models, include a robust long range Poincar\'e inequality and a very flexible multi-scale renormalisation tool, the Matryoska dolls. 
    \item In Chapter \ref{chap:universality} we examine the universality theory for KCM in one and two dimensions. It further elaborates our techniques and establishes a detailed map of the domain.
    \item In Chapter \ref{chap:out} we turn our attention to results on KCM out of equilibrium. Convergence to equilibrium and mixing times are investigated, using a set of tools completely separate from previous chapters.
    \item In Chapter \ref{chap:other} we briefly mention several settings, other than the one of Chapter~\ref{chap:models}, in which KCM have been studied. We also mention some closely related models, and provide more detailed references for the interested reader.
\end{itemize}
Dependencies between different chapters are shown in the next diagram.
\begin{center}
\begin{tikzpicture}[x=0.20\textwidth,y=0.1\textwidth]
\GraphInit[vstyle=Normal]
  \SetVertexNormal[Shape = rectangle,
  					LineWidth = 1pt]\tikzset{VertexStyle/.append style = {outer sep=1pt}}
  \SetVertexNormal[Shape      = circle,
                   LineWidth  = 1pt]\tikzset{VertexStyle/.append style = {outer sep=1pt}}
  \Vertex[x=0,y=2,L=\ref{chap:intro}]{intro}
  \Vertex[x=2,y=1,L=\ref{chap:1d}]{1d}
  \Vertex[x=1,y=1,L=\ref{chap:BP}]{BP}
  \Vertex[x=3,y=1,L=\ref{chap:FA2f}]{FA2f}
  \Vertex[x=1,y=2,L=\ref{chap:other}]{other}
  \Vertex[x=0,y=1,L=\ref{chap:models}]{models}
  \Vertex[x=2,y=2,L=\ref{chap:out}]{out}
  \Vertex[x=4,y=1,L=\ref{chap:universality}]{universality}
  
  \tikzset{EdgeStyle/.style = {->}}
  \Edge(1d)(FA2f)
  \Edge(FA2f)(universality)
  \Edge(BP)(out)
  \Edge(BP)(1d)
  \Edge(models)(BP)
\end{tikzpicture}
\end{center}
Consequently, Chapters~\ref{chap:intro} and \ref{chap:other} can be regarded as optional general knowledge. Chapters~\ref{chap:models} and~\ref{chap:BP} are indispensable core material. A graduate course on the subject could cover these two chapters and a selection of Chapters~\ref{chap:1d} and/or~\ref{chap:out}, which both introduce a large variety of techniques in an accessible setting. The remaining Chapters~\ref{chap:FA2f} and~\ref{chap:universality} are intended for a more expert audience, particularly for newcomers to the field, who have already covered the basics, but need some background and intuition before delving into the details of specific papers.

The more basic material (Chapters~\ref{chap:intro}-\ref{chap:1d}) is covered including full proofs or detailed sketches and featuring exercises to help assimilating the content. Subsequent chapters are less detailed and often refer to original papers for technical details. We apologise for the inevitable inaccuracies due to favouring simplicity over technical completeness. Indeed, we aim at highlighting the heuristic ideas and the guiding lines behind each method and result.

We have tried to keep the presentation as self-contained as possible, but there are some prerequisites. We do not assume, but hope for some familiarity with the basics of standard textbooks in the field such as \cite{Levin09} on Markov chains and \cites{Liggett99,Liggett05} on interacting particle systems, while some of the contents of an undergraduate course in probability theory, as covered for instance in \cites{Durrett19,Legall22,Feller68,Dudley02}, will be used without notice. While it is possible to only refer to these books as needed, it may be a good investment to first acquire some superficial experience with their content, which is excellent to have in any case.

\section*{Acknowledgements}

We thank Damiano de Gaspari, Fabio Martinelli, Quentin Moulard and Fabio Toninelli for helpful discussions. We thank Giulio Biroli, Oriane Blondel, Laure Mar\^ech\'e, Fabio Martinelli, Assaf Shapira, R\'eka Szab\'o for insightful comments and corrections on the presentation and Filippo Nuccio for further support. I.H.~was supported by the Austrian Science Fund (FWF): P35428-N. Most of his work on this manuscript was done in 2023 and 2024, when he was affiliated with TU Wien, whose hospitality is gratefully acknowledged. C.T.~was supported by ERC through Grant 680275 'MALIG'.

\mainmatter

	\chapter{Background from physics}

\label{chap:intro}
 
\abstract{In this chapter we discuss the physics background. We  recall the basic phenomenology of the liquid-glass transition and more general jamming transitions, and explain the role of KCM as models for glassy dynamics.}

\section{The liquid-glass transition}
\label{sec:vetri}

From the point of view of statistical physics, a key motivation behind the study of
KCM is the effort to understand the  \emph{liquid-glass transition}. Glass is widely present in our daily life: it is a very versatile material, easily produced and manipulated on an industrial scale by cooling different liquid mixtures (e.g.\ silica, sodium carbonate and calcium oxide). And yet a microscopic understanding of this state of matter (which, according to archaeological findings in Egypt and Eastern Mesopotamia, people have been manufacturing since 3000 B.C.) and of how glass forms is still an open challenge for condensed matter physicists (see e.g.\ the review
\cite{Arceri21}). In 1995,
Nobel prize Philip W.\ Anderson \cite{Anderson95} wrote: \emph{``The deepest and most interesting unsolved problem in solid state theory is probably the theory of the nature of glass and the glass transition.''}
He added, \emph{``This could be the next breakthrough in the coming decade.''} Thirty
years later, physicists still disagree about the nature of glass and on how it forms. 

At the heart of this puzzle lies the intriguing fact that the glasses display  properties of both solids and liquids. In fact, we could  either say ``glass is an extremely viscous liquid that does not flow'' or ``glass is an unstructured, amorphous 
solid''. Indeed, despite its macroscopic rigidity, the microscopic structure of a glass has the same disordered arrangement of molecules as a liquid. In other words,
based on a single snapshot, liquid and glass are essentially indistinguishable.

This lack of order might seem in contrast with the thermodynamics paradigm that, when a liquid is cooled below its melting temperature, an ordered  structure should form and
the liquid should become a crystal. The secret is to
perform the cooling sufficiently fast. This way the nucleation of the crystal is prevented and the liquid enters a long-lived metastable state, the \emph{super-cooled liquid} phase. 
Very roughly speaking, the liquid-crystal transition is avoided because molecules do not have enough time to reorganise and form the ordered crystal structure.
The molecules move slower and slower forming a thick syrup and eventually they get trapped in the structureless  glass state.
Since  the nucleation time of the stable crystal structure is out of reach for any reasonable experiment, though the state of super-cooled liquid is not thermodynamically stable, for all practical purposes it can be considered as an equilibrium system. 
In particular, one can define a relaxation time  (and measure it experimentally via  viscosity) and establish   fluctuation-dissipation relations connecting response to an external driving force and correlations functions. Essentially, we can forget about the crystal and just focus on the super-cooled phase.

\begin{figure}
\begin{minipage}{\textwidth}
\leftfigure[c]{\includegraphics[width=0.4\textwidth]{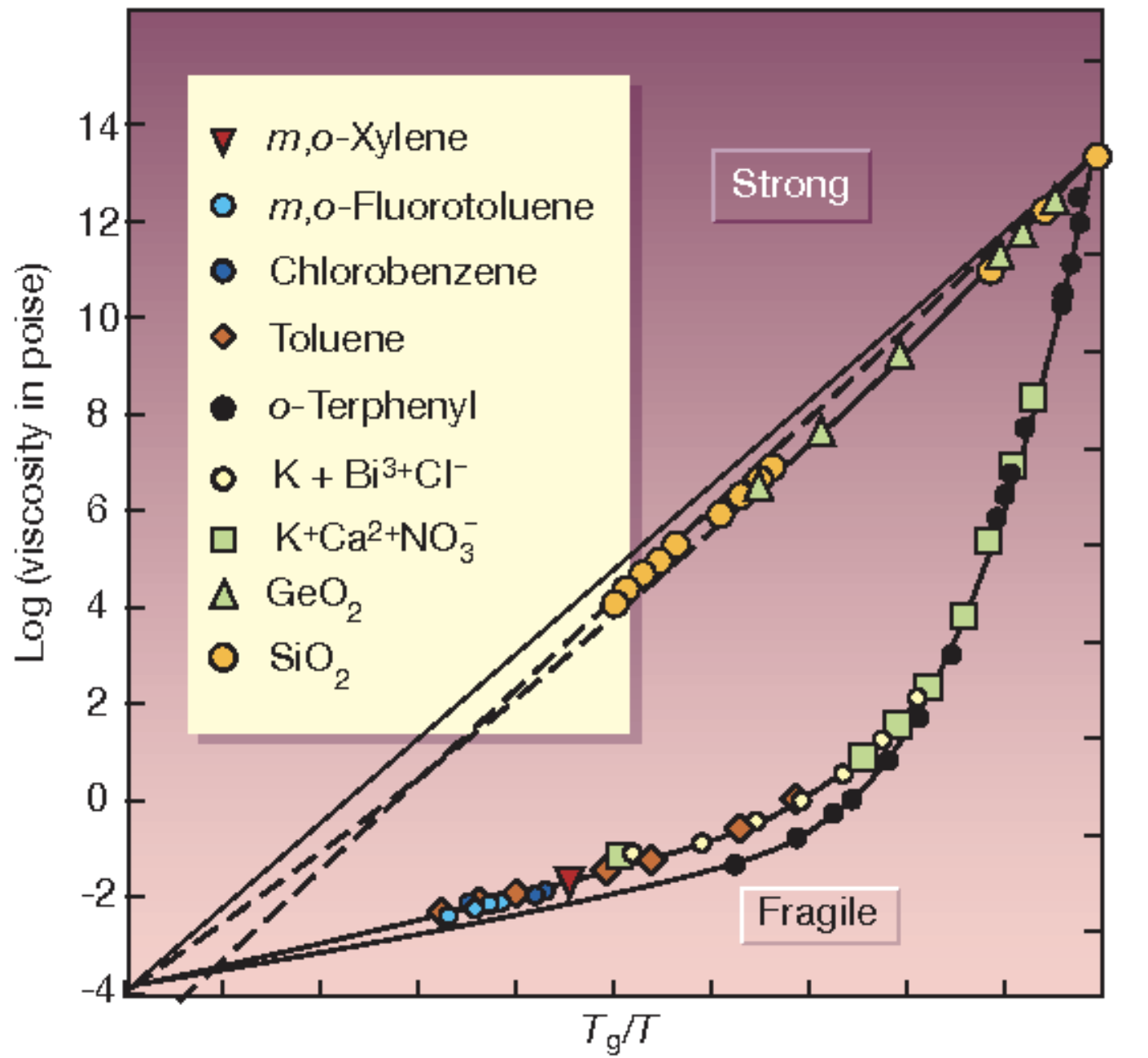}}\hspace{\fill}
\rightfigure[c]{\includegraphics[width=0.4\textwidth]{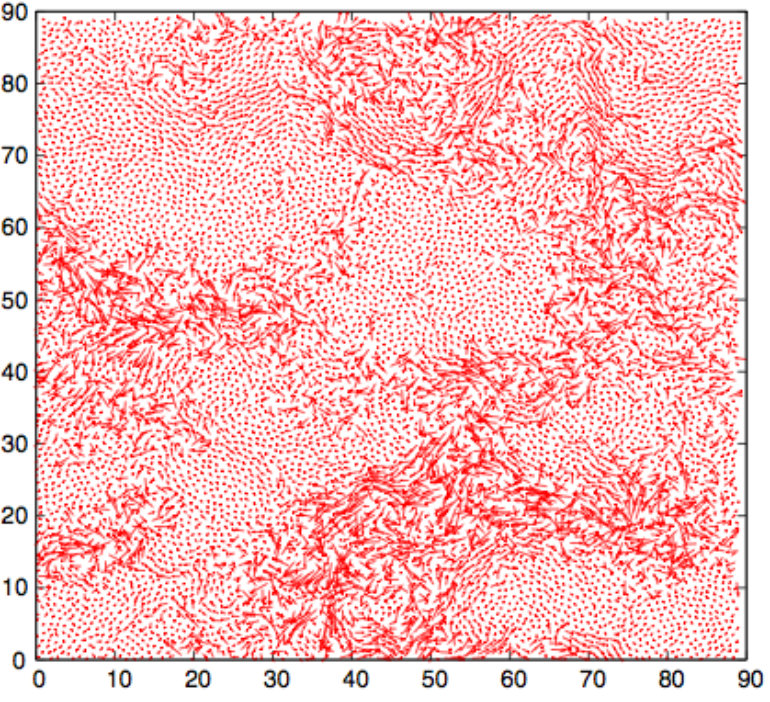}}
\end{minipage}
\leftcaption{\label{fig:vetri:1}Logarithm of the 
viscosity of several glass-forming
liquids plotted against the inverse temperature.  Here temperature is rescaled by  the empirical glass transition temperature,  $T_g$, defined  as the value at which the  viscosity  equals $10^{12}Pa\cdot s$. 
Reprinted from \cite{Debenedetti01} with permission. Copyright © 2001, Springer Nature Limited.
}
\rightcaption{\label{fig:vetri:2}Spatial map of single particle displacements in a molecular dynamics
simulation of a  super-cooled liquid in two dimensions. Arrows show the displacement of each particle
in a trajectory of length comparable to the relaxation
time. Reprinted from \cite{Berthier11} with permission. Copyright © 2011 American Physical Society.}
\end{figure}

A common feature of super-cooled liquids  around the glass transition is the sharp  slowdown of dynamics. As shown in Figure \ref{fig:vetri:1}, viscosity can increase by over 14 orders of magnitude upon a small decrease in temperature. 
It also highlights the fact that super-cooled liquids can be classified into two groups: \emph{strong} and  \emph{fragile}. If we let $\eta$ be the viscosity, $T$ be the temperature and we define the activation energy as $E:=T\log\eta$, strong liquids are characterised
by a temperature-independent $E$, while for fragile liquids $E$  increases as $T$ decreases. The corresponding scaling forms for $\eta$  are called \emph{Arrhenius} and \emph{super-Arrhenius}, respectively.

This dramatic growth of time scales\footnote{For example, for fragile glass-forming liquids time scales at the melting temperature are typically of the order of the picosecond  (which is also roughly the time scale of molecular motion), and are of the order of 100 seconds when temperature equals  $2/3$ of the melting temperature.} is related to the fact that when the temperature is lowered, the density increases: molecules tend to obstruct each other, blocked structures may arise, and  motion becomes very cooperative\footnote{{That is, a growing number of particles need to move in a coordinated way in order for relaxation to occur.}}. A key experimental observation of the cooperative motion is the fact that when a glass cools, the molecules do not slow down uniformly.  There is indeed a clear coexistence of fast and slow regions. This phenomenon is called  \emph{dynamical heterogeneity}: some regions of the liquid jam, while in other regions molecules continue to move around \cites{Berthier11a,Ediger00}. 
Thus, even if a change of structure does not occur when the glass is formed, an
underlying dynamical phase transition separating slow and fast trajectories seems to occur  (see Figure \ref{fig:vetri:2}). 
An indirect experimental evidence of dynamical heterogeneity is the decoupling of self-diffusion coefficient $D_s$ and viscosity $\eta$: super-cooled liquids violate the phenomenological Stokes--Einstein relation, $D_s\eta/T=\mathrm{constant}$, which holds in homogeneous liquids. 
Though both $D_s^{-1}$ and $\eta$
increase when the temperature is lowered, the former does not increase as fast as the latter. This leads to $D_s\eta/T$ increasing
by 2-3 orders of magnitude
when temperature is decreased towards the glass transition temperature (see e.g.\ \cite{Mapes06}). The reason why this decoupling of self-diffusion and viscosity is interpreted as a hallmark 
of dynamical heterogeneity is that the diffusion of the tracer  particle (see Section~\ref{sec:tracer} for a more formal description) should be dominated by the fastest regions and the structural relaxation (measured by the viscosity) by the slowest regions.
 
Despite a great deal of experimental and theoretical investigation, a complete understanding of this behaviour and of other
peculiar phenomena occurring in the vicinity of the glass transition (aging, hysteresis, rejuvenation, anomalous transport phenomena \dots) is still far out of reach. None of the numerous theories covers all the above phenomenology  and a common consensus  around ``the'' theory of the glass transition is still lacking in the physics community (see \cite{Arceri21} for a  review of various theories).  A central theoretical difficulty is the fact that from the point of view of critical phenomena the situation is very peculiar: the liquid/glass transition displays a \emph{mixed character}. Indeed, diverging time and length scales (typical of second order phase transitions), are accompanied by a discontinuous order parameter (typical of first order transitions). The jump of the order parameter corresponds to the discontinuous emergence of an amorphous density profile.
Furthermore, both from the experimental and the theoretical point of view, the  degeneracy of ground states complicates the problem. Thus, everybody agrees at least on one point: this is certainly not a standard type of ergodicity breaking transition!

The active research on the glass transition is also enhanced by the fact that a dynamical arrest towards an amorphous state displaying similar properties occurs in a large variety of systems upon tuning a proper external parameter (see e.g.\ \cite{Arceri21}*{Sec. IV}).
These  phenomena, that are generically dubbed  \emph{jamming transitions}, occur for several materials: grains in powders (granular media),  emulsions, foams, colloidal suspensions,  polymers, plastics, ceramics, etc. Last but not least, understanding the glass transition  would probably yield to  novel theoretical and numerical tools that could be useful in other fields  of science handling  
systems displaying a non trivial global collective behavior. This is the vast realm of systems that goes under the general name of 
 \emph{complex systems}.

\section{KCM as models for glassy dynamics}

Kinetically constrained models (KCM) are toy models for the liquid-glass and jamming transitions. They rely on the idea that these are dynamical phenomena in which static interactions play a minor role. The kinetic constraints are therefore devised to mimic the mechanism of local caging  which slows  the dynamics down at low temperature or high density. Originally motivated by free volume theories \cite{Glarum60}, KCM have been promoted in the last decades as a paradigm model for the so called {\emph {dynamical facilitation theory}} for the glass transition (see e.g.\ \cites{Biroli13,Chacko24}).
Indeed, despite their simplicity and their trivial statics, KCM display many key dynamical features of real materials that undergo glass or jamming transitions: anomalous ergodicity breaking transitions, percolation of blocked structures, dynamical arrest, non-trivial spatio-temporal fluctuations,  dynamical heterogeneities, aging\dots
Furthermore, depending on the choice of the constraints, they feature either a super-Arrhenius or an Arrhenius behavior for the relaxation time, thus we recover both the behavior of fragile and strong supercooled liquids.

On the other hand, a major criticism to KCM is  that
a convincing derivation of these toy models via a coarse graining from realistic molecular models of liquids is missing (though some attempts in this direction have been recently performed both in experiments of granular glasses \cites{Candelier09,Candelier10} and in numerical simulations of super-cooled liquids \cite{Ozawa23}).
In particular, it is not clear how one should identify at the molecular level the facilitating (empty) sites.
We refer the reader to  \cites{Biroli13,Arceri21,Ritort03,Garrahan11} for further comments on successes and limitations of  KCM as toy models for real glass forming liquids as well as for the illustration of alternative theories. What we can confidently say, adopting a sentence from \cite{Arceri21}, is that KCM  have been influential and very instructive in order to develop a theoretical understanding of glassy phenomena.

Regarding derivations of KCM, we should also mention a different class of models which have been introduced to prove that that kinetic constraints can emerge spontaneously from static interactions, the so-called plaquette models. These are spin models with the usual Glauber dynamics reversible w.r.t.\ a Gibbs measure corresponding to a particular Hamiltonian $\cH$. For certain choices of $\cH$, the relaxation at low temperature is dominated by the motion of  ``excited'' plaquettes. These excitations act as a source of mobility since the energetic barrier to flip a spin is smaller in their vicinity. Thus, at low temperature, their dynamics can be somehow mapped to a KCM. We will return to plaquette models in Section~\ref{sec:plaquettes}.

We also mention that
in recent years there has been an increasing amount of work on the quantum versions of KCM that have been proposed in the study of  Rydberg atoms \cite{Turner18} and as models for quantum many-body localization \cite{Garrahan18}. For example, in \cite{Pancotti20}, the quantum version of the East KCM has been analyzed and the occurrence of a first-order quantum transition at which the ground state becomes exponentially localized (with a consequent slowdown of dynamics) has been detected.


\chapter{Models}

\label{chap:models}
\abstract{
In this chapter we provide a general definition of KCM and associated quantities of interest. We introduce the choices of constraints corresponding to some of the most studied models, including East, Friedrickson-Andersen $j$-spin facilitated (FA-$j$f), Duarte, North-East (NE) and Spiral models.}

\section{Notation}\label{Section:Setting}
We denote the sets of nonnegative integers, integers and reals by $\bbN$, $\bbZ$ and $\bbR$ respectively. The models that we consider in this manuscript are most often defined on the infinite integer lattice $\mathbb Z^d$ or a subset thereof. We denote by $\vec x=(x_1,\dots x_d)$ the \emph{sites} of $\mathbb Z^d$, by $\vec e_1,\dots,\vec e_d$ its \emph{canonical basis vectors} and by $\Omega=\{0,1\}^{\mathbb Z^d}$ the \emph{configuration space}. Elements of $\Omega$ are called \emph{configurations} and denoted by Greek letters $\sigma,\eta,\omega,\dots$. For a configuration $\omega\in\Omega$ and a site $\vec x\in\bbZ^d$, $\omega_{\vec x}$ denotes the \emph{occupancy variable} (the value) of $\omega$ at $\vec x$. When $\omega_{\vec x}=1$, we say that site $\vec x$ is \emph{occupied}. When $\omega_{\vec x}=0$ we say that $\vec x$ is \emph{empty}. For $\omega\in\Omega$, we write $|\omega|=|\{\bx\in\bbZ^d:\omega_\bx=0\}|$ for the number of empty sites in $\omega$.

KCM are endowed with a parameter $q\in(0,1]$ called \emph{vacancy density} and corresponding to the inverse temperature $\beta$ (in the right units) via the relation
\begin{equation}\label{eq:temperature}
q=1/(1+e^{\beta}).\end{equation}
In particular, the limit $q\to 0$ corresponds to the zero temperature limit, $\beta\to\infty$. Given $q\in[0,1]$, we denote the product Bernoulli$(1-q)$ measure on $\Omega$, under which a site is empty with probability $q$, by $\mu_q$ or simply $\mu$, when $q$ is clear from the context. 

The mean with respect to a measure $\nu$ on $\Omega$ of a function $f:\Omega\to\bbR$ is denoted by $\nu(f)$, while its variance is denoted by $\var_\nu(f)$ or simply $\var(f)$ when $\nu=\mu$. Functions $f$ we are interested in only take a finite number of values, so no integrability issues arise.

We sometimes work on a subset $\Lambda\subset\bbZ^d$ of the lattice. Correspondingly, we set $\Omega_{\Lambda}=\{0,1\}^{\Lambda}$. The \emph{restriction} of a configuration $\omega\in\Omega$  to $\Lambda$ is denoted by $\omega_\Lambda\in\Omega_\Lambda$. For any measure $\nu$ that is the product of a measure on ${\Omega_\Lambda}$ and ${\Omega_{\Lambda^c}}$, we similarly denote by $\nu_{\Lambda}$ the restriction to $\Lambda$. When $\Lambda=\{\vec x\}$ is a singleton, we simply write $\nu_{\vec x}$ for $\nu_{\Lambda}$. We write $\var_\Lambda$ for $\var_{\mu_\Lambda}$, that is the variance with respect to the occupation variables in $\Lambda$. Given disjoint sets $\Lambda_1,\Lambda_2\subset \bbZ^d$, and $\omega^{(1)}\in \Omega_{\Lambda_1}$ and $\omega^{(2)}\in\Omega_{\Lambda_2}$, we write $\omega^{(1)}\cdot\omega^{(2)}\in \Omega_{\Lambda_1\cup\Lambda_2}$ for the configuration such that
\begin{equation}
\label{eq:def:concatenation}(\omega^{(1)}\cdot\omega^{(2)})_{\vec x}=\begin{cases}
\omega^{(1)}_{\vec x}&\vec x\in\Lambda_1,\\
\omega^{(2)}_{\vec x}&\vec x\in\Lambda_2.
\end{cases}\end{equation}
We denote the fully occupied (resp.\ empty) configurations by $\bone$ (resp.\ $\bzero$) and omit the domain $\Lambda$ in $\bone_\Lambda$ (resp.\ $\bzero_\Lambda$), when it is clear from the context.

\section{Update families}
\label{sec:rules}
KCM are Glauber type Markov processes on $\Omega$ (or $\Omega_{\Lambda}$) reversible w.r.t.\ $\mu$ ($\mu_{\Lambda}$, respectively).  We give the general definition introduced in \cite{Cancrini08} to cover all the models studied in physics. Each KCM is characterised by its \emph{update family}, namely a finite collection $\mathcal U=\{U_1,\dots,U_m\}$ of finite subsets of $\mathbb Z^d\setminus \{0\}$ called \emph{update rules}. Given a vertex ${\vec x}\in \bbZ^d$, we will say that the \emph{constraint} at $\vec x$ is satisfied by the configuration $\omega\in\Omega$ if at least one of the update rule translated at ${\vec x}$ is completely empty, namely,
\begin{equation}
c_{{\vec x}}(\omega)=c_{\vec x}^{\cU}(\omega)=
\begin{cases}
1&\exists U\in\cU, \forall \vec u\in U,\omega_{{\vec x}+\vec u}=0,\\
0&\text{otherwise}.
\end{cases}\label{eq:def:cx}
\end{equation}
Observe that constraints $c_{\vec x}$ are non-increasing w.r.t.\ the product partial order on configurations in $\Omega$ given by $\omega\le\omega'$ if $\omega(\vec y)\le \omega'(\vec y)$ for all $\vec y\in\bbZ^d$. Inversely, any non-increasing function $c_0:\Omega\to\{0,1\}$ depending only on the restriction of the configuration to a finite subset of $\bbZ^d\setminus\{0\}$ can be written in the form $c_0^\cU$ for some update family $\cU$. Correspondingly, there is a natural partial order on update families: $\cU_1\le \cU_2$, if $c_0^{\cU_1}(\omega)\le c_0^{\cU_2}(\omega)$ for every $\omega\in\Omega$.\footnote{Note that since we only care about the update family through the constraint it induces, we identify update families yielding the same constraints, e.g.\ $\{\{1\}\}$ and $\{\{1\},\{1,2\}\}$ in one dimension.} The maximal and minimal update families correspond to $\cU=\varnothing$ (i.e.\ $c_0^\cU\equiv0$) and $\cU=\{\varnothing\}$ (i.e.\ $c_0^\cU\equiv1$). Since the goal of KCM is to explore the effect of constraints on the dynamics, we discard these cases, by systematically assuming update families and update rules to be nonempty. Moreover, we usually drop $\cU$ from the notation, as $\cU$ is fixed or arbitrary.

\begin{figure}
\centering
\begin{subfigure}{0.3\textwidth}
\begin{center}
\begin{tikzpicture}
\draw (-1,-1.5) circle (1 mm) ;
\draw (0,-1) circle (1 mm);
\draw[step=0.5cm,gray,very thin](-2,-2)grid(-1,-1);
\draw[step=0.5cm,gray,very thin](-0.5,-2)grid(0.5,-1);
\draw (-1.5cm,-1.5cm) node[cross=3pt,rotate=0] {};
\draw (0,-1.5cm) node[cross=3pt,rotate=0] {};
\draw  (0.5 cm,1.8cm) node {};\end{tikzpicture}
\caption{East\label{fig:East}}
\end{center}
\end{subfigure}
\begin{subfigure}{0.3\textwidth}
\begin{center}
\begin{tikzpicture}
\draw (-1.5,-1.5) circle (1 mm) ;
\draw (0.5,-2) circle (1 mm);
\draw (0.5,0.5) circle (1 mm);
\draw (-0.5,0) circle (1 mm);
\draw[step=0.5cm,gray,very thin](-1.5,-2)grid(-0.5,-1);
\draw[step=0.5cm,gray,very thin](0,-2)grid(1,-1);
 \draw[step=0.5cm,gray,very thin](-1.5,-0.5)grid(-0.5,0.5);
\draw[step=0.5cm,gray,very thin](0,-0.5)grid(1,0.5);
\draw (-1 cm,-1.5cm) node[cross=3pt,rotate=0] {};
\draw (0.5 cm,-1.5cm) node[cross=3pt,rotate=0] {};
\draw (-1 cm,0.cm) node[cross=3pt,rotate=0] {};
\draw (0.5 cm,0cm) node[cross=3pt,rotate=0] {};
\end{tikzpicture}
\caption{FA-$1$f\label{fig:FA1f}}
\end{center}
\end{subfigure}
 \begin{subfigure}{0.3\textwidth}
\begin{center}
 \begin{tikzpicture}
 \draw[step=0.5cm,gray,very thin](-1.5,-2)grid(-0.5,-1);
\draw[step=0.5cm,gray,very thin](0,-2)grid(1,-1);
 \draw[step=0.5cm,gray,very thin](-1.5,-0.5)grid(-0.5,0.5);
\draw[step=0.5cm,gray,very thin](0,-0.5)grid(1,0.5);
 \draw[step=0.5cm,gray,very thin](-1.5,0.999)grid(-0.5,2);
\draw[step=0.5cm,gray,very thin](-0.0001,0.999)grid(1,2);
\draw (-1cm,1.5cm) node[cross=3pt,rotate=0] {};
\draw (0,1.5) circle (1 mm) ;
\draw (0.5,2) circle (1 mm) ;
\draw (-1.5,1.5) circle (1 mm) ;
\draw (-1,1) circle (1 mm) ;
\draw (-1,0.5) circle (1 mm) ;
\draw (-0.5,0) circle (1 mm) ;
\draw (1,0) circle (1 mm) ;
\draw (0.5,-0.5) circle (1 mm) ;
\draw (0.5,-1) circle (1 mm) ;
\draw (0.5,-2) circle (1 mm) ;
\draw (-0.5,-1.5) circle (1 mm) ;
\draw (-1.5,-1.5) circle (1 mm) ;
\draw (-1cm,0cm) node[cross=3pt,rotate=0] {};
\draw (0.5cm,1.5cm) node[cross=3pt,rotate=0] {};
\draw (0.5cm,0cm) node[cross=3pt,rotate=0] {};
\draw (0.5cm,-1.5cm) node[cross=3pt,rotate=0] {};
\draw (-1 cm,-1.5cm) node[cross=3pt,rotate=0] {};
 \end{tikzpicture}
\caption{FA-$2$f\label{fig:FA2f}}
\end{center}
 \end{subfigure}
 \begin{subfigure}{0.3\textwidth}
\begin{center}
 \begin{tikzpicture}
 \draw[step=0.5cm,gray,very thin](-1.5,-2)grid(-0.5,-1);
\draw[step=0.5cm,gray,very thin](0,-2)grid(1,-1);
 \draw[step=0.5cm,gray,very thin](-1.5,-0.5)grid(-0.5,0.5);
\draw (-1,0.5) circle (1 mm) ;
\draw (-1,-0.5) circle (1 mm) ;
\draw (-1.5,-1.5) circle (1 mm) ;
\draw (-1,-2) circle (1 mm) ;
\draw (0,-1.5) circle (1 mm) ;
\draw (0.5,-1) circle (1 mm) ;
\draw (-1cm,0cm) node[cross=3pt,rotate=0] {};
\draw (0.5cm,-1.5cm) node[cross=3pt,rotate=0] {};
\draw (-1 cm,-1.5cm) node[cross=3pt,rotate=0] {};
 \end{tikzpicture}
\caption{Duarte\label{fig:Duarte}}
\end{center}
 \end{subfigure}
\begin{subfigure}{0.3\textwidth}
\begin{center}
 \begin{tikzpicture}
 \draw[step=0.5cm,gray,very thin](-2,-2)grid(-1,-1);
\draw (-1,-1.5) circle (1 mm) ;
\draw (-1.5,-1) circle (1 mm) ;
\draw (-1.5,-1.5) node[cross=3pt,rotate=0] {};
 \end{tikzpicture}
\caption{North-East\label{fig:NE}}
\end{center}
 \end{subfigure}
 \begin{subfigure}{0.3\textwidth}
\begin{center}
 \begin{tikzpicture}
 \draw[step=0.5cm,gray,very thin](-1.5,-2)grid(-0.5,-1);
\draw[step=0.5cm,gray,very thin](0,-2)grid(1,-1);
 \draw[step=0.5cm,gray,very thin](-1.5,-0.5)grid(-0.5,0.5);
\draw[step=0.5cm,gray,very thin](0,-0.5)grid(1,0.5);
\draw (-1,0.5) circle (1 mm) ;
\draw (-0.5,0) circle (1 mm) ;
\draw (0,0.5) circle (1 mm) ;
\draw (1,0.5) circle (1 mm) ;
\draw (0,0) circle (1 mm) ;
\draw (-0.5,0.5) circle (1 mm) ;
\draw (0.5,0.5) circle (1 mm) ;
\draw (-0.5,-0.5) circle (1 mm) ;
\draw (0,-1) circle (1 mm) ;
\draw (0.5,-2) circle (1 mm) ;
\draw (-0.5,-1.5) circle (1 mm) ;
\draw (-1,-2) circle (1 mm) ;
\draw (-1.5,-2) circle (1 mm) ;
\draw (-0.5,-2) circle (1 mm) ;
\draw (0,-2) circle (1 mm) ;
\draw (0,-1.5) circle (1 mm) ;
\draw (-1cm,0cm) node[cross=3pt,rotate=0] {};
\draw (0.5cm,0cm) node[cross=3pt,rotate=0] {};
\draw (0.5cm,-1.5cm) node[cross=3pt,rotate=0] {};
\draw (-1 cm,-1.5cm) node[cross=3pt,rotate=0] {};
\end{tikzpicture}
\caption{Spiral\label{fig:spiral}}
\end{center}
 \end{subfigure}
\caption{The two-dimensional update families introduced in Section \ref{sec:rules}.\label{fig:rules}}
\end{figure}
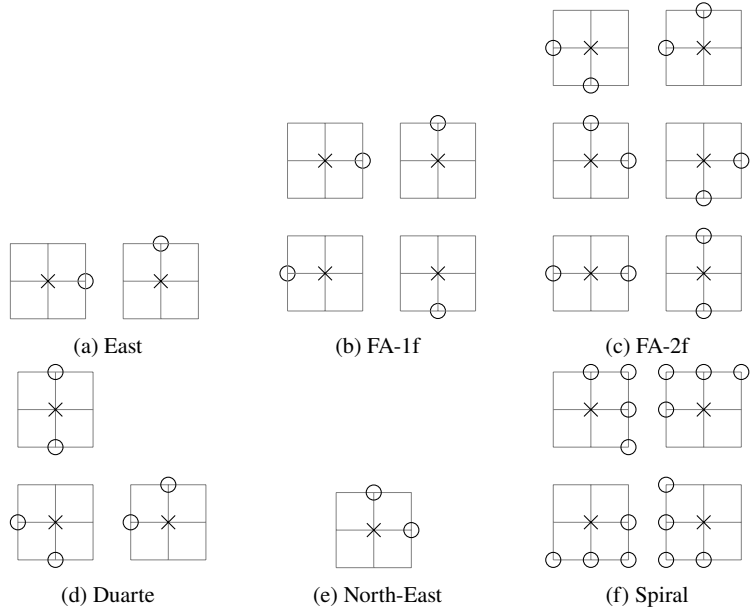
We now present some of the most commonly considered update families (see Figure~\ref{fig:rules} for their two-dimensional representation). While at this point they may seem arbitrary, in Chapter~\ref{chap:universality}, we will see that they are representatives of different universality classes displaying very different behaviour.
\begin{itemize}
\item The \textbf{East} \cites{Jackle91,Ashton06} update family is $\cU=\{\{\vec e_1\},\dots,\{\vec e_d\}\}$. That is, the East constraint at $\vec x\in\bbZ^d$ is satisfied, if $\vec x$ has an empty neighbour in some of the positive coordinate directions.
\item The \textbf{Frederickson--Andersen $j$-spin facilitated (FA-$j$f or $j$-neighbour)} \cites{Fredrickson84,Fredrickson85} update family is $\cU=\{U\subset\{\vec e_1,\dots,\vec e_d,-\vec e_1,\dots,-\vec e_d\}:|U|=j\}$, where $j\in\{1,\dots,2d\}$ is a parameter. That is, the FA-$j$f constraint at $\vec x\in\bbZ^d$ is satisfied, if $\vec x$ has at least $j$ empty neighbours. 
\item The \textbf{Duarte} \cite{Duarte89} update family is $\cU=\{\{-\vec e_1,\vec e_2\},\{-\vec e_1,-\vec e_2\},\{-\vec e_2,\vec e_2\}\}$ in two dimensions. That is, the Duarte constraint at $\vec x\in\bbZ^2$ is satisfied, if $\vec x$ has at least two empty neighbours other than $\vec x+\vec e_1$.
\item The \textbf{North-East (NE)} \cite{Reiter92} update family is $\cU=\{\{\vec e_1,\dots,\vec e_d\}\}$ for $d\ge 2$ (for $d=1$ we get the East model). That is, the constraint at $\vec x\in\bbZ^d$ is satisfied, if all the neighbours of $\vec x$ in the positive coordinate directions are empty.
\item The \textbf{Spiral} \cite{Toninelli08} update family is
$\cU=\{U_1,U_2,U_3,U_4\}$ in two dimensions, where $U_1=\{(1,-1),(1,0),(1,1),(0,1)\}$ and $U_2,U_3,U_4$ are obtained by rotating $U_1$ by $\pi/2$, $\pi$ and $3\pi/2$ around the origin respectively.
\end{itemize}

\section{The Markov process}
\label{subsec:markov}
Let us first informally describe KCM on $\mathbb Z^d$
via its so-called \emph{graphical representation}. 
Each vertex ${\vec x}\in\bbZ^d$ is equipped with a unit intensity Poisson process, whose atoms $t_{\vec x,k}$ for $k\in\bbN$ are the \emph{clock rings}. We are further given independent Bernoulli random variables $s_{\vec x,k}$ with parameter $1-q$ called \emph{coin tosses}. At the clock ring $t_{\vec x,k}$, if the current configuration $\omega$ satisfies the constraint at ${\vec x}$, we set the occupation variable $\omega_{\vec x}$ to $s_{\vec x,k}$. Such updates are called \emph{legal}. If, on the contrary, the constraint is not satisfied, the configuration remains unchanged at time $t_{\vec x,k}$.

More formally, the Markov process can be constructed via its self-adjoint Markov \emph{semigroup} $P_t:=e^{t\cL}$ on $\mathrm{L}^2(\mu)$, where the \emph{generator} $\cL$ is a non-negative self-adjoint operator with domain $\mathrm{Dom}(\mathcal L)$ that can be constructed in a  standard way (see e.g.\ \cite{Liggett05}*{Sections~I.3, IV.4}) starting from its action on \emph{local functions} (i.e.\ functions depending on the occupancy variables on  a finite number sites):
\begin{equation}
\label{eq:generator}
\cL f=\sum_{\vec x\in \bbZ^d}c_{\vec x}\cdot \left(\mu_{\vec x}(f)-f\right).
\end{equation}
Spelling out the definition of $\mu_{\vec x}$, the generator can be equivalently rewritten as
\begin{equation}
\label{eq:generatorbis}
\cL f(\omega)=\sum_{\vec x\in \bbZ^d}c_{\vec x}(\omega)\left((1-\omega_{\vec x})(1-q)+\omega_{\vec x}q\right) (f(\omega^{\vec x})-f(\omega))
\end{equation}
with $\omega^{\vec x}$  the configuration
obtained from $\omega$ by flipping its value
at $\vec x$,  i.e.
\begin{equation}
\label{eq:def:omega:flip}
\omega^{\vec x}_{\vec y}=
\begin{cases}
\omega_{\vec y} &\vec y \neq \vec x,\\
1-\omega_{\vec x} &\vec y=\vec x.
\end{cases}\end{equation}
We further introduce the Dirichlet form 
$\cD:{\mbox{Dom}}(\cL)\to\mathbb R$ defined as  
 \begin{equation}\label{eq:dirichlet}
\cD(f)=- \mu(f \cdot \cL f)=\sum_{\vec x\in \bbZ^d}\mu\left(c_{\vec x} \Var_{\vec x}(f)\right).
\end{equation}
Using the formulation \eqref{eq:generatorbis} it is not hard to verify that
\[\cD(f)=\frac{1}{2}\int\mu(\mathrm{d}\omega)\sum_{\vec x\in \bbZ^d}c_{\vec x}(\omega)\left((1-\omega_{\vec x})(1-q)+\omega_{\vec x}q\right)\left(f(\omega^{\vec x})-f(\omega)\right)^2.\]

When the initial distribution at time $t=0$ is $\nu$, the law and expectation of the KCM process on the Skorokhod space $D([0,\infty),\Omega)$ of c\`adl\`ag functions are denoted by $\bbP_{\nu}$ and $\bbE_{\nu}$ respectively (see \cite{Billingsley99}*{Chapter III} for background). If $\nu$ is concentrated over a single configuration, $\nu=\delta_{\sigma}$, we
write simply $\bbP_\sigma$ and $\bbE_\sigma$ for $\bbP_\nu$ and $\bbE_\nu$, while if $\nu=\mu$, we simply write $\bbP$ and $\bbE$. We use $\omega(t)$ to denote the state of the KCM at time $t\ge 0$.

We next discuss an important property of KCM---reversibility (see \cite{Liggett05}*{Section II.5} for background). It is not hard to verify that, since the constraint $c_{\vec x}(\omega)$ does not depend on $\omega_{\vec x}$, the dynamics satisfies detailed balance w.r.t.\ the product measure $\mu$. Therefore, $\mu$ is reversible (i.e.\ $\mu(f\cdot P_t g)=\mu(g\cdot P_t f)$ for all $f,g\in \mathrm{L}^2(\mu)$ and $t\ge0$) and therefore it is an invariant measure for the process (i.e.\ $\mu P_t=\mu$ for all $t\ge0$). However, $\mu$ is \emph{not} the unique invariant measure---e.g.\ the Dirac measure on the fully occupied configuration is clearly invariant in view of \eqref{eq:def:cx} and \eqref{eq:generator}. We nevertheless refer to the KCM with initial condition $\mu$ as the \emph{stationary process}.
  
\section{Boundary conditions}
\label{sec:finite_vol}
KCM can also be defined on finite or infinite subsets $\Lambda\subset \mathbb Z^d$ (we write $\Lambda\Subset\bbZ^d$ when we assume that $\Lambda$ is finite). In this case the most natural choice is to imagine that the configuration  is defined also outside $\Lambda$, where it is frozen and equal to some reference configuration $\sigma\in\Omega_{\bbZ^d\setminus\Lambda}$,  the \emph{boundary condition}. Then, for $\vec x\in\Lambda$, $\omega\in\Omega_{\Lambda}$, the constraint is defined as 
\begin{equation}
\label{eq:def:cx:finite}
c^{\sigma}_{\vec x}(\omega)= c_{\vec x}(\omega\cdot\sigma)
\end{equation}
(recall \eqref{eq:def:concatenation} and \eqref{eq:def:cx}).  
We denote by $\cL_{\sigma}$ and
$\mathcal{D}_{\sigma}$ the generator and  Dirichlet form of this process on $\Omega_{\Lambda}$ that are obtained by restricting the sums in \eqref{eq:generator} and \eqref{eq:dirichlet} to sites in $\Lambda$ and substituting $c_{\vec x}$ with $c^{\sigma}_{\vec x}$. We similarly denote by $\bbP_{\nu}^{\sigma}$ and $\bbE_{\nu}^\sigma$ ($\bbP_\zeta^{\sigma}$ and $\bbE_\zeta^\sigma$, if $\nu=\delta_\zeta$, and $\bbP^\sigma$ and $\bbE^\sigma$, if $\nu=\mu$) the law and expectation of the process with initial distribution $\nu$ and by $\omega^\sigma(t)$ the process at time $t$. Note that $\mu_\Lambda$ is reversible for this process.

\section{Characteristic times and critical parameters}\label{sec:time_scales}
Having defined KCM, we next discuss the kinds of questions and quantities that we seek to tackle for them. Let us start with  three natural observables. The \emph{emptying}, \emph{occupation} and \emph{persistence} times of a KCM $(\omega(t))_{t\ge 0}$ are given by
\begin{align}
\label{eq:def:tau}
\tau_0&=\inf\left\{t\ge 0:\omega_0(t)=0\right\},&\tau_1&=\inf\left\{t\ge 0:\omega_0(t)=1\right\},&\tau_\vee&=\max(\tau_0,\tau_1)
\end{align}
respectively. Turning to more analytic quantities, the \emph{relaxation time} (also known as the inverse of the \emph{spectral gap} of $\cL$) is defined as 
\begin{align}
  \label{eq:def:Trel}
\trel={}&\frac{1}{\mbox{gap}},& {\mbox{gap}}={}&\inf_{\substack{f\in \mathrm{Dom}(\cL)\\ \Var(f)\neq 0}}\frac{\cD(f)}{\Var(f)}
\end{align}
where $\cD$ is the Dirichlet form \eqref{eq:dirichlet}.
This definition is equivalent to saying that $\trel$ is the smallest constant $C\ge 0$ such that the {\emph{Poincar\'e inequality}}
\begin{equation}
  \label{eq:poincare}
  \Var(f)\leq C \cD(f)
  \end{equation}
is satisfied for any $f\in\mathrm{Dom}(\cL)$. A finite relaxation time is equivalent to the fact that the measure $\mu$ is mixing for the semigroup $P_t$ with exponentially decaying correlations (see e.g.\ \cite{Saloff-Coste97}), namely 
for all $f\in \mathrm{L}^2(\mu)$ it holds that\footnote{Indeed, by reversibility, $
    \frac{\mathrm d}{\mathrm{d}t}\var (P_t f)=-2\mathcal{D}(P_t f)$.}
\begin{equation}\label{eq:expo}\Var\left(P_t f\right)\leq e^{-2 t/\trel}\Var(f).\end{equation}
 Thus, the relaxation time controls the decay of correlations in the stationary process.
 Of course, the above time scales have no reason to be finite, so it is natural to consider the corresponding critical parameters. We define the \emph{ergodicity} and \emph{exponential decay} critical parameters
 \begin{align}
\label{eq:def:qc}\qc&=\qc(\cU)=\inf\left\{q>0:\bbP(\tau_0<\infty)=1\right\},\\
\label{eq:def:qct}\qct&=\qct(\cU)=\inf\left\{q>0:\trel<\infty\right\}.
\end{align}

Our main goals are to determine (or characterise) \qc and \qct and to study the asymptotics of $\tau_0$ and $\trel$ for the stationary process as $q\to\qc+$, as well as the behaviour of the KCM out of equilibrium (with initial condition different from $\mu$) for any $q>0$. 

Once a Poincar\'e inequality \eqref{eq:poincare} is established 
it is natural to ask whether it would be possible to prove also a stronger coercive inequality for the generator. In particular one can investigate whether a logarithmic or modified logarithmic Sobolev inequality holds (see \cite{Saloff-Coste97} for background). These correspond, respectively, to the existence of a finite constant  $C_{\mathrm{LS}}$ and $C_{\mathrm{MLS}}$ such that for any nonnegative $f\in\mathrm{Dom}(\cL)$ we have
\begin{align}
\label{eq:logsob}
\Ent(f)&{}\leq C_{\mathrm{LS}}\mathcal{E}\left(\sqrt{f},\sqrt{f}\right),\\
\label{eq:modlogsob}\Ent(f) &{}\leq C_\mathrm{MLS}\mathcal{E} (f,\log f),\end{align}
where for any two functions $f,g$ we set 
\begin{align*}
\Ent(f)&{}=\mu\left(f\log(f/\mu(f))\right),&\cE(f,g)&{}=-\mu(f\mathcal L g).\end{align*}
The (stronger) inequality \eqref{eq:logsob} is known to be equivalent to the hypercontractivity property of the semigroup $P_t$ \cite{Diaconis96}. Instead, \eqref{eq:modlogsob}
is equivalent to exponential decay of the relative entropy for the generator $P_t$ (see \cite{Bobkov06}), namely for each probability measure $\nu$ on $\Omega$  it holds that\footnote{Indeed, if we let $f_t$ be the relative density of $\nu P_t$ w.r.t.\ $\mu$ it holds that $\frac{\mathrm d}{\mathrm{d}t} H(\nu P_t ||\mu)= \frac{\mathrm d}{\mathrm{d}t} \Ent(f_t)=-\cE(f_t,\log f_t)$.}
\begin{equation}\label{eq:rel_ent_dec_1}H(\nu P_t ||\mu) \leq e^{-t/C_{\mathrm{MLS}}}H(\nu ||\mu) \end{equation}
where for any two measures $\nu_1,\nu_2$ on $\Omega$ we denote by $H(\nu_1||\nu_2)$ the relative entropy (or Kullback--Leibler divergence) of $\nu_1$ w.r.t.\ $\nu_2$. 

The definitions in \eqref{eq:def:tau}, \eqref{eq:def:Trel},  \eqref{eq:logsob}, \eqref{eq:modlogsob} naturally extend to the finite volume setting. We denote the relaxation time (logarithmic or modified logarithmic Sobolev constants) for the KCM  with boundary condition $\sigma$ by $T^{\sigma}_{\rm{rel}}$ (resp.\ $C_\mathrm{MLS}^{\sigma}$ and $C_\mathrm{MLS}^{\sigma}$). When $\sigma=\bzero_{\bbZ^d\setminus\Lambda}$ for $\Lambda\subset\bbZ^d$, we simplify the notation to $T^{\Lambda}_{\rm{rel}}$ (resp.\ $C_\mathrm{MLS}^{\Lambda}$ and $C_\mathrm{MLS}^{\Lambda}$).

We finally define another natural time scale known as mixing time of the KCM on $\Lambda\subset\bbZ^d$ with boundary condition $\sigma\in\Omega_{\bbZ^d\setminus\Lambda}$ (see \cite{Levin09} for background). The total variation distance $d_{\rm{TV}}$ between two measures $\mu_1,\mu_2$ on $\Omega_\Lambda$ is given by
\begin{equation}
\label{eq:def:dTV}
d_{\rm{TV}} (\mu_1,\mu_2)=\sup_{A\in\mathcal F} \left|\mu_1(A)-\mu_2(A)\right|,
\end{equation}
where $\mathcal F$ denotes the Borel $\sigma$-field generated by the open sets of $\Omega_\Lambda$. For $\varepsilon\in(0,1)$, the mixing time is
\begin{equation}\label{eq:def:tmix}t_{\rm{mix}}^{\sigma}(\varepsilon)=\inf\left\{t>0, \ \max_{\omega\in
\Omega_{\Lambda}}d_{\rm{TV}} \left(\delta_\omega P_t,\mu\right)\le \varepsilon\right\}.
\end{equation}
If $\sigma=\bzero_{\bbZ^d\setminus\Lambda}$, we simplify the notation to $t_{\rm{mix}}^{\sigma}=t_{\rm{mix}}^{\Lambda}$. The mixing time has a natural probabilistic interpretation: it is the time when it is impossible to distinguish (e.g.\ by a statistical test) the law of the process at that time from the equilibrium measure, regardless of the initial state.

	\chapter{From bootstrap percolation to kinetically constrained models}
\label{chap:BP}
\abstract{
In this chapter we introduce bootstrap percolation cellular automata, which are instrumental for the study of KCM. We first state some relevant known results for these automata. We then show that several fundamental properties of KCM, including ergodicity, mixing and exponential relaxation, can be directly related to their bootstrap percolation counterparts.}

\section{Bootstrap percolation}
Bootstrap percolation (BP) is a family of monotone cellular automata. They may be viewed as the discrete time synchronous monotone analogue of KCM. Specific instances of BP have been studied since the 1970s \cites{Chalupa79,Kogut81,Pollak75}, but it is convenient to directly introduce them in greater generality as considered in \cites{Schonmann90,Gravner99,Bollobas15}. Like KCM, BP is defined by an update family $\cU$ (recall Section~\ref{sec:rules}). Given $\omega\in\Omega$, we define $\sB_\cU(\omega)\in\Omega$ by\footnote{Note that, in most of BP literature, the roles of the 0 and 1 states are reversed.}
\begin{equation}
\label{eq:def:BP}
(\sB_\cU(\omega))_{\vec x}=\begin{cases}
    0&\omega_{\vec x}=0\text{ or }\exists U\in\cU,\forall \vec u\in U,\omega_{\vec x+\vec u}=0,\\
    1&\text{otherwise}
\end{cases}
\end{equation}
for all $\vec x\in\bbZ^d$. In words, in one discrete time step we empty all sites for which the constraint (recall \eqref{eq:def:cx}) is satisfied. In BP, empty sites remain empty. In view of this monotonicity, it is natural to define the \emph{closure}
\begin{equation}
    \label{eq:def:closure}
    [\omega]=[\omega]_\cU=\inf_{t\in\bbN}\sB_\cU^{\circ t}(\omega)=\lim_{t\to\infty}\sB_\cU^{\circ t}(\omega)\in\Omega,
\end{equation}
where $\sB_\cU^{\circ t}$ denotes the $t$-fold iteration of $\sB_\cU$ and the infimum and limit are taken with respect to the product partial order and product topology respectively. That is, $[\omega]$ is the configuration obtained upon iterating the bootstrap percolation map of \eqref{eq:def:BP}. Similarly to \eqref{eq:def:tau}, we define the \emph{BP emptying time}
\begin{equation}\tbp=\inf\left\{t\in\bbN:\left(\sB_\cU^{\circ t}(\omega)\right)_0=0\right\}\in\bbN\cup\{\infty\}.\end{equation}
On a domain $\Lambda\subset\bbZ^d$ with boundary condition $\sigma\in\Omega_{\bbZ^d\setminus\Lambda}$, we set
\begin{equation}
\label{eq:def:BP:finite}
\left(\sB_\cU^\sigma(\omega)\right)_{\vec x}=\begin{cases}
    0&\omega_{\vec x}=0\text{ or }\exists U\in\cU,\forall \vec u\in U,(\sigma\cdot\omega)_{\vec x+\vec u}=0,\\
    1&\text{otherwise}
\end{cases}
\end{equation}
for $\omega\in\Omega_\Lambda$ and $\vec x\in\Lambda$. We further define $[\omega]^{\sigma}=\lim_{t\to\infty}(\sB_\cU^\sigma)^{\circ t}(\omega)$ pointwise.

So far BP is completely deterministic. We next introduce randomness by considering an initial condition $\omega$ distributed according to the product Bernoulli measure $\mu_q$ with parameter (density of initially occupied sites) $1-q\in[0,1]$ (recall Section \ref{Section:Setting}). Following \eqref{eq:def:qc} and \eqref{eq:def:qct}, we define the \emph{emptying} and \emph{exponential decay} critical thresholds
\begin{align}
    \label{eq:def:qcbp}q_{\mathrm c}^{\mathrm{BP}}&=\inf\left\{q>0:\mu_q\left(\tbp<\infty\right)=1\right\},\\
    \label{eq:def:qctbp}\qct^{\mathrm {BP}}&=\inf\left\{q>0:\liminf_{t\to\infty}\frac{-\log\mu_q(\tbp>t)}{t}>0\right\}.
\end{align}
Note that by ergodicity of the product measure (with respect to translations, see e.g.\ \cite{Keller98} for background), we have $\mu_q(\tbp<\infty)=1$ if and only if $\mu_q([\omega]=\bzero_{\bbZ^d})=1$. As for KCM, one is primarily interested in determining these thresholds and the asymptotics of $\tbp$ as $q\to q_{\mathrm c}^{\mathrm{BP}}+$.

In the present text our focus is on KCM, so we take BP results for granted. Let us therefore gather a few facts about the BP models corresponding to the update families introduced in Section~\ref{sec:rules}.
\begin{theorem}[BP background]\label{th:BP}\leavevmode
\begin{itemize}
    \item For East BP in dimension $d\ge 1$, $q_{\mathrm c}^{\mathrm{BP}}=\qct^{\mathrm{BP}}=0$ and $q^{1/d}\tbp$ converges to a Weibull distribution, as $q\to0$: $\mu(q^{1/d}\tbp\ge t)\to e^{-t^d/d!}$ for $t\ge 0$.
    \item For 1-neighbour BP in dimension $d\ge 1$, $q_{\mathrm c}^{\mathrm{BP}}=\qct^{\mathrm{BP}}=0$ and $q^{1/d}\tbp$ converges to a Weibull distribution, as $q\to0$: $\mu(q^{1/d}\tbp\ge t)\to e^{-(2t)^d/d!}$ for $t\ge 0$.
    \item For $j$-neighbour BP in dimension $d$ with $d\ge j\ge 2$ we have $q_{\mathrm c}^{\mathrm{BP}}=\qct^{\mathrm{BP}}=0$ and $q^{1/(d-j+1)}\log^{\circ (j-1)}\tbp$ converges in probability to a constant\footnote{See \cite{Balogh09a}*{(1)-(3)} for an explicit expression of $\lambda(d,j)$.} $\lambda(d,j)>0$, as $q\to0$. 
    \item For $2$-neighbour BP in $d=2$ dimensions, there exist positive constants\footnote{We have $\lambda=\pi^2/18$ and $\lambda_2\approx7.0545$, see \cite{Hartarsky24locality}*{Section~A.1.2} for an explicit expression.} $\lambda,\lambda_2$ such that, as $q\to0$, we have
    \begin{equation}
        \label{eq:2nBP}\mu_q\left(\left|\log \tbp-\frac{\lambda}{q}+\frac{\lambda_2}{\sqrt q}\right|\le \frac{\log^2(1/q)}{\sqrt[3]q}   \right)\to1.
    \end{equation}
    \item For Duarte BP we have $q_{\mathrm c}^{\mathrm{BP}}=\qct^{\mathrm{BP}}=0$ and $q\log \tbp/\log^2(1/q)$ converges in probability to a positive constant, as $q\to0$.
    \item For North-East BP in $d\ge 2$ dimensions we have $q_{\mathrm c}^{\mathrm{BP}}=\qct^{\mathrm{BP}}=1-p_{\mathrm c}^{\mathrm{OP},d}$, where $p_{\mathrm c}^{\mathrm{OP},d}\in(0,1)$ is the critical probability of oriented site percolation in $d$ dimensions (see e.g.\ \cites{Durrett84,Liggett05,Liggett99,Hartarsky22GOSP} for background) and $\mu_{q_{\mathrm c}^{\mathrm{BP}}}(\tbp=\infty)=0$.
    \item For Spiral BP we have $q_{\mathrm c}^{\mathrm{BP}}=\qct^{\mathrm{BP}}=1-p_{\mathrm c}^{\mathrm{OP},2}$ and $\mu_{q_{\mathrm c}^{\mathrm{BP}}}(\tbp=\infty)>0$.
    \item For $j$-neighbour BP in dimension $d$ with $j>d\ge 1$ we have $q_{\mathrm c}^{\mathrm{BP}}=\qct^{\mathrm{BP}}=1$.
\end{itemize}
\end{theorem}
\begin{proof}
For East BP, we only sketch the argument and leave the details as an exercise to the reader. First, we verify that for any $t\in\bbN$, we have $\mu_q(\tbp>t)=(1-q)^{N_d(t)}$, where $N_d(t)=\{(x_1,\dots,x_d)\in\bbN^d:\sum_{i=1}^dx_i\le t\}$. One can check that $N_d(t)=\binom{t+d}{d}=t^d/d!+O_d(t^{d-1})$, where the implicit constant may depend on the dimension $d$. In particular, for any $q>0$, $\mu_q(\tbp=\infty)=0$ and $\mu_q(\tbp>t)$ decays at least exponentially, so $q_{\mathrm{c}}^{\mathrm{BP}}=\qct^{\mathrm{BP}}=0$. Moreover, the convergence in distribution follows from the asymptotics of $N_d(t)$.

The proof for 1-neighbour BP is analogous, replacing $N_d(t)$ by the volume of the discrete $\ell^1$ ball of radius $t$, whose cardinal is asymptotically equivalent to $2^dN_d(t)$.

For $j$-neighbour BP with general $2\le j\le d$, the asymptotics of $\log^{\circ(j-1)}\tbp$ is due to Balogh, Bollob\'as, Duminil-Copin and Morris \cite{Balogh12} (the case $j=2$ was established by Holroyd \cites{Holroyd03,Holroyd06}), while the identification of the critical value is due to Schonmann \cite{Schonmann92}*{Theorem~3.1}.

The result for 2-neighbour BP is due to Hartarsky and Teixeira \cite{Hartarsky24locality}.

The quantitative result for Duarte BP is due to Bollob\'as, Duminil-Copin, Morris and Smith \cite{Bollobas17}, while the qualitative one is due to Schonmann \cites{Schonmann92,Schonmann90}.

The result for North-East BP follows from the fact that $\tbp=\infty$ if and only if there is an infinite oriented path of occupied sites from the origin (see \cite{Schonmann90}), together with classical results on oriented site percolation \cites{Aizenman87,Menshikov86,Bezuidenhout90}.

The result for Spiral BP is due to Toninelli and Biroli \cite{Toninelli08} (also see \cites{Hartarsky21,Toninelli07,Toninelli07a,Toninelli06,Jeng08}).

For $j$-neighbour BP with $j>d\ge 1$ it suffices to observe that for $q<1$, we have 
\[\mu_q\left(\tbp=\infty\right)\ge\mu_q\left(A_0\cap\{0,1\}^d=\varnothing\right)=(1-q)^{2^d}>0.\]
This concludes the proof for all models.
\end{proof}
The reader may have noticed the following pattern \cite{Hartarsky21}*{Conjecture~8.1} (also see \cite{Schonmann92}) in Theorem~\ref{th:BP}.
\begin{conjecture}[Sharp phase transition]
\label{conj:qc:qct}
For any update family $\cU$ in any dimension it holds that $q_{\mathrm{c}}^{\mathrm{BP}}=\qct^{\mathrm{BP}}$.
\end{conjecture}
This is an important open problem in BP theory, which has so far been resolved for update families contained in an open half-space with the origin on its boundary \cite{Hartarsky22sharpness}*{Theorem~1.6}, as well as those with $q_{\mathrm{c}}^{\mathrm{BP}}=0$. For the latter assertion, note that \cite{BalisterNaNb} proves a stretched exponential bound on the tail of $\tbp$, but a standard renormalisation argument \cite{Schonmann92} can be used to recover an exponential decay.

Concerning the continuity of the phase transition, there is not even a guess what the answer should be, leaving the following problem wide open.
\begin{problem}
    Determine which update families satisfy $\mu_{q_{\mathrm c}^{\mathrm{BP}}}(\tbp=\infty)=0$ like North-East BP and unlike Spiral BP.
\end{problem}

\section{Legal paths}
\label{sec:legal:paths}
We next introduce  the notion of legal path that will be instrumental in several proofs. In words, a legal path is a sequence of configurations differing by a legal update (recall Section~\ref{subsec:markov}). 
\begin{definition}[Legal paths]
\label{def:legal:path}
Given a domain $\Lambda\subset\bbZ^d$, boundary condition $\sigma\in\Omega_{\bbZ^d\setminus\Lambda}$ and two configurations $\omega,\omega'\in\Omega_\Lambda$, a \emph{legal path from $\omega$ to $\omega'$ in $\Lambda$} is a finite sequence $(\omega^{(i)})_{i=0}^n$ of configurations in $\Omega_\Lambda$ such that $\omega^{(0)}=\omega$, $\omega^{(n)}=\omega'$ and for each $i\in\{1,\dots,n\}$ it holds that either $\omega^{(i)}=\omega^{(i-1)}$ or there exists $\vec x^{(i)}\in\Lambda$ such that $\omega^{(i)}=(\omega^{(i-1)})^{\vec x^{(i)}}$ (recall \eqref{eq:def:omega:flip}) and $c^\sigma_{\vec x^{(i)}}(\omega^{(i)})=1$ (recall \eqref{eq:def:cx:finite}). The \emph{length} of the legal path $(\omega^{(i)})_{i=0}^n$ is $n$. Notice that if $(\omega^{(i)})_{i=0}^n$ is a legal path from $\omega$ to $\omega'$ in $\Lambda$, then its \emph{inverse} $(\omega^{(n-i)})_{i=0}^n$ is a legal path from $\omega'$ to $\omega$ in $\Lambda$.
\end{definition}
It turns out that BP provides a simple way to know when a legal path exists and to construct it. We start by observing that the BP closure is invariant along legal paths.
\begin{lemma}[Invariance of closure]
\label{lem:closure}
Let $\omega\in\Omega$ and $\vec x\in\bbZ^d$ be such that $c_{\vec x}(\omega)=1$. Then $[\omega]=[\omega^{\vec x}]$ and for any $\Lambda\subset\bbZ^d$ such that $\vec x\in\Lambda$ we have $[\omega_\Lambda]^{\omega_{\bbZ^d\setminus\Lambda}}=[\omega^{\vec x}_\Lambda]^{\omega_{\bbZ^d\setminus\Lambda}}$. Consequently, for $\omega',\omega''\in\Omega_\Lambda$ connected by a legal path in $\Lambda$ with boundary condition $\sigma\in\Omega_{\bbZ^d\setminus\Lambda}$ we have $[\omega'\cdot\sigma]=[\omega''\cdot\sigma]$ and $[\omega']^\sigma=[\omega'']^{\sigma}$.
\end{lemma}
\begin{proof}
Assume without loss of generality that $\omega_{\vec x}=0$. By \eqref{eq:def:cx} and \eqref{eq:def:BP}, we have that $(\sB_\cU(\omega^{\vec x}))_{\vec x}=0=\omega_{\vec x}$, so $\sB_\cU(\omega^{\vec x})\le \omega\le \omega^{\vec x}$. By \eqref{eq:def:closure}, since $[\cdot]$ is non-increasing, this yields $[\omega]=[\omega^{\vec x}]$. The proof of $[\omega_\Lambda]^{\omega_{\bbZ^d\setminus\Lambda}}=[\omega^{\vec x}_\Lambda]^{\omega_{\bbZ^d\setminus\Lambda}}$ is analogous. The statement on legal paths follows by induction on the length.
\end{proof}
\begin{lemma}
\label{lem:legal:path}
Let $\Lambda\subset\bbZ^d$, $\sigma\in\Omega_{\bbZ^d\setminus\Lambda}$ and $\omega\in\Omega_\Lambda$. Let $\vec x\in\Lambda$ be such that $\omega_{\vec x}=1$. There exists a legal path from $\omega$ to $\omega^{\vec x}$ if and only if $[\omega]^\sigma_{\vec x}=0$. Moreover, if it exists, the legal path can be chosen with length at most $2|\Lambda|$.
\end{lemma}
\begin{proof}
The only if direction was proved in Lemma~\ref{lem:closure}. Assume that $[\omega]_{\vec x}^\sigma=0$. Then there exists a finite set $\Lambda'\subset\Lambda$ such that $[\omega_{\Lambda'}]^{\sigma\cdot\omega_{\Lambda\setminus\Lambda'}}_{\vec x}=0$. Up to replacing $\Lambda$ by $\Lambda'$, we may assume that $\Lambda$ is finite. We construct a legal path in $\Lambda$ of length at most $|\Lambda|$ from $\omega$ to $[\omega]^\sigma$. To achieve this, we empty an arbitrary occupied site whose constraint with boundary condition $\sigma$ is satisfied. Such a vertex always exist until we reach $[\omega]^\sigma$. We similarly obtain a legal path from $\omega^{\vec x}$ to $[\omega^{\vec x}]^\sigma=[\omega]^\sigma$ and conclude by concatenating its inverse with the legal path from $\omega$.
\end{proof}
Deducing or proving the following corollary is left as an exercise to the reader.
\begin{corollary}[Legal paths and closure]
\label{cor:legal:path}
Let $\Lambda\subset\bbZ^d$, $\sigma\in\Omega_{\bbZ^d\setminus\Lambda}$ and $\omega,\omega'\in\Omega_\Lambda$ be such that 
$\sum_{\vec x\in \Lambda}|\omega_{\vec x}-\omega'_{\vec x}|<\infty$. Then there exists a legal path from $\omega$ to $\omega'$ in $\Lambda$ with boundary condition $\sigma$, if and only if $[\omega]^\sigma=[\omega']^{\sigma}$ (recall \eqref{eq:def:BP:finite}).
\end{corollary}

\begin{definition}[Ergodic boundary condition]\label{def:ergoBC}
Let $\Lambda\Subset\bbZ^d$. We say that a boundary condition $\sigma\in\Omega_{\mathbb Z^d\setminus \Lambda}$ is \emph{ergodic}, if  
$[\bone_\Lambda]^\sigma=\bzero_\Lambda$. By Corollary~\ref{cor:legal:path}, this is equivalent to $\mathcal L^{\sigma}$ defining an ergodic process on $\Omega_{\Lambda}$.\end{definition}
\section{Ergodicity}
\label{sec:ergo}
The simple deterministic statements of Section~\ref{sec:legal:paths} entail the following result of Cancrini, Martinelli, Roberto and Toninelli \cite{Cancrini08}*{Proposition~2.4}, whose fundamental importance for KCM is apparent in view of Theorem~\ref{th:BP}. 
\begin{theorem}[Ergodicity]
\label{th:ergo}
For any update family $\cU$ and $q\in(0,1)$, the following are equivalent
\begin{enumerate}
    \item \label{item:ergo:5}$\bbP_{\mu_q}(\tau_\vee<\infty)=1$,
    \item \label{item:ergo:4}$\bbP_{\mu_q}(\tau_0<\infty)=1$,
    \item \label{item:ergo:1}$\mu_q(\tbp<\infty)=1$,
    \item \label{item:ergo:2}$0$ is a simple eigenvalue of $\cL$ on $\mathrm{L}^2(\mu_q)$, that is, the dynamics is \emph{ergodic},
    \item \label{item:ergo:3}for all $f\in\mathrm{L}^2(\mu_q)$ we have $\lim_{t\to\infty} P_t f=\mu_q(f)$, that is, the dynamics is \emph{mixing}.
\end{enumerate}
In particular, $\qc=q_{\mathrm c}^{\mathrm{BP}}$ (recall \eqref{eq:def:qc} and \eqref{eq:def:qcbp}). 
\end{theorem}
\begin{proof}
\textbf{\ref{item:ergo:5} implies \ref{item:ergo:4}.} This follows directly from the definition \eqref{eq:def:tau}.\\
\textbf{\ref{item:ergo:4} implies \ref{item:ergo:1}.} Consider the function $f:\Omega\to\{0,1\}:\omega\mapsto[\omega]_0$ in $\mathrm{L}^2(\mu_q)$. By Lemma~\ref{lem:closure}, we have that $\cL_{n}f=0$, where $\cL_n$ is the KCM generator on $\Omega$ defined by restricting the sum in \eqref{eq:generator} to $\Lambda_n=\{-n,\dots,n\}^d$. Moreover, clearly $\cL_{n}g\to\cL g$ for any local function $g$, so $\cL_{n}f\to\cL f$ (see \cite{Liggett05}*{Corollary I.3.14}\footnote{As explained in \cite{Liggett05}*{Section IV.4}, the corollary and other results apply in $\mathrm L^2(\mu_q)$ instead of the space of continuous functions for the product topology.}). Thus, $\cL f=0$, so that for any $t\ge 0$,  $\bbP_{\mu_q}$-a.s.\ $f(\omega(t))=f(\omega(0))$. Consequently, \begin{equation}
\label{eq:ergo:Q}\bbP_{\mu_q}\left(\forall t\in\bbQ\cap[0,\infty), f(\omega(t))=f(\omega(0))\right)=1.\end{equation} 
Assume that $\bbP_{\mu_q}(\tau_0<\infty)=1$. Then $\omega_0(t)=0$ for all $t\in[\tau_0,\tau_0+\varepsilon]$ for some random $\varepsilon>0$. Thus, \eqref{eq:ergo:Q} gives $\bbP_{\mu_q}(f(\omega(0))=f(\omega(\tau_0))=0)=1$, which concludes the proof, since $\mu_q(\tbp<\infty)=1-\mu_q(f)=1$.\\
\textbf{\ref{item:ergo:1} implies \ref{item:ergo:2}.} Fix $f\in\mathrm{L}^2(\mu_q)$ such that $\cL f=0$. Then, recalling \eqref{eq:dirichlet}, we have $\sum_{\vec x\in\bbZ^d}\mu_q(c_{\vec x}\var_{\vec x}(f))=\cD(f)=0$, so each of the (non-negative) summands is 0. We seek to prove that $\var(f)=0$ and our starting point is the unconstrained Poincar\'e inequality\footnote{In other words, the spectral gap of the generator of the product KCM corresponding to $\cU=\{\varnothing\}$ is $1$, since it is the tensor product of irreducible 2-state Markov processes with total rate $1$. This classical fact also follows e.g.\ by taking $\cX=\bbX_1$ in Lemma~\ref{lem:two-block}.}
\begin{equation}
\label{eq:Poincare:unconstrained}\var(f)\le \sum_{\vec x\in\bbZ^d}\mu(\var_{\vec x}(f)).\end{equation}
Assume further that \ref{item:ergo:1} holds and fix some $\vec x\in\bbZ^d$. Then
\[\mu\left(\var_{\vec x}(f)\right)\le \mu\left(\omega_0\left(f\left(\omega^{\vec x}\right)-f(\omega)\right)^2\right)=\sum_{n=1}^\infty \mu\left(\1_{\cE_n\setminus\cE_{n-1}}(\omega)\left(f\left(\omega^{\vec x}\right)-f(\omega)\right)^2\right)\]
where $\cE_n=\{\omega\in\Omega:[\omega_{\Lambda_n}]^{\bone_{\bbZ^d\setminus\Lambda_n}}_{{\bx}}=0\}$ and $\Lambda_n=\vec x+\{-n,\dots,n\}^d$. But by Lemma \ref{lem:legal:path} for any $\omega\in\cE_n\setminus\cE_{0}$ there exists a legal path $(\omega^{(i)})_{i=0}^{2|\Lambda_{n}|}$ from $\omega$ to $\omega^{\vec x}$ in $\Lambda_{n}$. Writing $f(\omega^{\vec x})-f(\omega)$ telescopically along this legal path and using the Cauchy--Schwarz inequality yields
\[\1_{\cE_n\setminus\cE_{n-1}}(\omega)\left(f\left(\omega^{\vec x}\right)-f(\omega)\right)^2\le 2|\Lambda_n|\sum_{i=1}^{2|\Lambda_n|}c_{\vec x^{(i)}}\left(\omega^{(i)}\right)\left(f\left(\omega^{(i)}\right)-f\left(\omega^{(i-1)}\right)\right)^2,\]
where $\vec x^{(i)}$ is the location of the legal update from $\omega^{(i-1)}$ to $\omega^{(i)}$.
In order to recover $\mu(c_{\vec x^{(i)}}\var_{\vec x^{(i)}}(f))$ from the $\mu$-average of the last summand, we only need to perform a change of measure and observe that $\mu(\omega)/\mu(\omega^{(i)})\le (q(1-q))^{-|\Lambda_n|}$. Putting everything together, we obtain
\[\mu(\var_{\vec x}(f))\le \sum_{n=1}^\infty 4|\Lambda_n|^2(q(1-q))^{-|\Lambda_n|-1}\sum_{\vec y\in\Lambda_n}\mu\left(c_{\vec y}\var_{\vec y}(f)\right)=0.\]
Recalling \eqref{eq:Poincare:unconstrained}, this gives that $f$ is $\mu$-a.s.\ constant as desired.\\
\textbf{\ref{item:ergo:2} implies \ref{item:ergo:5}.} This follows from the ergodic theorem (see e.g.\ \cite{Billingsley65}) applied to the function $\omega\mapsto\omega_{0}$.\\
\textbf{\ref{item:ergo:2} is equivalent to \ref{item:ergo:3}.} This is \cite{Liggett05}*{Theorem~IV.4.13} (based on operator spectral theory).
\end{proof}

Before moving on, let us comment on the proof, which showcases two important ideas. Firstly, in order to empty the origin, we need to be able to do so in BP. Secondly, if BP is able to empty some site, we can turn that into a Poincar\'e inequality for the corresponding KCM. Building on these two insights, one can go surprisingly far. The first one is useful for obtaining lower bounds on $\tau_0$, while the second one enables upper bounds. Finally, let us mention that the proof of the implication from \ref{item:ergo:1} to \ref{item:ergo:2} is our first encounter with the canonical path technique for Markov chains originating in \cites{Sinclair89,Lawler88} (also see e.g.\ \cite{Levin09}*{Section 13.5}).

\section{Exponential decay}
\label{sec:exponential:decay}
In Section~\ref{sec:ergo} we saw that ergodicity and mixing of KCM are equivalent to BP a.s.\ emptying $\bbZ^d$. All these results are purely qualitative and provide no quantitative control whatsoever, e.g.\ on the tails of emptying times. Our next task is to transfer the tail behaviour of BP to KCM, by proving the following result adapted from \cites{Hartarsky21,Cancrini08,Martinelli19,Cancrini09}. Once again, its interest is made clear by Theorem~\ref{th:BP} (also recall Conjecture~\ref{conj:qc:qct} and Theorem~\ref{th:ergo}).
\begin{theorem}[Exponential decay]
\label{th:exponential}
For any update family $\cU$ and $q\in(0,1)$, the following are equivalent.
\begin{enumerate}
    \item \label{item:expo:2}$\liminf_{t\to\infty}-\log\bbP_{\mu_q}(\tau_\vee>t)/t>0$,
    \item \label{item:expo:4}$\liminf_{t\to\infty}-\log\bbP_{\mu_q}(\tau_0>t)/t>0$,
    \item \label{item:expo:6}$\liminf_{t\to\infty}-\log\mu_q(\tbp>t)/t>0$,
    \item \label{item:expo:7}$\trel<\infty$ (recall \eqref{eq:def:Trel}).
\end{enumerate}
In particular, $\qct=\qct^{\mathrm{BP}}$ (recall \eqref{eq:def:qct} and \eqref{eq:def:qctbp}).
\end{theorem}
\begin{proof}[Sketch]
\textbf{\ref{item:expo:2} implies \ref{item:expo:4}.} This follows directly from the definition \eqref{eq:def:tau}.\\
\textbf{\ref{item:expo:4} implies \ref{item:expo:6}.} We claim that there exists a constant $\delta>0$ (depending on $\cU$ and $q$) such that if we run BP and KCM from the same initial configuration $\omega\in\Omega$, 
\begin{equation}
\label{eq:t0t0bp}\bbP_\omega\left(\tau_0\le \delta \tbp\right)\le e^{-\tbp(\omega)}.\end{equation}
The idea is that, in order to empty the origin, legal updates have to occur in the right order along some path of length $\tbp(\omega)$, starting at the origin. Consecutive vertices of the path are allowed to be at distance at most $\max_{U\in\cU,\vec u\in U} \|\vec u\|$. However, the number of such paths is at most $e^{C\tbp(\omega)}$  for some $C>0$. Moreover, the probability that the sum of $N$ exponential random variables of mean 1 is smaller than $\delta N$ is at most $e^{-2CN}$, choosing $\delta$ sufficiently small depending on $C$. We conclude the proof of the claim by a union bound on the possible paths. See \cite{Martinelli19}*{Lemma~4.3} for more details.\\
\textbf{\ref{item:expo:6} implies \ref{item:expo:7}.} For simplicity of the presentation, we focus on the two-dimensional case and assume that $\{(2,1),(1,2)\}\in\cU$. The proof proceeds in three steps. First, we prove that for $q=1-\varepsilon_0$ sufficiently close to 1, $\trel<\infty$ for the update family $\cU'=\{U'\}=\{\{(1,0),(1,1),(0,1)\}\}$. Second, we perform a renormalisation,\footnote{In physics, one would rather speak of coarse-graining, but we adopt the mathematical jargon.} by tessellating space into large square boxes, which are deemed good if the $\cU$-KCM restricted to the box is able to empty most of the bottom and left boundaries of the box in the current configuration. Third, we show how to completely empty a given (possibly non-good) box, assuming that the three neighbouring boxes corresponding to $\cU'$ are good. We postpone the discussion of the first step to Section~\ref{sec:bisection:2d}, where the bisection technique is presented (for the full details, refer to \cite{Cancrini08}*{Section~4}).

\begin{figure}
\begin{minipage}{\textwidth}
\leftfigure[c]{\begin{tikzpicture}[line cap=round,line join=round,>=triangle 45,x=4.0cm,y=4.0cm]
\fill[fill=black,fill opacity=0.25] (0,0) -- (0.7,0) -- (0.7,0.7) -- (0,0.7) -- cycle;
\draw (0,0)-- (1,0);
\draw (1,0)-- (1,1);
\draw (1,1)-- (0,1);
\draw (0,1)-- (0,0);
\draw [very thick] (0.1,0.1)-- (0.9,0.1);
\draw [very thick] (0.9,0.1)-- (0.9,0.9);
\draw [very thick] (0.9,0.9)-- (0.1,0.9);
\draw [very thick] (0.1,0.9)-- (0.1,0.1);
\draw (0,0)-- (0.7,0) node[below]{$(1-3\varepsilon)n$};;
\draw (0.7,0)-- (0.7,0.7);
\draw (0.7,0.7)-- (0,0.7)
;
\draw (0,0.7)-- (0,0);
\draw [dashed] (0.9,0.1)-- (0.7,0);
\draw [dashed] (0.1,0.9)-- (0,0.7);
\draw (0.1,0) node[below]{$\varepsilon n$};
\draw (1,1) node[below right]{$\Lambda$};
\end{tikzpicture}}\hspace{\fill}
\rightfigure[c]{\begin{tikzpicture}[line cap=round,line join=round,>=triangle 45,x=2cm,y=2cm]
\fill[fill=black,fill opacity=0.25] (1,0) -- (1.7,0) -- (1.7,0.7) -- (1,0.7) -- cycle;
\fill[fill=black,fill opacity=0.25] (1,1) -- (1.7,1) -- (1.7,1.7) -- (1,1.7) -- cycle;
\fill[fill=black,fill opacity=0.25] (0,1) -- (0.7,1) -- (0.7,1.7) -- (0,1.7) -- cycle;
\draw [very thick] (0,0)-- (1,0)-- (1,1)-- (0,1)--cycle;
\draw (1,0)-- (1.7,0);
\draw (1.7,0)-- (1.7,0.7);
\draw (1.7,0.7)-- (1,0.7);
\draw (1,0.7)-- (1,0);
\draw (1,1)-- (1.7,1);
\draw (1.7,1)-- (1.7,1.7);
\draw (1.7,1.7)-- (1,1.7);
\draw (1,1.7)-- (1,1);
\draw (0,1)-- (0.7,1);
\draw (0.7,1)-- (0.7,1.7);
\draw (0.7,1.7)-- (0,1.7);
\draw (0,1.7)-- (0,1);
\draw (1,0)-- (2,0);
\draw (2,0)-- (2,1);
\draw (2,1)-- (1,1);
\draw (1,1)-- (1,0);
\draw (1,1)-- (2,1);
\draw (2,1)-- (2,2);
\draw (2,2)-- (1,2);
\draw (1,2)-- (1,1);
\draw (0,1)-- (1,1);
\draw (1,1)-- (1,2);
\draw (1,2)-- (0,2);
\draw (0,2)-- (0,1);
\draw [dashed] (1,1.7)-- (0.7,1.1);
\draw [dashed] (1.7,1)-- (1.1,0.7);
\draw (0.5,0.5) node{$\Lambda$};
\draw (1.7,0) node[below]{$(2-3\varepsilon)n$};
\draw (1,0) node[below]{$n$};
\end{tikzpicture}}
\end{minipage}
\leftcaption{\label{fig:exponential:renormalisation:1}A good box in the renormalisation of the proof of Theorem~\ref{th:exponential}. The thick box is emptied in \eqref{eq:expo:renorm:1}. This allows emptying the shaded one in \eqref{eq:expo:renorm:2}. The dashed lines with slopes $1/2$ and $2$ indicate how empty sites propagate via the rule $\{(2,1),(1,2)\}$.}\rightcaption{\label{fig:exponential:renormalisation:2}If the three shaded boxes are empty, we are able to empty the thick box $\Lambda$, as indicated by the dashed lines with slopes $1/2$ and $2$, thanks to the update rule $\{(2,1),(1,2)\}$.}
\end{figure}

We turn to the second step. Fix $\varepsilon>0$ small enough depending on $\varepsilon_0$ and then take $n\in\bbN$ large enough. Let $\Lambda=\{0,\dots,n-1\}^2$ be the renormalisation box. Then the exponential decay provided by \ref{item:expo:6} and a union bound give
\begin{equation}
\label{eq:expo:renorm:1}\mu_q\left(\left[\omega_\Lambda\right]^{\bone_{\bbZ^{{2}}\setminus\Lambda}}_{[\varepsilon n,(1-\varepsilon)n)^2}=\bzero_{[\varepsilon n,(1-\varepsilon)n)^2}\right)\ge 1-\varepsilon_0.\end{equation}
In words, it is likely that BP in $\Lambda$ empties all of $\Lambda$ except a thin frame (see Figure~\ref{fig:exponential:renormalisation:1}), in which case we say that $\Lambda$ is \emph{good}. Thus, \eqref{eq:expo:renorm:1} states that the probability that a box is good is at least $1-\varepsilon_0$. However, since $\{(2,1),(1,2)\}\in\cU$, we have
\begin{equation}
\label{eq:expo:renorm:2}\left[\bzero_{[\varepsilon n,(1-\varepsilon)n)^2}\cdot\bone_{\Lambda\setminus[\varepsilon n,(1-\varepsilon)n)^2}\right]^{\bone_{\bbZ^{{2}}\setminus\Lambda}}_{[0,\dots,(1-3\varepsilon)n)^2}=\bzero_{[0,\dots,(1-3\varepsilon)n)^2},\end{equation}
that is, BP empties all but a thin strip along the top and right boundaries of $\Lambda$, as desired (see Figure~\ref{fig:exponential:renormalisation:1}). See \cite{Hartarsky21}*{Section~7.4} for more details on the second step.

We move on to the third step. Again using $\{(2,1),(1,2)\}\in\cU$, we have
\begin{equation}
\label{eq:expo:renorm:3}\left[\bzero_{[0,(1-3\varepsilon)n)+nU'}\cdot\bone_{[0,2n)^2\setminus([0,\dots,(1-3\varepsilon)n)+nU')}\right]_\Lambda^{[0,2n)^2}=\bzero_\Lambda\end{equation}
(see Figure~\ref{fig:exponential:renormalisation:2}).
Consequently, if the three boxes of the form $\Lambda+n\vec u$ with $\vec u\in U'$ are good, then $\Lambda$ can be emptied by BP in the union of these four boxes. Then Lemma~\ref{lem:legal:path} provides a legal path in $[0,2n)^2$ from $\omega_{[0,2n)^2}$ to $\omega_{[0,2n)^2}^{\vec x}$ for any $\vec x\in\Lambda$. Using this path and the canonical path approach as in the proof of Theorem~\ref{th:ergo}, this yields that for some $\gamma<\infty$ depending on $n$ and $q$ and any local function $f$,
\begin{align*}
\mu\left(\1_{\forall\vec u\in U',\Lambda+n\vec u\text{ is good}}\var_{\Lambda}(f)\right)\le& \sum_{\vec x\in\Lambda}\mu\left(\1_{\forall\vec u\in U',\Lambda+n\vec u\text{ is good}}\var_{\vec x}(f)\right)\\
\le&\gamma\sum_{\vec x\in[0,2n)^2}\mu(c_{\vec x}\var_{\vec x}(f)).\end{align*}
Finally, it remains to see that for some $\gamma'<\infty$ and any local function $f$,
\[\var(f)\le \gamma'\sum_{\vec x\in\bbZ^{{2}}}\mu\left(\1_{\forall\vec u\in U',n\vec x+\Lambda+n\vec u\text{ is good}}\var_{n\vec x+\Lambda}(f)\right).\]

Recalling \eqref{eq:def:cx} and \eqref{eq:poincare}, we recognise exactly the Poincar\'e inequality for the $\cU'$-KCM except that each site (box) now has more than two states and the parameter $q$ is replaced by $\mu_q(\Lambda\text{ is good})\ge 1-\varepsilon_0$. We will frequently encounter such general state-space KCM arising from renormalisation procedures (see Section~\ref{sec:general:KCM}) and will see that usually the corresponding Poincar\'e inequalities are proved just like the ones for ordinary KCM. Since $\trel<\infty$ for the $\cU'$-KCM at $q\ge 1-\varepsilon_0$, by the first step, this completes the proof. See \cite{Cancrini08}*{Section~5} for more details.\\
\textbf{\ref{item:expo:7} implies \ref{item:expo:2}.}
Define 
\[\lambda_0=\inf\left\{\cD(f):f\in\mathrm{L}^2(\mu),\mu\left(f^2\right)=1,\forall \omega\in\Omega,\omega_0=0\Rightarrow f(\omega)=0\right\}.\]
Observe that any $f$ as above satisfies $\var(f)\ge q$, since $\mu^2(f)=(1-q)^2\mu^2(f|\omega_0=1)\le (1-q)^2\mu(f^2|\omega_0=1)=(1-q)$. Recalling \eqref{eq:def:Trel}, this implies 
$\lambda_0\ge q/\trel$. Finally, a general result on hitting times for Markov chains \cite{Aldous02a}*{Proposition~3.21} yields
\begin{equation}
\label{eq:F0F1}
\bbP(\tau_0>t)\le e^{-t\lambda_0}\le e^{-qt/\trel}\end{equation}
for any $t\ge 0$. This proves \ref{item:expo:4} and one can recover \ref{item:expo:2}, by proceeding analogously for $\tau_1$ (recall \eqref{eq:def:tau}). See \cite{Cancrini08}*{Theorem~3.6} for a different proof.
\end{proof}

The next result \cite{Cancrini08}*{Lemma 2.11, Proposition~2.13} reduces the relaxation time of a KCM in infinite volume to the one in finite volume (recall Section~\ref{sec:finite_vol}).
\begin{proposition}[Finite volume reduction]
    \label{prop:infinite:to:finite}Let $\cU$ be an update family and $q\in(0,1)$. Given $\Lambda\Subset\bbZ^d$, recall that $\trel^{\Lambda}$ denotes the relaxation time of the finite-volume generator $\cL_{\bzero_{\bbZ^d\setminus\Lambda}}$. 
    Then
    \[\trel=\lim_{\Lambda\Subset\bbZ^d,\Lambda\to\bbZ^d}\trel^\Lambda
    .\]
    The same holds for $C_{\mathrm{LS}}$ and $C_{\mathrm{MLS}}$.
\end{proposition}
\begin{proof}We prove the chain of inequalities
\begin{equation}
\label{eq:finite:volume:chain}
\trel\le\sup_{\Lambda\Subset\bbZ^d} \trel^\Lambda\le 
\lim_{\Lambda\Subset\bbZ^d,\Lambda\to\bbZ^d}\trel^\Lambda\le \trel.
\end{equation}
For the first one, it suffices to prove the Poincar\'e inequality \eqref{eq:poincare} with $C=\sup_{\Lambda\Subset\bbZ^d}\trel^\Lambda$ for any local function $f$. Since $f$ is local, for $\Lambda$ large enough we have $\var(f)=\var_{\Lambda}(f)$ and $\cD(f)=\cD_{\bzero_{\bbZ^d\setminus\Lambda}}(f)$. So we are done by \eqref{eq:def:Trel}.

For the second and third inequalities and the existence of the limit in \eqref{eq:finite:volume:chain}, it suffices to show that for $V\subset\Lambda\subset
\bbZ^d$ we have
\begin{equation}
\label{eq:monotonicity}
\trel^V\le \trel^\Lambda.
\end{equation}
Consider $f:\Omega_V\to\bbR$ and extend $f$ to $\tilde f:\Omega_\Lambda\to \bbR$ by $\tilde f(\omega)=f(\omega_V)$. Clearly, $\var(f)=\var(\tilde f)$, since $\mu_\Lambda=\mu_V\otimes\mu_{\Lambda\setminus V}$, so it remains to check that \[\cD_{\bzero_{\bbZ^d\setminus\Lambda}}(\tilde f)\le\cD_{\bzero_{\bbZ^d\setminus V}}(f).\]
But, recalling \eqref{eq:dirichlet}, this follows from $c_\bx^{\bzero_{\bbZ^d\setminus\Lambda}}\le c_\bx^{\bzero_{\bbZ^d\setminus V}}$. The latter is true, using \eqref{eq:def:cx:finite} and the fact that $c_\bx$ is non-increasing.

The proof for the (modified) logarithmic Sobolev constant is identical.
\end{proof}

\section{Stronger functional inequalities}
We next consider the time scales corresponding to mixing, modified logarithmic Sobolev and logarithmic Sobolev inequalities (recall Section~\ref{sec:time_scales}). We start with general bounds on the mixing time.
\begin{proposition}[{Basic mixing time bounds}]
\label{prop:tmix:easy}
For any update family $\cU$ and $q>0$ there exists $C>0$ such that the following holds. If $q>\qct$, then for all $\varepsilon\in(0,1)$ and $\Lambda\subset\bbZ^d$,
\begin{equation}
\label{eq:tmix:upper:volume}
t_{\mathrm{mix}}^{\bzero_{\bbZ^d\setminus\Lambda}}(\varepsilon)\le C(\log(1/\varepsilon)+|\Lambda|).\end{equation}
For all $\varepsilon\in(0,1)$ and all $n$ large enough, for any $\sigma\in\Omega_{\bbZ^d\setminus\{1,\dots,n\}^d}$,
\begin{equation}\label{eq:tmix:lower:linear}
t_{\mathrm{mix}}^{\sigma}(\varepsilon)\ge n/C.\end{equation}
\end{proposition}
\begin{proof}
For \eqref{eq:tmix:upper:volume}, we use a general inequality for reversible Markov chains (see e.g.\ \cite{Levin09}*{Theorem~20.6})
\[t_{\mathrm{mix}}^{\bzero_{\bbZ^d\setminus\Lambda}}(\varepsilon)\le \log\left(\frac{1}{\varepsilon\mu_*}\right)\trel^\Lambda,\]
with $\mu_*=\min_{\omega\in\Omega_{\Lambda}}\mu(\omega)=\min(q,1-q)^{|\Lambda|}$. By \eqref{eq:def:qct}, for $q>\qct$, we have $\trel<\infty$. Moreover, by \eqref{eq:monotonicity}, $\trel\ge \trel^{\bzero_{\bbZ^d\setminus\Lambda}}$, concluding the proof of \eqref{eq:tmix:upper:volume}.

For \eqref{eq:tmix:lower:linear}, we choose $\omega=\bone_\Lambda$ in \eqref{eq:def:tmix} with $\Lambda=\{1,\dots,n\}^d$. We denote the $\cU$-KCM with parameter $q$, boundary condition $\sigma\in\Omega_{\bbZ^d\setminus\Lambda}$ and initial condition $\omega$ by $(\eta(t))_{t\ge0}$. Fix $C>0$ large enough and set $\Lambda'=\{\lceil n/\sqrt C\rceil,\dots,n-\lceil n/\sqrt C\rceil\}^d$. It suffices to show that for any $q\in[0,1]$
\begin{equation}
\label{eq:speed:propagation}
\lim_{n\to\infty}\bbP\left(\eta_{\Lambda'}(n/C)=\bone_{\Lambda'}\right)=1.\end{equation}

Let us denote $\|\cU\|=\max\{\|u\|:\exists U\in\cU,u\in U\}$. A \emph{path from $\bx\in\Lambda'$ to $\partial\Lambda$} is a finite sequence of sites $(\bx_i)_{i=1}^N$ in $\Lambda$ with $\bx_1=\bx$, $\bx_N$ at distance at most $\|\cU\|$ from $\bbZ^d\setminus\Lambda$, and, for each $i\in\{2,\dots,N\}$, $\|\bx_i-\bx_{i-1}\|\le\|\cU\|$. Assume that $\eta_t(\bx)=0$ for some $t>0$. Then we can construct a path from $\bx$ to $\partial\Lambda$ inductively as follows, assuming that $\bx_i$ is chosen (by definition, $\bx_1=\bx$). If $d(\bx_i,\bbZ^d\setminus\Lambda)\le \|\cU\|$, we set $N=i$ and stop. Otherwise, consider the first time $t'\le t$ such that $\eta_{t'}(\bx_i)=0$ and let $\bx_{i+1}$ be arbitrarily chosen in $\bx+\bigcup_{U\in\cU}U$ so that $\eta_{\bx_{i+1}}(t'-)=0$. We call this a \emph{decreasing path of $\bx$}, since the emptying times of its sites are decreasing. For a given path from $\bx$ to $\partial\Lambda$ of length $k$, the probability that it is decreasing is at most the probability that the sum of $k$ i.i.d.\ standard exponential random variables $E_i$ is  at most $t$. Indeed, this follows from the graphical construction of Section~\ref{subsec:markov}. Moreover, $k\ge (n/\sqrt C)/\|\cU\|$, so, setting $t=n/C$, for $n$ large enough we obtain
\begin{align*}
\bbP(\eta_t(\bx)=0)&{}\le \sum_{k\ge \frac{n}{\|\cU\|\sqrt C}}\left|\{(\bx_i)_i\text{ path from $\bx$ to $\partial\Lambda$ of length $k$}\}\right|\cdot\bbP\left(\sum_{l=1}^kE_l\le t\right)\\
&{}\le \sum_{k\ge \frac{n}{\|\cU\|\sqrt C}} (2\|\cU\|+1)^{kd}\cdot\bbP\left(\sum_{l=1}^kE_l\le \frac{k}{C^{1/3}}\right)\\
&{}\le \sum_{k\ge \frac{n}{\|\cU\|\sqrt C}}(2\|\cU\|+1)^{kd}\cdot\frac{(\bbE(e^{-C^{1/3}E_1}))^k}{e^{-k}}\le e^{-n/(\sqrt C\|\cU\|)},
\end{align*}
where we used the exponential Markov inequality in the last line and took into account that $C$ can be chosen large enough depending on $\|\cU\|$ and $d$. The desired \eqref{eq:speed:propagation} then follows by a union bound over $\bx\in\Lambda'$.
The argument used in the proof \eqref{eq:tmix:lower:linear} is commonly referred to as \emph{finite speed of propagation} and applies equally well to any finite-range interacting particle system \cite{Liggett99}*{Section~I.1}. In fact, this is a way to show that the graphical construction of the KCM produces a well-defined Markov process.
\end{proof}

\begin{corollary}[{Infinite logarithmic Sobolev constants}]
\label{cor:LS}
    For any update family $\cU$ and $q\in(0,1)$, we have \begin{equation}\label{inf}C_{\mathrm{LS}}=C_{\mathrm{MLS}}=\infty.\end{equation}
    Furthermore, there exists $C=C(q)>0$ such that the following holds \begin{itemize}\item[(i)]   
Fix  $n$ and let  $\Lambda=\{-n,\dots,n\}^d$. For any $\sigma\in\Omega_{\bbZ^d\setminus\Lambda}$ and for $n$ large enough, we have
\begin{equation}\label{eq:ls:lower:linear}
C_{\mathrm{LS}}^{\sigma}\ge C_{\mathrm{MLS}}^{\sigma}\ge n/C.\end{equation}
\item[(ii)]
If $q>\qct$, then, for any $\Lambda\Subset\bbZ^d$, we have
\begin{equation}
\label{eq:ls:upper:volume}
C_{\mathrm{LS}}^{\Lambda}\leq C|\Lambda|.
\end{equation}
\end{itemize}
\end{corollary}
\begin{proof}
The first inequality in \eqref{eq:ls:lower:linear} is a general and classical result
(see \cite{Diaconis96}*{Lemma 2.7}). To prove the second inequality we set 
\[f(\omega)=(1-q)^{-|\Lambda|}\1_{\omega_\Lambda=\bone_\Lambda},\]
so that $\mu(f)=1$. Let $\mu^f$ be the probability measure with density $f$ w.r.t.\ $\mu$. Then, by \eqref{eq:speed:propagation}, for a properly chosen $C>0$, we have
\begin{equation}
\label{cont}
\lim_{n\to\infty}\mathbb E_{\mu^f}(\eta_{0}(Cn))=1.
\end{equation}
Combining \eqref{eq:rel_ent_dec_1}, \eqref{cont} and
Pinsker's inequality implies the second inequality of \eqref{eq:ls:lower:linear}.

Inequality \eqref{eq:ls:upper:volume}
follows from Theorem~\ref{th:exponential}\ref{item:expo:7} together with \eqref{eq:finite:volume:chain} and the general and classical bound (see \cite{Saloff-Coste97} Corollary 2.2.10)
\begin{equation}
C_{\mathrm{LS}}^{\sigma}    \le \frac{\log((1/\mu^*)-1)}{1-2\mu^*}\trel^{\sigma}
\end{equation}
where $\mu^*:=\min_{\omega\in\Omega_{\Lambda}}
\mu(\omega)$.
Finally, \eqref{inf} follows from \eqref{eq:ls:lower:linear} and Proposition
~\ref{prop:infinite:to:finite}.
\end{proof}

\section{Conclusion}
\label{subsec:tools:BP}
In Theorems~\ref{th:ergo} and~\ref{th:exponential} we completely reduced the critical values $\qc$ and $\qct$ to their BP counterparts. As we saw in Theorem~\ref{th:BP} and will see more generally in Chapter~\ref{chap:universality}, modulo Conjecture~\ref{conj:qc:qct}, this gives complete information about $\qc$ and $\qct$. Our next goal is to find the asymptotics of $\tau_0$ as $q\to\qc+$, as in the BP Theorem~\ref{th:BP}. In reality, even in the simplest BP models with $\qc\in(0,1)$, such as North-East BP, the asymptotics of $\tbp$ as $q\to \qc+$ and of $\mu_{\qc}(\tbp>t)$ as $t\to\infty$ remain inaccessible, even though it is classically conjectured that they should be governed by certain critical exponents. Therefore, in the remainder of the manuscript we mostly focus on update families for which $\qc=0$.

Before we direct our efforts to seeking asymptotics as $q\to0$, let us summarise the techniques we learned in this chapter, as several of them will be of further use. 

\noindent\textbf{BP lower bound.} In \eqref{eq:t0t0bp} we have seen that BP provides a natural lower bound for KCM time scales.  Let us also record the general quantitative bound obtained by combining \eqref{eq:t0t0bp} and \eqref{eq:F0F1}:
\begin{equation}\label{eq:tau0}
\delta\mu_q\left(\tbp\right)\le \bbE_{\mu_q}(\tau_0)\leq \frac{\trel}{q}
\end{equation}
for some $\delta>0$ depending only on $\cU$ (not necessarily equal to the one in \eqref{eq:t0t0bp}).

\noindent\textbf{Test functions.} Thanks to the variational definition \eqref{eq:def:Trel} of $\trel$, one can obtain lower bounds by plugging well-chosen functions. We used this in the proof of the implication \ref{item:ergo:4} to \ref{item:ergo:1} of Theorem~\ref{th:ergo}. The test function we used here simply reflected BP, but in the next chapter we will need a more subtle choice.

\noindent\textbf{Canonical paths.} This method for proving upper bounds on $\trel$ was used in the implication \ref{item:ergo:1} to \ref{item:ergo:2} of Theorem~\ref{th:ergo}. In this instance, the canonical paths were simply granted by BP, but in the next chapter, paths will incorporate more refined heuristics on the KCM dynamics.

\noindent\textbf{Renormalisation.} In the proof of the implication \ref{item:expo:6} to \ref{item:expo:7} we used the idea of renormalisation. We regarded large boxes as single sites in an auxiliary (generalised) KCM dynamics. This splits the problem of proving upper bounds on $\trel$ in two. First, we prove an upper bound on the relaxation time of the auxiliary dynamics. Second, we show (e.g.\ by canonical paths) how to ``locally'' reconstruct the original Dirichlet form from the auxiliary one, taking advantage of the auxiliary constraint. We will make extensive use of this technique in the next chapters.

\noindent \textbf{Finite speed of propagation.} This observation allows us to show that with very high probability no information about the state of the process at a given place or its boundary condition can travel faster than linearly. The proof of the lower bound of Proposition~\ref{prop:tmix:easy} was an application of this fact.

\noindent \textbf{Stronger functional inequalities.} 
In Corollary~\ref{cor:LS}, we established that, contrary to Poincar\'e inequalities, logarithmic and modified logarithmic Sobolev inequalities in infinite volume are \emph{not} a suitable tool for studying KCM. Nonetheless, such techniques can be of use on suitably chosen finite volumes.


	\chapter{One-dimensional models}

\label{chap:1d}

\abstract{In this chapter we investigate one-dimensional KCM. Most notably, these include the FA-1f and East models. We present the techniques used to determine the scaling of their characteristic times as $q\to 0$. We familiarise ourselves with the use of test functions, non BP-based canonical paths, combinatorial bottlenecks and bisection in the simplest possible setting. We then move on to FA-2f and general KCM, still in one dimension. These one-dimensional models are not only interesting in their own right, but will also serve as tools for the study of higher-dimensional KCM via renormalisation.\\}

In one dimension there are three nearest neighbour KCM corresponding to the update families $\{\{-1\},\{1\}\}$, $\{\{1\}\}$ and $\{\{-1,1\}\}$ (excluding the trivial cases $\{\varnothing\}$ and $\varnothing$, as usual). We recognise the FA-1f, East and FA-2f models respectively. The latter is not very interesting from our viewpoint, but we consider it for completeness. FA-1f and East on the other hand are not only of interest themselves, but also provide fundamental building blocks for renormalisation arguments for more complex models (as in the proof of Theorem~\ref{th:exponential}).

\section{FA-1f}
\label{sec:FA1f}
In this section, we consider $\cU=\{\{-1\},\{1\}\}$. From Theorem~\ref{th:BP}, we have $q_c=0$. We therefore look for asymptotics as $q\to0$, which are provided by the work of Cancrini, Martinelli, Roberto and Toninelli \cite{Cancrini08} and Shapira \cite{Shapira20} (and its arXiv version).
\begin{theorem}[FA-1f asymptotics]
\label{th:FA1f}
For FA-1f in $d=1$ dimension there exists $C>0$ such that for $q$ small enough
\begin{align*}1/C\le q^{3}\trel\le{}&C,&1/C\le q^3\bbE_{\mu_q}[\tau_0]\le{}&C,&\bbP_\mu(\tau_0>t)\le{}& e^{-Cq^3t}.\end{align*}
\end{theorem}
\begin{remark}[{Arrhenius law}]
\label{rem:Arrhenius}
Note that if we rewrite Theorem~\ref{th:FA1f} in terms of the inverse temperature, using \eqref{eq:temperature}, the above scaling corresponds to an Arrhenius divergence, $\tau_0,T_{\rm rel}\approx \exp(c\beta)$, as for strong glass forming liquids (see Figure~\ref{fig:vetri:1}).
\end{remark}
Before embarking on the proof, let us explain the heuristics behind Theorem~\ref{th:FA1f}. For $q$ small empty sites are typically isolated (and therefore cannot be immediately removed). However, if $x\in\bbZ$ is empty, at rate $q$ we can empty $x-1$ or $x+1$.  Suppose the first event occurs. At this point the constraint at $x$ is also satisfied and, in a time of order one, we will (with equal probability) either occupy $x-1$ or $x$. In the latter case, the  net result is that we have ``moved'' the empty site from $x$ to $x-1$. So, we intuitively expect empty sites to behave like  random walks moving at rate $q$ until they meet. When they meet, they typically coalesce. This explains the scaling $1/q^3$: it is the time required to overcome the typical distance $\ell=1/q$ between two consecutive empty sites (inverse rate times distance squared, that is $1/q \times 1/q^2$).

\begin{proof}[Theorem~\ref{th:FA1f}]
In order to show that $\trel\ge q^{-3}/C$, we recall  \eqref{eq:def:Trel}. Define the test function $f(\omega)=\min\{k\ge 1:\omega_{k}\omega_{-k+1}=0\}$, that is, the distance from $1/2$ to the nearest empty site rounded up. As in the proof of Theorem~\ref{th:BP}, one can check that $\mu(f(\omega)=k+1)=(2q-q^2)(1-q)^{2k}$ for $k\ge 0$. This geometric random variable has $\var(f)=(1-q)^2/(2q-q^2)^2$. On the other hand, by \eqref{eq:dirichlet},
\begin{align*}\cD(f)={}&2\sum_{x\ge 1} \mu(c_x\var_x(f))=2q(1-q)\sum_{x\ge 1}\mu\left((1-\omega_{x+1})\prod_{y=-x+1}^{x-1}\omega_y\right)\\={}&
2q^2\sum_{x\ge 1}(1-q)^{2x}=\frac{q(1-q)^2}{1-q/2}.\end{align*}
Hence, by \eqref{eq:def:Trel},
\[\trel\ge \frac{\var(f)}{\cD(f)}=\frac{1}{4q^3(1-q/2)}\ge \frac{1}{4q^{3}}.\]

The inequality $\bbE_{\mu}[\tau_0]\ge 1/(Cq^3)$ is proved in a similar way, but using a more subtle variant of \eqref{eq:def:Trel} regarding hitting times. See \cite{Shapira20}*{Section~4.2} and \cite{Shapira20a}.

The proof of $\trel\le C/q^3$ proceeds similarly to the implication from \ref{item:expo:6} to \ref{item:expo:7} in Theorem~\ref{th:exponential}. However, we use canonical paths reflecting the heuristic mechanism discussed above rather than the brutal legal paths provided by bootstrap percolation in Lemma~\ref{lem:legal:path}. The first step is proving a Poincar\'e inequality for the generalised East model, which will be discussed in Section~\ref{sec:general:KCM}. It implies that there exists a constant $C<\infty$ independent of $q$ such that for any local function $f$,
\begin{equation}
\label{eq:East:Poincare:renorm}\var(f)\le C\sum_{x\in\bbZ}\mu(\tilde c_x\var_{\Lambda_x}(f)),\end{equation}
where $\Lambda_x=\{x\lceil 1/q\rceil,\dots,(x+1)\lceil1/q\rceil-1\}$ and $\tilde c_x=1-\prod_{y\in\Lambda_{x+1}} \omega_{y}$. This is precisely the unidimensional analogue of the renormalisation of Figure \ref{fig:exponential:renormalisation:2}.

In view of \eqref{eq:East:Poincare:renorm}, we seek to bound $\mu(\tilde c_x\var_{\Lambda_x}(f))$ with terms $\mu(c_y\var_y(f))$ of the Dirichlet form from \eqref{eq:dirichlet}. To do so, for each $x\in\bbZ$, $\omega\in\Omega$ such that $\omega_{\Lambda_{x+1}}\neq \bone$ and $y\in\Lambda_x$ with $\omega_y=1$, we define a legal path $(\omega^{(i)})_{i=0}^l$ from $\omega$ to $\omega^y$ as follows. Set $\xi=\min\{z>y:\omega_z=0\}$, $l=2(\xi-y-1)+1$ and define $\omega^{(0)}=\omega$, $\omega^{(l)}=\omega^y$ and 
\begin{align}
\label{eq:path:FA1f:1d}
\omega^{(2i-1)}={}&\omega^{\xi-i},&
\omega^{(2i)}={}&(\omega^{\xi-i})^{\xi-i-1}
\end{align}
for $i\in \{1,\dots,\xi-y-1\}$ ($\omega^{(0)}=\omega$ and $\omega^{(l)}=\omega^y$ by definition). In words, this is the legal path sending an empty interval of length oscillating between one and two sites from $\xi$ to $y$. Let $x^{(i)}$ be the site such that $\omega^{(i)}=(\omega^{(i-1)})^{x^{(i)}}$. Observe that for odd $i\in\{3,\dots,l\}$ we have $\omega_{x^{(i)}}^{(i)}=1=\omega^{(i)}_{x^{(i+1)}}$, so it is convenient to set $j_i=2\lceil i/2\rceil-1$ (that is, the odd number in $\{i,i-1\}$) for $i\in\{2,\dots,l\}$ and $j_1=0$. 

Thus, by the Cauchy--Schwarz inequality, any $\omega$ such that $\tilde c_x(\omega)=1$ satisfies
\begin{align*}
q\omega_y\left(f(\omega^y)-f(\omega)\right)^2&{}\le ql\sum_{i=1}^{l}c_{x^{(i)}}\left(\omega^{(i)}\right)\left(f\left(\omega^{(i)}\right)-f\left(\omega^{(i-1)}\right)\right)^2,\\
&{}=ql\sum_{i=1}^l\omega^{(j_i)}_{x^{(i)}}c_{x^{(i)}}\left(\omega^{(j_i)}\right)\left(f\left(\omega^{(j_i)}\right)-f\left(\left(\omega^{(j_i)}\right)^{x^{(i)}}\right)\right)^2.\end{align*}
Next note that, given $y$, $\omega^{(j_i)}$ and whether $i=1$ or not, we can recover $\omega$ and that $\mu(\omega)/\mu(\omega^{(j_i)})\le (1-q)/q$ (equality holds except for $i=1$). Integrating the last display over $\omega$, we obtain
\[\mu(\tilde c_x\var_y(f))\le 4\lceil 1/q\rceil\sum_{z\in\Lambda_x\cup\Lambda_{x+1}}\mu\left(\tilde c_x(\omega)c_{z}(\omega)\omega_z\left(f(\omega)-f(\omega^z)\right)^2\right),\]
since $l\le 2\lceil 1/q\rceil$. Summing the last result over $y\in \Lambda_x$ and then $x\in\bbZ$, we obtain
\[\sum_{x\in\bbZ}\mu(\tilde c_x\var_{\Lambda_x}(f))\le \sum_{x\in\bbZ}\mu\left(\tilde c_x\sum_{y\in\Lambda_x}\var_y(f)\right)\le 8\lceil1/q\rceil^3\cD(f),\]
recalling \eqref{eq:dirichlet} and \eqref{eq:Poincare:unconstrained}. Plugging this into \eqref{eq:East:Poincare:renorm} and recalling \eqref{eq:def:Trel} concludes the proof.

The proof that $\bbP_\mu(\tau_0>t)\le e^{-Cq^3t}$ follows similar but more delicate lines (see \cite{Shapira20}*{Section~4.1} for the details). Finally the upper bound on $\bbE_{\mu}[\tau_0]$ follows directly from the last inequality.
\end{proof}
We note that higher dimensional analogues of Theorem~\ref{th:FA1f} are known and the scaling is $\log(1/q)/q^{2}$ in $d=2$ and $q^{-2}$ in $d\ge 3$ \cites{Cancrini08,Shapira20} with the exception of the following conjecture, whose lower bound remains open.
\begin{conjecture}[{Relaxation time in two dimensions}]
\label{conj:FA1f}
For FA-1f in $d=2$ dimensions, there exists $C>0$ such that for $q$ small enough $1/C\le q^2\trel/\log(1/q)\le C$.
\end{conjecture}

\section{East}
Recall that the East model in one dimension corresponds to the update family $\cU=\{\{1\}\}$. It is, in a sense, the simplest non-trivial KCM and lies at the base of the theory. Indeed, the very first rigorous results on KCM were proved for the East model around the turn of the century \cites{Chung01,Aldous02}. Furthermore, the East model appears also in other contexts including the study of  random walks on the group of  upper triangular matrices \cites{Peres13,Ganguly19}.

Once again, by Theorem~\ref{th:BP}, we have $\qc=0$, so we are interested in asymptotics as $q\to 0$. We address upper and lower bounds separately, showcasing two important techniques. Together they give the following result of Aldous and Diaconis, and Cancrini, Martinelli, Roberto and Toninelli \cites{Aldous02, Cancrini08}.

\begin{theorem}[East asymptotics]
\label{th:East}
For the East KCM in $d=1$ dimension we have
\begin{align*}
\lim_{q\to0}\frac{\log \trel}{(\log(1/q))^2}&{}=\frac{1}{2\log 2},
&\lim_{q\to 0}\bbP_\mu\left(\left|\frac{2\log 2\cdot \log \tau_0}{(\log(1/q))^{2}}-1\right|<\varepsilon\right)&{}=1
\end{align*}
for any $\varepsilon>0$.
\end{theorem}
We refer the reader to \cite{Chleboun14} for finer results including the scaling of $T_{\rm rel}^{\Lambda}$ as a function of $q$ and $|\Lambda|$ and for the equivalence, up to a length scale $|\Lambda|=O(1/q)$, of the relaxation and mixing time of the East model; to \cite{Faggionato13} for a survey dedicated to this model; to \cite{Chleboun16} for higher dimensions; to \cite{Couzinie22b} for a multicolour version.

\begin{remark}[{Super-Arrhenius law}]
\label{rem:superArrhenius}
If we rewrite Theorem~\ref{th:East} in terms of temperature using \eqref{eq:temperature}, we get $\tau_0,\trel\approx \exp(\beta^2/(2\log 2))$. This super-Arrhenius divergence is reminiscent of the scaling for fragile super-cooled liquids (see Figure~\ref{fig:vetri:1}), which is an important reason for interest in this model among physicists. Also note that this scaling diverges much faster than the emptying time for the corresponding BP (recall Theorem~\ref{th:BP}) or FA-1f (recall Theorem~\ref{th:FA1f} and Remark~\ref{rem:Arrhenius}).
\end{remark}

\subsection{Lower bound: combinatorial bottleneck}
\label{subsec:East:lower}
The lower bound is proved via a \emph{combinatorial bottleneck}. In rough terms, the strategy is as follows. We consider the stationary KCM started at a typical configuration under $\mu_q$ with $q$ small. We identify a set of configurations around the origin which (deterministically) cannot be avoided if the origin is to be infected. For instance, this could be having an atypically large number of empty sites in the vicinity of the origin, but will usually be more subtle. We then seek to evaluate the probability and the number (entropy) of these \emph{bottleneck configurations}. Finally, we use stationarity and a union bound on time to show that if these configurations are unlikely and there are few of them, thus one needs to wait a lot of time to observe any of them close to the origin. The hard part of such arguments is identifying the correct bottleneck.

\subsubsection{Combinatorics for the East model}
The key ingredient to understand the behaviour of  the one-dimensional East model is the following combinatorial result of Chung, Diaconis and Graham \cite{Chung01}.

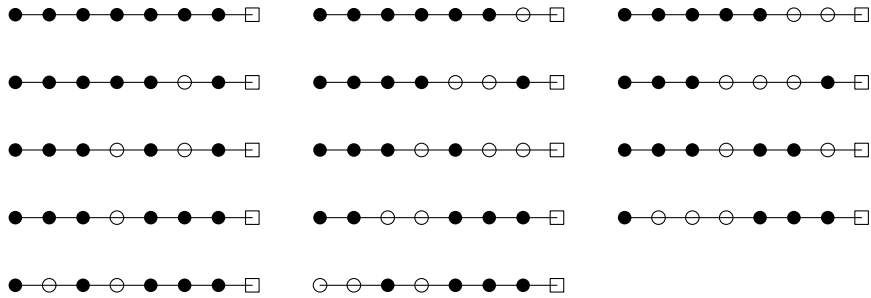
\begin{figure}
\begin{subfigure}{0.25 \textwidth}
\begin{tikzpicture}[line cap=round,line join=round,x=2.0cm,y=2.0cm, scale=0.45]
\draw(0,0) -- (3.5,0);
\fill  (0, 0) circle (2mm);
\fill (0.5, 0) circle (2mm);
\fill (1, 0) circle (2mm);
\fill (1.5, 0) circle (2mm);
\fill  (2, 0) circle (2mm);
\fill (2.5, 0) circle (2mm);
\fill (3, 0) circle (2mm);
\draw (3.4, -0.1) rectangle (3.6,0.1);
\draw(4.5,0) -- (8,0);
\fill  (4.5, 0) circle (2mm);
\fill (5, 0) circle (2mm);
\fill (5.5, 0) circle (2mm);
\fill (6, 0) circle (2mm);
\fill  (6.5, 0) circle (2mm);
\fill (7, 0) circle (2mm);
\draw (7.5, 0) circle (2mm);
\draw (7.9, -0.1) rectangle (8.1,0.1);
\draw(9,0) -- (12.5,0);
\fill  (9, 0) circle (2mm);
\fill (9.5, 0) circle (2mm);
\fill (10, 0) circle (2mm);
\fill (10.5, 0) circle (2mm);
\fill (11, 0) circle (2mm);
\draw (11.5, 0) circle (2mm);
\draw(12, 0) circle (2mm);
\draw (12.4, -0.1) rectangle (12.6,0.1);
\draw(0,-1) -- (3.5,-1);
\fill  (0, -1) circle (2mm);
\fill (0.5, -1) circle (2mm);
\fill (1, -1) circle (2mm);
\fill (1.5, -1) circle (2mm);
\fill  (2, -1) circle (2mm);
\draw (2.5, -1) circle (2mm);
\fill (3, -1) circle (2mm);
\draw (3.4, -1.1) rectangle (3.6,-0.9);
\draw(4.5,-1) -- (8,-1);
\fill  (4.5, -1) circle (2mm);
\fill (5, -1) circle (2mm);
\fill (5.5, -1) circle (2mm);
\fill (6, -1) circle (2mm);
\draw  (6.5, -1) circle (2mm);
\draw (7, -1) circle (2mm);
\fill (7.5, -1) circle (2mm);
\draw (7.9, -1.1) rectangle (8.1,-0.9);
\draw(9,-1) -- (12.5,-1);
\fill  (9, -1) circle (2mm);
\fill (9.5, -1) circle (2mm);
\fill (10, -1) circle (2mm);
\draw (10.5, -1) circle (2mm);
\draw (11, -1) circle (2mm);
\draw (11.5, -1) circle (2mm);
\fill(12, -1) circle (2mm);
\draw (12.4, -1.1) rectangle (12.6,-.9);
\draw(0,-2) -- (3.5,-2);
\fill  (0, -2) circle (2mm);
\fill (0.5, -2) circle (2mm);
\fill (1, -2) circle (2mm);
\draw (1.5, -2) circle (2mm);
\fill  (2, -2) circle (2mm);
\draw (2.5, -2) circle (2mm);
\fill (3, -2) circle (2mm);
\draw (3.4, -2.1) rectangle (3.6,-1.9);
\draw(4.5,-2) -- (8,-2);
\fill  (4.5, -2) circle (2mm);
\fill (5, -2) circle (2mm);
\fill (5.5, -2) circle (2mm);
\draw (6, -2) circle (2mm);
\fill  (6.5, -2) circle (2mm);
\draw (7, -2) circle (2mm);
\draw (7.5, -2) circle (2mm);
\draw (7.9, -2.1) rectangle (8.1,-1.9);
\draw(9,-2) -- (12.5,-2);
\fill  (9, -2) circle (2mm);
\fill (9.5, -2) circle (2mm);
\fill (10, -2) circle (2mm);
\draw (10.5, -2) circle (2mm);
\fill  (11, -2) circle (2mm);
\fill(11.5, -2) circle (2mm);
\draw(12, -2) circle (2mm);
\draw (12.4, -2.1) rectangle (12.6,-1.9);
\draw(0,-3) -- (3.5,-3);
\fill  (0, -3) circle (2mm);
\fill (0.5, -3) circle (2mm);
\fill (1, -3) circle (2mm);
\draw (1.5, -3) circle (2mm);
\fill  (2, -3) circle (2mm);
\fill (2.5, -3) circle (2mm);
\fill (3, -3) circle (2mm);
\draw (3.4, -3.1) rectangle (3.6,-2.9);
\draw(4.5,-3) -- (8,-3);
\fill  (4.5, -3) circle (2mm);
\fill (5, -3) circle (2mm);
\draw (5.5, -3) circle (2mm);
\draw(6, -3) circle (2mm);
\fill  (6.5, -3) circle (2mm);
\fill (7, -3) circle (2mm);
\fill (7.5, -3) circle (2mm);
\draw (7.9, -3.1) rectangle (8.1,-2.9);
\draw(9,-3) -- (12.5,-3);
\fill  (9, -3) circle (2mm);
\draw (9.5, -3) circle (2mm);
\draw (10, -3) circle (2mm);
\draw (10.5, -3) circle (2mm);
\fill  (11, -3) circle (2mm);
\fill(11.5, -3) circle (2mm);
\fill(12, -3) circle (2mm);
\draw (12.4, -3.1) rectangle (12.6,-2.9);
\draw(0,-4) -- (3.5,-4);
\fill  (0, -4) circle (2mm);
\draw (0.5, -4) circle (2mm);
\fill (1, -4) circle (2mm);
\draw (1.5, -4) circle (2mm);
\fill  (2, -4) circle (2mm);
\fill (2.5, -4) circle (2mm);
\fill (3, -4) circle (2mm);
\draw (3.4, -4.1) rectangle (3.6,-3.9);
\draw(4.5,-4) -- (8,-4);
\draw  (4.5, -4) circle (2mm);
\draw (5, -4) circle (2mm);
\fill (5.5, -4) circle (2mm);
\draw (6, -4) circle (2mm);
\fill  (6.5, -4) circle (2mm);
\fill (7, -4) circle (2mm);
\fill (7.5, -4) circle (2mm);
\draw (7.9, -4.1) rectangle (8.1,-3.9);
\end{tikzpicture}
\end{subfigure}
\caption{\label{fig:East_path} A legal path on $\{\dots,-2-1\}$ with at most $3$ simultaneous empty sites starting from a completely occupied configuration and creating a vacancy at $-(2^3-1)$. Successive steps of the path should be read from left to right and top to bottom. The empty square at site $0$ stands for the empty boundary condition.}\end{figure}

\begin{proposition}[{Combinatorial bottleneck for East}]
\label{prop:comb}
Consider the East model on $\Lambda=\{\dots,-2,-1\}\subset \mathbb Z$ with boundary condition $\bzero_{\{0,1,\dots\}}$. For $n\ge 0$, let $V(n)$ be the set of all configurations that the process can reach from $\bone_\Lambda$ via a legal path (recall Definition~\ref{def:legal:path}) in which all configurations contain 
at most $n$ empty sites. For $k\in\{0,\dots,n\}$, let $V(n,k)=\{\omega\in V(n):\sum_{x\in\Lambda}(1-\omega_x)=k\}$ be the configurations in $V(n)$ with $k$ empty sites. Finally, let $\ell(n)$ be the largest distance of an empty site (in $\Lambda$) from $0$ in $V(n)$, that is
\[
\ell(n)=\sup\left\{-y: y \in \Lambda, \exists \omega\in V(n),\omega_y=0\right\}
\]
with the convention $\sup\varnothing=0$. For $n \geq 0$,
\begin{align}
\label{eq:combi:1}
\ell(n)&{}=2^{n}-1,\\
\label{eq:combi:2} \left|V(n,n)\right|&{}\leq n!2^{\binom{n}{2}} .
\end{align}
\end{proposition}
The first part, \eqref{eq:combi:1}, was already observed in \cites{Sollich03,Sollich99}. The second statement, \eqref{eq:combi:2}, is not very far from being tight. Indeed, \cite{Chung01} proved that $c^n\le |V(n)|/(2^{\binom{n}{2}}n!)\le C^n$ with $c\approx 0.67$ the largest root of $384x^3-336x^2+54x-1$ and $C=1/\log 4\approx 0.72$. Proving Proposition \ref{prop:comb} is an excellent exercise, which we invite the reader to do before moving on. We provide a full proof, as it is very instructive of the mechanism governing the relaxation of the East model (see Figure~\ref{fig:East_path} for an illustration).

\begin{proof}[Proposition~\ref{prop:comb}]
We prove \eqref{eq:combi:1} by induction on $n$. The statement is trivial for $n=0$. Fix $n\ge 1$ and assume that for all $m< n$ we have $\ell(m)= 2^m-1$. Given $\omega\in V(n)\setminus\{\bone_\Lambda\}$, let $k\le n$, $x_k<\dots<x_1<x_0=0$ and $X=\{x_1,\dots,x_k\}$ be such that $\omega=\bzero_X\cdot\bone_{\Lambda\setminus X}$.

\begin{claim}
There exists $i\in\{1,\dots,k\}$ such that $x_{i}-x_{i-1}\le 2^{n-k}$.
\end{claim}
\begin{proof}
Assume the contrary. It is important to recall that the inverse of a legal path is legal, so there exists a legal path from $\omega$ to $\bone$ via configurations with at most $n$ empty sites. If $k=n$, this immediately gives that there exists $i\in\{1,\dots,k\}$ with $x_i=x_{i-1}-1$, since otherwise no legal move can remove empty sites and there are already $n$ of them. 

Assume $k<n$ and let $\omega^{(j)}$ be a legal path from $\omega$. We prove by a further induction on $j\ge 0$ that $\omega^{(j)}_X=\bzero_X$. The statement is trivial for $j=0$. If it is true for all $t<j\ge1$, then we can decompose the dynamics into the intervals $\{x_{i}+1,\dots,x_{i-1}-1\}$ for $i\in\{1,\dots,k\}$. Each such interval starts with $\bone$ initial condition and has $\bzero$ boundary condition. Thus, by the first induction hypothesis for $m=n-k$, we necessarily have $\omega_{x_i+1}^{(j-1)}=1$. Recalling that $\cU=\{\{1\}\}$ and \eqref{eq:def:cx}, this implies that $\omega_{x_i}^{(j)}=\omega_{x_i}^{(j-1)}=0$, completing the second induction. Since there is a legal path from $\omega$ to $\bone$, but the empty sites of $\omega$ cannot be removed, we obtain the desired contradiction proving the claim.
\end{proof}

Returning to the first induction, let $i\in\{1,\dots,k\}$ such that $x_i-x_{i-1}\le 2^{n-k}$. By the induction hypothesis for $m=n-k$, we can find a legal path $\gamma$ in $\Lambda_i=\{x_i+1,\dots,x_{i-1}-1\}$ with boundary condition $\omega_{\Lambda\setminus\Lambda_i}\cdot\bzero_{\bbZ\setminus\Lambda}$ from $\omega_{\Lambda_i}=\bone_{\Lambda_i}$ to a configuration $\omega'\in\Omega_{\Lambda_i}$ with $\omega'_{x_i+1}=0$ and $\gamma$ features at most $n-k$ empty sites in $\Lambda_i$ simultaneously. To see this, consider a legal path placing an empty site at $-l(n-k)\le x_i+1-x_{i-1}$, truncate it at the first step when an empty site is placed at $x_i+1-x_{i-1}$ and shift this path by $x_{i-1}$. 

We can now form the path from $\omega$ to $\omega^{x_i}$ by performing $\gamma$, then occupying $x_i$ and then performing the inverse of $\gamma$ (which is still legal, because only the boundary condition at $x_{i-1}$ is used due to the orientation of the East update family). By construction this path never creates more than $k+(n-k)=n$ empty sites in $\Lambda$, so $\omega^{x_i}\in V(n,k-1)$. Let $\omega^{x_i}=\bzero_{X'}\cdot\bone_{\Lambda\setminus X'}$ with $X'=\{x'_1,\dots,x'_{k-1}\}$ and $0>x'_1>\dots>x'_{k-1}$. Then
\[-x_k=\sum_{j=1}^k(x_{j-1}-x_j)\le 2^{n-k}+\sum_{j=1}^{k-1}(x'_{j-1}-x_j')=2^{n-k}-x'_{k-1}.\]
Iterating this inequality, we obtain
\[-x_k\le \sum_{j=n-k}^{n}2^j\le \sum_{j=0}^{n-1}2^j=2^{n}-1,\]
so $l(n)\le 2^n-1$ as desired. To see that this is an equality it suffices to follow the equalities above, which naturally leads to the path depicted in Figure~\ref{fig:East_path}.

Finally, \eqref{eq:combi:2} also follows easily from the above. Namely, each configuration $\omega\in V(n,k)$ can be encoded by the index $i\in\{1,\dots,k\}$, the distance $x_{i-1}-x_i\in\{1,\dots,2^{n-k}\}$ and the configuration $\omega^{x_i}\in V(n,k-1)$. Iterating this encoding gives
\[|V(n,n)|\le \prod_{k=1}^{n}(k2^{n-k})=n!2^{\binom n2},\]
which concludes the proof.
\end{proof}

\subsubsection{From the combinatorial result to the emptying time}
\label{subsubsec:bottleneck:deduction}
We now deduce the lower bound of Theorem~\ref{th:East} from Proposition~\ref{prop:comb}. This was done in somewhat different ways in \cites{Aldous02,Cancrini10}, but we rather present a proof along the lines of \cite{Hartarsky22univlower}, which is more adapted to generalisations. Let $n=\lfloor (\log(1/q)-\log\log(1/q))/\log2\rfloor$ and $\Lambda_n=\{0,\dots,2^n-1\}$. In view of Proposition~\ref{prop:comb}, we identify configurations in $V(n)$ and $V(n,n)$ with their restriction to $\Lambda_n$. Let $\cA=\{\omega\in\Omega:\omega_{\Lambda_n}=\bone_{\Lambda_n}\}$. By \eqref{eq:combi:1}, if $\cA$ occurs at time $0$, then there exists $t\le \tau_0$ such that $\omega_{\Lambda_n}(t)\not\in V(n)$. But, exiting $V(n)\times\Omega_{\bbZ\setminus\Lambda_n}$ requires visiting $V(n,n)\times\Omega_{\bbZ\setminus\Lambda_n}$. 

Let $T=(n^2q)^{-n/2}$ and let $N$ denote the number of updates (legal or illegal) at sites $x\in\Lambda_n$ up to time $T$. Then $N$ has the Poisson law with parameter $nT$. For $i\in\{1,\dots,N\}$, let $\theta_i$ denote the times of these updates. Observe that, by stationarity, $\omega(\theta_{{i}})$ is distributed according to $\mu$ (see \cite{Hartarsky22univlower}*{Claim~3.11} for a formal proof). Putting this together, we get
\begin{align}
\nonumber\bbP_\mu(\tau_0\le T)&{{}\le 1-\mu(\cA)+\bbP_\mu\left(\bigcup_{i=1}^N\left\{\omega_{\Lambda_n}(\theta_i)\in V(n,n)\right\}\right)}\\
&{}\le 1-\mu(\cA)+\bbP(N\ge 2nT)+2nT\mu_{\Lambda_n}(V(n,n)),\label{eq:East:union:bound}
\end{align}
that is, if $\tau_0\le T$, then we start outside $\cA$, or there are many updates, or at some update the configuration is in $V(n,n)$.

It remains to bound \eqref{eq:East:union:bound}. Firstly, for $q\to0$ we have
\[\mu(\cA)=(1-q)^{2^n}\ge (1-q)^{1/(q\log(1/q))}\to 1.\]
Furthermore, by the Bienaym\'e-Chebyshev inequality, $\bbP(N\ge 2nT)\to 0$, as $nT\to \infty$, which is the case when $q\to 0$. Finally, by \eqref{eq:combi:2}, \begin{align*}
2nT\mu_{\Lambda_n}(V(n,n))&{}= 2nT|V(n,n)|q^n(1-q)^{2^n-1-n}\le 2nTn!2^{\binom n2}q^n\\
&{}\le 2n\left(n^2q\right)^{-n/2} en(n/e)^n2^{n^2/2}q^n\le 2en^2e^{-n}\to 0.
\end{align*}
Inserting these bound in \eqref{eq:East:union:bound}, we obtain that $\bbP_\mu(\tau_0>T)\to 1$ as $q\to 0$. This concludes the proof for the emptyting time, since
\[\frac{\log T}{(\log(1/q))^2}\to\frac{1}{2\log 2}.\]
The analogous lower bound for $\trel$ follows directly from \eqref{eq:F0F1}.

\subsection{Upper bound: the bisection technique}
\label{subsec:East:upper}
The upper bound of Theorem~\ref{th:East} is our first encounter with the bisection technique. It was introduced in \cite{Cancrini08}, drawing inspiration from \cite{Martinelli99}*{Proposition 3.5}. It is not only very useful for the study of KCM, but has also been applied in other settings \cites{Caputo12,Bhatnagar07}.

\subsubsection{Two-block dynamics}The reader may have noticed that up to now we have not really proved any Poincar\'e inequality. We have only been reducing one inequality to another one we already know via renormalisation and canonical paths. The next lemma is, in a sense, the only Poincar\'e inequality we prove from scratch in this monograph. Morally, it deals with the East model on only two sites with empty boundary condition. This being a Markov process with only 4 states, which even happens to be a birth-death chain, one could compute the spectrum of its generator explicitly by hand. However, we state the result directly for the generalised version of the East model, as this is the version that is useful in renormalisation arguments. While the result is originally from \cite{Cancrini08}*{Proposition~4.5}, we rather give the proof from \cite{Hartarsky22phd}*{Lemma~1.3.8}, which is more probabilistic.
\begin{lemma}[Two-block dynamics]
\label{lem:two-block}
Let $(\bbX,\pi)$ be the product of two finite probability spaces $(\bbX_1,\pi_1)$ and $(\bbX_2,\pi_2)$. Let $\var_{1}(f)=\var_\pi(f(X_1,X_2)|X_2)$ and similarly for $\var_{2}(f)$. Fix a nonempty event $\cX\subset \bbX_1$. Then for any $f:\bbX\to\bbR$
\[\var_{\pi}(f)\le \frac{\bbE_\pi(\var_{1}(f)+\1_{\cX}\var_{2}(f))}{1-\sqrt{1-\pi(\cX)}}\le \frac{2}{\pi(\cX)}\bbE_\pi\left(\var_{1}(f)+\1_{\cX}\var_{2}(f)\right).\]
\end{lemma}
A way to interpret this is as a Poincar\'e inequality (i.e.\ bound on the relaxation time) for a continuous time Markov chain which updates $X_1$ at rate $1$ and updates $X_2$ at rate $1$, provided that $\cX$ occurs. In fact, the relaxation time bound $1/(1-\sqrt{1-\pi(\cX)})$ is optimal.
\begin{proof}
Couple two copies of the chain described above, by attempting the same updates in both. For this, use a graphical representation as in Section~\ref{subsec:markov} attempting updates at $X_1$ and $X_2$ with rate $1$, but deeming those in $X_2$ illegal if $\cX$ does not occur. The two chains clearly coalesce as soon as we update $X_1$ so that $\cX$ occurs and then immediately update $X_2$. Consider (legal or illegal) updates on $X_2$ preceded by an update at $X_1$. Their number up to time $T$ is $\lfloor N/2\rfloor$ with $N$ a Poisson random variable with mean $T$. Each one succeeds in coupling the chains independently with probability $\pi(\cX)$. It is elementary to check that $\bbE(\lambda^N)=e^{-T(1-\lambda)}$ for any $\lambda\in(0,\infty)$. Thus, the probability that the two chains are not equal at time $T$ is at most
\[\bbE\left[(1-\pi(\cX))^{\lfloor N/2\rfloor}\right]\le \frac{1}{(1-\pi(\cX))\exp(T(1-\sqrt{1-\pi(\cX)}))}.\]
Classical results on Markov chains \cite{Levin09}*{Proposition~4.7, Corollary~12.6, Remark~13.13}\footnote{For continuous time Markov chains the spectral radius in \cite{Levin09}*{Corollary~12.6} is replaced by $e^{-1/\trel}$.} then give that $\trel\le 1/(1-\sqrt{1-\pi(\cX)})$, as desired. Finally, for the second inequality we use the Bernoulli inequality: $x\le 2(1-\sqrt{1-x})$ for any $x\in[0,1]$.
\end{proof}

If we apply Lemma \ref{lem:two-block} to $\bbX_1=\bbX_2=\{0,1\}$, $\pi_1=\pi_2=\mu_q$ and $\cX=\{0\}$, we obtain that the relaxation time of the East KCM with boundary condition $\bzero_{\bbZ\setminus\{0,1\}}$ is at most $1/(1-\sqrt{1-q})$. When $q\to0$ this is approximately $2/q$. However, the great advantage of the first bound in Lemma~\ref{lem:two-block} is that its prefactor tends to $1$ as $\pi(\cX)\to 1$. Thus, we can hope to apply this result to progressively larger volumes and more likely events $\cX$ and obtain a relaxation time bound uniform in the volume. This is the main idea behind the upper bound in Theorem~\ref{th:East}, which we discuss next.

\subsubsection{Bisection technique}
\label{subsubsec:bisection:1d}
Recall that by Proposition~\ref{prop:infinite:to:finite}, in order to prove the upper bound of Theorem~\ref{th:East}, it suffices to bound the finite volume relaxation times $\trel^\Lambda$ uniformly as the volume diverges. The bisection technique consists in an iterative application of Lemma~\ref{lem:two-block}. Rather than presenting the somewhat technical and artificial-looking proof of \cite{Cancrini08}, let us take a more instructive approach to see how the proof is conceived. 

\paragraph{Basic idea}
Let $\Lambda_k=\{1,\dots,2^k\}$ and $\Lambda'_k=\Lambda_{k+1}\setminus\Lambda_k$ for any $k\ge 0$. Clearly, $\trel^{\Lambda_0}=1$. We further seek to relate $\trel^{\Lambda_k}$ and $\trel^{\Lambda_{k+1}}$. Fix $k\ge 0$ and apply Lemma~\ref{lem:two-block} with $\bbX_1=\Omega_{\Lambda'_k}$, $\bbX_2=\Omega_{\Lambda_{k}}$, $\pi_1=\mu_{\Lambda'_k}$, $\pi_2=\mu_{\Lambda_k}$ and $\cX=\{\omega_{2^k+1}=0\}\subset\bbX_1$. This gives
\begin{equation}
\label{eq:East:bisection:1}\var(f)\le \frac{\mu(\var_{\Lambda'_k}(f)+\1_\cX\var_{\Lambda_k}(f))}{1-\sqrt{\varepsilon_k}}\end{equation}
for any $f:\Omega_{\Lambda_{k+1}}\to\bbR$, where $\varepsilon_k=1-\pi_1(\cX)=1-q$. For the first term above, \eqref{eq:poincare} and translation invariance directly give
\begin{equation}
\label{eq:East:bisection:2}\var_{\Lambda'_k}(f)\le \trel^{\Lambda_k}\cD_{\bzero_{\bbZ\setminus\Lambda'_k}}(f).\end{equation}
Yet, the fact that the East update family $\cU=\{\{1\}\}$ only looks to the right gives that for any $x\in\Lambda'_k$ we have $c_x^{\bzero_{\bbZ\setminus\Lambda'_k}}=c_x^{\bzero_{\bbZ\setminus\Lambda_{k+1}}}$, so \begin{equation}
\label{eq:East:bisection:3}\cD_{\bzero_{\bbZ\setminus\Lambda'_k}}(f)=\sum_{x\in\Lambda'_k}\mu\left(c_x^{\bzero_{\bbZ\setminus\Lambda_{k+1}}}\var_x(f)\right)\end{equation}
by \eqref{eq:dirichlet}. For the second term in \eqref{eq:East:bisection:1}, we similarly have
\begin{equation}
\label{eq:East:bisection:4}
\1_\cX\cD_{\bzero_{\bbZ\setminus\Lambda_k}}(f)\le \sum_{x\in\Lambda_k}\mu\left(c_x^{\bzero_{\bbZ\setminus\Lambda_{k+1}}}(f)\right),\end{equation}
since $\cX$ guarantees precisely the presence of an empty site as in the boundary condition $\bzero_{\bbZ\setminus\Lambda_k}$. Combining \eqref{eq:East:bisection:1}-\eqref{eq:East:bisection:4} and \eqref{eq:def:Trel}, we obtain
\begin{equation}
\label{eq:East:bisection:recurrence}\trel^{\Lambda_{k+1}}\le \frac{\trel^{\Lambda_k}}{1-\sqrt{\varepsilon_k}}.
\end{equation}
Iterating the above relation, and recalling Proposition~\ref{prop:infinite:to:finite}, we get
\begin{equation}
\label{eq:East:product}\trel\le \frac2q\prod_{k=0}^\infty\frac{1}{1-\sqrt{\varepsilon_k}}.\end{equation}

\paragraph{Spread the boundary condition}
Unfortunately, this first attempt fails, because $\varepsilon_k=1-q$ does not decay with $k$, so the product in \eqref{eq:East:product} is infinite. In order to fix this problem, we should define the event $\cX$ differently. Namely, fix some $\delta>0$ small enough and let $\delta_k=\lfloor 2^{k(1-\delta)}\rfloor$. Then let $\cX=\{\omega_{2^k+\{1,\dots,\delta_k\}}\neq \bone\}$, that is, there is an empty site among the first $\delta_k$ sites to the East of $\Lambda_k$. With this choice $\varepsilon_k=(1-q)^{\delta_k}$ does decay sufficiently fast for the right hand side of \eqref{eq:East:product} to be finite. In fact, one can compute that the product is at most $q^{-C}(1/q)^{\log_2(1/q)/(2-2\delta)}$ for some constant $C>0$ (see \cite{Cancrini08}*{Section 6.1} for more details). Taking $\delta$ small, this is exactly the upper bound we want.

However, with this choice of $\cX$, \eqref{eq:East:bisection:4} is no longer valid. To deal with this issue, on $\cX$ we can define the random variable
\begin{equation}
\label{eq:def:xi}\xi(\omega)=\left\{\max i\le \delta_k:\omega_{2^k+i}=0\right\}\ge 1\end{equation}
indicating the position of the rightmost empty site in $\Lambda'_k$ at distance at most $\delta_k$ from the boundary of $\Lambda_k$. Then we can rewrite
\begin{align}
\label{eq:East:muLk+1}\mu_{\Lambda_{k+1}}\left(\1_\cX\var_{\Lambda_k}(f)\right)&{}=\sum_{i=1}^{\delta_k}\mu_{\Lambda_{k+1}}\left(\1_{\xi=i}\var_{\Lambda_k}(f)\right)\\
&{}=\sum_{i=1}^{\delta_k}\mu_{\Lambda_{k+1}}\left(\1_{\xi=i}\mu_{\{1,\dots,2^k+i-1\}}\left(\var_{\Lambda_k}(f)\right)\right),\nonumber\end{align}
since $\1_{\xi=i}$ is independent of $\omega_{\{1,\dots,2^k+i-1\}}$. But, setting $V_i=\{1,\dots,2^k+i-1\}$, we have
\begin{equation}
\1_{\xi=i}\mu_{V_i}(\var_{\Lambda_k}(f))\le \1_{\xi=i}\var_{V_i}(f)\le \trel^{V_i}\1_{\xi=i}\cD_{\bzero_{\bbZ\setminus V_i}}(f),
\label{eq:East:1xii}
\end{equation}
since $\xi=i$ guarantees that $\omega_{2^k+i}=0$, so the boundary condition is indeed empty. Note that in the first inequality above, we used the convexity of the variance 
 that implies that for any volumes $A,B$ it holds
 \begin{equation}\label{eq:convexity:var}
 \mu_{A}(\var_{B}(f))\leq\var_{A\cup B}(f).
 \end{equation}

Combining \eqref{eq:East:muLk+1} and \eqref{eq:East:1xii} together with the monotonicity property \eqref{eq:monotonicity}, we get
\begin{align*}\mu_{\Lambda_{k+1}}\left(\1_\cX\var_{\Lambda_k}(f)\right)&{}\le {\trel^{V_{\delta_k}}\sum_{i=1}^{\delta_k}\mu_{\Lambda_{k+1}}\left(\1_{\xi=i}\cD_{\bzero_{\bbZ\setminus V_i}}(f)\right)}\\
&{}\le \trel^{V_{\delta_k}}\sum_{x\in V_{\delta_k}}\mu_{{\Lambda_{k+1}}}\left(c_x\var_x(f)\right),\end{align*}
using $\sum_{i=1}^{\delta_k}\1_{\xi=i}\le 1$ and the fact that $\1_{\xi=i}\cdot c_x^{\bzero_{\bbZ\setminus V_i}}\le \1_{\omega_{2^{k}+i}=0}\cdot c_x^{\bzero_{\bbZ\setminus V_i}}\le c_x$ for any $i\in\{1,\dots,\delta_k\}$ and $x\in V_i$.
 Further combining this with \eqref{eq:East:bisection:1}-\eqref{eq:East:bisection:3}, we obtain
 \begin{equation}
 \label{eq:East:bisection:almost}\var(f)\le \frac{\trel^{V_{\delta_k}}}{1-\sqrt{\varepsilon_k}}\mu\left(\cD_{\Lambda_{k+1}}(f)+\sum_{x\in V_{\delta_k}\cap\Lambda'_k}\left(c_x^{\bzero_{\bbZ\setminus\Lambda_{k+1}}}\var_x(f)\right)\right).\end{equation}

Recalling \eqref{eq:def:Trel}, we see that \eqref{eq:East:bisection:almost} is almost the result we seek to prove, \eqref{eq:East:bisection:recurrence}. 

\paragraph{Final adjustments}
We are only left with mending two technical problems with the previous argument. Firstly, in the right hand side of \eqref{eq:East:bisection:almost}, the terms corresponding to $x\in V_{\delta_k}\cap\Lambda'_k$ appear twice (once in the Dirichlet form). In order to solve this, we consider many possible choices of the partition of $\Lambda_{k+1}$ into $\Lambda_k$ and $\Lambda'_k$, keeping the total volume $2^{k+1}$ and the overlap $\delta_k$ fixed. We then average \eqref{eq:East:bisection:almost} over these choices. This yields an additional factor $1+1/s_k$ in \eqref{eq:East:bisection:recurrence}, where $s_k=\lfloor2^{k\varepsilon/3}\rfloor$ is the number of choices we consider.

Secondly, $V_{\delta_k}$ is slightly larger than $\Lambda_k$ and so is the corresponding relaxation time. This issue is solved by choosing the sizes of all $\Lambda_k$ growing as $2^k+2^{k(1-\varepsilon/3)}$ rather than $2^k$. Once these problems are solved, we get \eqref{eq:East:bisection:recurrence} (with the additional factor $1+1/s_k$) and conclude the proof of the upper bound in Theorem~\ref{th:East}. The upper bound on $\tau_0$ follows from \eqref{eq:F0F1}.

Interestingly, prior to \cite{Cancrini08}, the conjecture in the physics literature on the exponent was $T_{\rm rel} \sim q^{{\log_2(1/q)}}$, with an exponent off by a factor $2$. In order to get the correct scaling one has to take into account a subtle balance between the energetic and entropic contributions that, atypically, lie at the same level for the one-dimensional East model. Remarkably, the bisection technique is able to automatically  take into account this subtle balance and provide a tight result correcting the conjectured exponent.  

\section{FA-2f}
\label{sec:FA2f:1d}
In this section, we briefly discuss the one-dimensional FA-2f update family $\cU=\{\{-1,1\}\}$. From Theorem~\ref{th:BP}, we have $\qc=1$, because two neighbouring occupied sites remain occupied at all times. Normally, our study of the model would end here, because the phase $q<\qc$ is rather complicated, but in the one-dimensional setting, we are able to say more. 

Specifically for FA-2f, one can check that the BP transformation (recall \eqref{eq:def:BP}) satisfies $[\omega]=\sB_\cU(\omega)$ for any $\omega\in\Omega$.That is, the BP process becomes stationary after one step. Taking Corollary~\ref{cor:legal:path} into account, the KCM dynamics can be decomposed into independent dynamics on intervals delimited by two occupied sites. On each such interval, we recover what is known as the hard-core Glauber dynamics with fugacity $\lambda=(1-q)/q$. That is because occupied sites cannot appear next to other occupied sites, while emptying is always possible (within an interval delimited by two occupied sites). There is a rich literature on Glauber dynamics of the hard-core model, particularly on general graphs with bounded degree, but also lattices of higher dimension (see e.g.\ \cites{Levin09,Sly10}). However, the one-dimensional lattice is somewhat degenerate from the standard viewpoint and does not appear to have been the subject of much study.

One natural question one could ask is how the system behaves on its \emph{ergodic component}, that is, the set of configurations such that $[\omega]=\bzero$. One can then still study the $q\to 0$ regime, which can be viewed as quenching the model from inverse temperature $\beta=-\infty$ to $\beta$ large. In this setting, it is possible to prove that the relaxation time (this time with respect to the Gibbs measure of the hard-core model, see \cite{Georgii11} for background, rather than the plain $\mu_q$ product measure) is finite for any $q>0$. This can be obtained, for example using classical techniques such as block dynamics and strong spatial mixing \cites{Martinelli94,Martinelli94a,Martinelli94b}. However, as we will see in Section~\ref{sec:general:KCM}, this can also be achieved via the bisection technique, which applies more broadly.

\section{General KCM}
\label{sec:general:KCM}
We have so far seen two ways in which it may be desirable to generalise KCM. Namely, allowing a state space larger than $\{0,1\}$ for each site, and working with the dynamics restricted to an ergodic component. Furthermore, for the purposes of studying higher-dimensional models, it is also useful to consider inhomogeneous KCM with site-dependent update families. We next define general KCM incorporating all these features, following \cite{Hartarsky21b}.

Fix $R>0$, $q\in(0,1)$ and $\Lambda\Subset\bbZ$. For each $x\in\Lambda$, fix a probability space $(\Omega_x,\pi_x)$ with $|\Omega_x|<\infty$ and an event $\cI_x\subset\Omega_x$ with $\pi_x(\cI_x)\ge q$. Let $(\Omega,\pi)=(\prod_{x\in\Lambda}\Omega_x,\bigotimes_{x\in\Lambda}\pi_x)$ be the corresponding product space. A boundary condition is any configuration $\eta\in\{0,1\}^{\bbZ\setminus\Lambda}$. Further fix an update family $\cU_x$ for each $x\in\Lambda$ so that for any $U\in\cU_x$ and $u\in U$ we have $|u|\le R$. The constraint at $x\in\Lambda$ is defined by
\[c_x^\eta(\omega)=\max_{U\in\cU_x}\prod_{u\in U,x+u\in\Lambda}\1_{\omega_{x+u}\in\cI_{x+u}}\prod_{u\in U,x+u\in\bbZ\setminus\Lambda}\left(1-\eta_{x+u}\right)\]
for a configuration $\omega\in\Omega$ and a boundary condition $\eta\in\{0,1\}^{\bbZ\setminus\Lambda}$.

Consider the Markov process such that for each site $x\in\Lambda$ such that $c_x=1$, the state of site $x$ is updated to an independent random variable with law $\pi_x$. That is, the process with generator
\[\cL(f)=\sum_{x\in\Lambda}c_x^\eta\cdot(\pi_x(f)-f).\]
Let $\cC$ be an arbitrarily chosen irreducible component of $\Omega$ for this dynamics and $\mu=\pi(\cdot|\cC)$. Then define $\trel$ via  \eqref{eq:dirichlet} and \eqref{eq:def:Trel}. We refer to this Markov process restricted to $\cC$ as a \emph{general KCM with range $R$ and facilitating parameter $q$}. The following result was proved by Hartarsky \cite{Hartarsky21b} via the bisection technique adapted for going back and forth several times between the two blocks in the two-block Lemma~\ref{lem:two-block}.

\begin{theorem}[General KCM upper bound]
\label{th:general:1d}
    There exists $C>0$ depending only on $R$ such that for every $q\in(0,1)$,
    \[\trel\le (2/q)^{C\log\min(|\Lambda|,2/q)}.\]
\end{theorem}

In words, Theorem~\ref{th:general:1d} states that, for any one-dimensional general KCM with uniformly bounded update rule range and probability of the facilitating state uniformly bounded away from 0 has a finite relaxation time scaling at most like the one of the East model (recall Theorem~\ref{th:East}). Note that the minimum reflects the fact that the product in \eqref{eq:East:product} approaches its limiting value for scales $k\approx\log(1/q)$. Theorem~\ref{th:general:1d} can also be extended to infinite volume along the lines of Proposition~\ref{prop:infinite:to:finite}, but one needs to be careful in defining the irreducible components (see \cite{Hartarsky21b}*{Observation 3}).

Let us discuss a few useful applications of Theorem~\ref{th:general:1d}. Firstly, FA-2f on its ergodic component is covered, just like any 1-dimensional (homogeneous binary) KCM. More importantly, we have the following bounds for the generalised FA-1f and East KCM. They were both derived in \cite{Martinelli19a}*{Proposition~3.4}, using the methods discussed in Sections~\ref{sec:FA1f} and~\ref{subsec:East:upper}, the second one also following from Theorem~\ref{th:general:1d}.

\begin{proposition}[Generalised FA-1f and East upper bound]
    \label{prop:FA1f:East:generalised}
    Let $\Lambda$ be a segment and $\eta\in\{0,1\}^{\bbZ\setminus\Lambda}$ with $\eta_x=0$ if $x>\max\Lambda$ and $\eta_x=1$ if $x<\min\Lambda$. Consider a general KCM on $\Lambda$ and assume it to be homogeneous with $\cU_x=\cU$ for all $x\in\Lambda$, with range 1, facilitating parameter $q\in(0,1)$ and boundary condition $\eta$. Then for some absolute $C>0$ we have
    \[\trel\le \begin{cases}
        (2/q)^{C}&\cU=\{\{-1\},\{1\}\},\\
        (2/q)^{C\log\min(|\Lambda|,2/q)}&\cU=\{\{1\}\}.
    \end{cases}\]
\end{proposition}
Let us note that one can also prove a polynomial bound on the relaxation time of (generalised) FA-1f on a segment with $\bone$ boundary condition on its ergodic component (that is, all configurations except $\bone$), see \cite{Blondel13}.

\section{Conclusion}
\label{subsec:tools:1d}
Let us review the state of our toolbox after the developments of this chapter (recall Section~\ref{subsec:tools:BP}). 

\noindent\textbf{Test functions.} The lower bound of Theorem~\ref{th:FA1f} was proved by guessing a non-trivial test function. This technique will not take us any further in the sequel for more sophisticated models, because guessing a suitable function and being able to compute the variance and Dirichlet form are quite implausible.

\noindent\textbf{Canonical paths.} The upper bound of Theorem~\ref{th:FA1f} relied on more subtle canonical paths than the ones provided by BP in Section~\ref{sec:legal:paths}. They reflect the heuristic view of the dynamics of FA-1f. Once we have a good intuition about the dominant relaxation mechanism of a KCM, we could, in principle try to implement it in a canonical path. Unfortunately, this approach quickly goes out of hand as the models get more complicated, since explicitly defining and analysing the paths involved becomes very laborious and quite tricky. We therefore avoid further recourse to canonical paths.

\noindent\textbf{Renormalisation.} In the proof of Theorem~\ref{th:FA1f}, we saw the details of the 1-dimensional renormalisation we already saw in Section~\ref{sec:exponential:decay}. This technique for proving upper bounds will be developed much further in the next chapter. Now that we have some simple KCM to build on, it will become our bread and butter tool for proving upper bounds.

\noindent\textbf{Combinatorial bottlenecks.} This method discussed in Section~\ref{subsec:East:lower} will be our method of choice for proving lower bounds on time scales in what follows. The content of Section~\ref{subsubsec:bottleneck:deduction} will require essentially no adaptation. The main difficulty in implementing this approach lies in identifying what needs to happen before the origin can be updated and proving that it is indeed necessary. Finding the correct bottleneck is usually guided by heuristics of the dominant relaxation mechanism (from upper bounds), estimating the probability of the bottleneck will usually be done using ideas from BP, while entropy tends not to pose problems. Thus, the main issue is proving rough analogues of \eqref{eq:combi:1} in more advanced settings.

\noindent\textbf{Bisection.} The idea used in Section~\ref{subsec:East:upper} was to iterate the simple two-block dynamics of Lemma~\ref{lem:two-block}. Bisection is our primary technique for proving directly that a KCM has finite relaxation time in infinite volume. In the next chapter we will see how to do the same in higher dimensions.

\noindent\textbf{General KCM.} Thanks to bisection, we were able to treat one-dimensional KCM in great generality. They are ready to use in renormalisation schemes. Although we will need to introduce some higher-dimensional models with general state space, one-dimensional general KCM will be sufficient for most of our purposes.

	\chapter{Fredrickson-Andersen 2-spin facilitated model}
\label{chap:FA2f}
  
\abstract{This chapter uses the setting of the FA-2f model to develop several new tools for determining the emptying time of the origin with high precision as $q\to 0$. We begin by using bisection in higher dimensions to show that this time scale is finite. We then discuss a robust long range Poincar\'e inequality approach for proving upper bounds. Finally, we assess the sharp threshold of FA-2f, which relies on a robust relation with bootstrap percolation for the lower bound and on the very flexible method of matryoshka dolls for the upper bound. All of these methods generalise to treat various other models.\\}

Recall from Section~\ref{sec:rules} that the FA-2f model's constraint requires at least two empty neighbours in order to change the state of a site. Throughout this section we work in two dimensions for simplicity of notation, but the arguments apply equally well in any dimension. For FA-2f (in two dimensions), the natural geometry is rectangular. It is therefore, convenient to denote by
\begin{equation}
\label{eq:def:rectangle}
R(a,b)=\{0,\dots,a-1\}\times\{0,\dots,b-1\}\subset\bbZ^2\end{equation}
the rectangle of side lengths $a,b\in\bbN$.

\section{Bisection in higher dimensions}
\label{sec:bisection:2d}
Let us begin by paying our debt by proving that $\trel<\infty$ for any $q>0$ for FA-2f. Recall that in the proof of Theorem~\ref{th:exponential}, the only fact whose proof was postponed until now is that for $q$ close enough to 1, the KCM with update family 
\[\cU'=\{U'\}=\left\{\{(1,0),(1,1),(0,1)\}\right\}\]
has a finite relaxation time. In fact, we rather used this result for the version of this KCM with general state space, as in Section~\ref{sec:general:KCM}, but the proof is very similar, so we focus on the binary case. In Section~\ref{subsec:East:upper} we discussed how to prove this in one dimension via the bisection technique. We next explain how the method adapts to higher dimensions, following \cite{Cancrini08}. In fact, the proof works for any $q>\qc(\cU')$.

We start with a simple observation, which could be attributed to Schonmann~\cite{Schonmann90}.
\begin{observation}[Oriented percolation correspondence]
\label{obs:GOSP}
    Endow $\bbZ^2$ with the oriented graph structure defined by the edge set $E=\{(\bx,\bx+\bu):\bx\in\bbZ^2,\bu\in U'\}$. Then for any configuration $\omega\in\Omega$, in $\cU'$-BP, the emptying time $\tbp$ is given by the number of sites in the longest (oriented) path from $0$ whose sites are all occupied.
\end{observation}
The proof of this fact by induction on the number of iterations of the BP map $\sB_\cU$ from \eqref{eq:def:BP} is left as an exercise to the reader. In view of Observation~\ref{obs:GOSP}, for $\bx\in A\subset\bbZ^2$, $B\subset\bbZ^2$ and $\omega\in\Omega_A$, we write $\bx\to[A]B$ if there exists a sequence of sites $(\bx_i)_{i=0}^{l}\in A^{l+1}$ occupied in $\omega$ with $\bx_0=\bx$, $\bx_l\in B$ and $(\bx_{i-1},\bx_i)\in E$ for all $i\in\{1,\dots,l\}$.

\begin{figure}
\centering
\begin{tikzpicture}[x=0.4cm,y=0.4cm]
\draw (0,0)--(20,0)--(20,10)--(0,10)--cycle;
\draw (10,0)--(10,10);
\draw (15,0)--(15,10);
\foreach \i in {0,1,2,3,4,5,7,8}
    \draw (10.5,0.5+\i) circle (0.25);
\foreach \i in {1,2,3,5,8}
    \fill (11.5,0.5+\i) circle (0.25);
\foreach \i in {0,4,6,7}
    \draw (11.5,0.5+\i) circle (0.25);
\foreach \i in {0,4,6,7,8}
    \fill (12.5,0.5+\i) circle (0.25);
\foreach \i in {3,5}
    \draw (12.5,0.5+\i) circle (0.25);
\foreach \i in {0,4,6,8}
    \fill (13.5,0.5+\i) circle (0.25);
\foreach \i in {2,3,5,7,9}
    \draw (13.5,0.5+\i) circle (0.25);
\foreach \i in {0,3,4,7,9}
    \fill (14.5,0.5+\i) circle (0.25);
\foreach \i in {1,2,5,6,8}
    \draw (14.5,0.5+\i) circle (0.25);
\draw (14.5,9.5)--(12.5,7.5);
\draw (13.5,8.5)--(11.5,8.5);
\draw (12.5,8.5)--(12.5,6.5);
\draw (14.5,7.5)--(13.5,6.5)--(12.5,6.5)--(11.5,5.5);
\draw (14.5,4.5)--(12.5,4.5)--(11.5,3.5)--(11.5,1.5);
\draw (14.5,4.5)--(14.5,3.5);
\draw (14.5,0.5)--(12.5,0.5);
\draw [decorate,decoration={brace,amplitude=10pt}]   (10,0) -- (0,0) node [midway,yshift= -0.5cm] {$\Lambda_k$};
\draw [decorate,decoration={brace,amplitude=10pt}]   (20,0) -- (10,0) node [midway,yshift= -0.5cm] {$\Lambda'_k$};
\draw [decorate,decoration={brace,amplitude=10pt}]   (0,10) -- (15,10) node [midway,yshift=0.5cm] {$V_k$};
\fill[opacity=0.25] (10,0)--(15,0)--(15,10)--(10,10);
\end{tikzpicture}
\caption{\label{fig:bisection:2d}Illustration of the bisection in Section~\ref{sec:bisection:2d}. The rectangle $\Lambda_{k+1}$ is partitioned into $\Lambda_k\sqcup\Lambda'_k$. The shaded rectangle $V_k'$ has width $\delta_k$. The full dots in it represent occupied sites in $\Gamma_k$, while $\partial\Gamma_k$ consists of the empty sites. In the figure the event $\cX_k$ occurs, since the paths do not reach the leftmost column of $V'_k$.}
\end{figure}
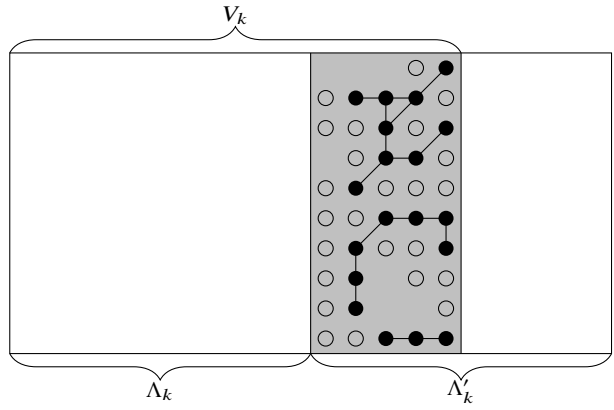

For any $k\ge 0$, consider the rectangles 
\begin{align*}
\Lambda_{2k}&{}=R\left(2^k,2^k\right)&\Lambda_{2k+1}&{}=R\left(2^{k+1},2^k\right)\end{align*}
and $\Lambda'_k=\Lambda_{k+1}\setminus\Lambda_k$, which is a translate of $\Lambda_k$. Note that the rectangles $(\Lambda_k)_{k\ge 0}$ are nested in such a way that each is obtained by stretching the previous rectangle twice either horizontally or vertically (see Figure~\ref{fig:bisection:2d}). These rectangles will play the role of the intervals $\Lambda_k$ in Section~\ref{subsubsec:bisection:1d}. As in one dimension, for any $k\ge 0$, we set $\delta_k=\lfloor 2^{k(1-\delta)}\rfloor$ with some fixed $\delta>0$ small enough. Further define
\begin{align*}
V'_{2k}&{}=\left(2^k,0\right)+R\left(\delta_k,2^k\right)&V'_{2k+1}&{}=\left(0,2^k\right)+R\left(2^{k+1},\delta_k\right)\\
\partial_+ V'_{2k}&{}=\left(2^{k}+\delta_k-1,0\right)+R\left(1,2^k\right)&\partial_+ V_{2k+1}&{}=\left(0,2^k+\delta_k-1\right)+R\left(2^{k+1},1\right)\\
\partial_- V'_{2k}&{}=\left(2^{k},0\right)+R\left(1,2^k\right)&\partial_- V_{2k+1}&{}=\left(0,2^k\right)+R\left(2^{k+1},1\right)\end{align*}
and $V_k=\Lambda_k\cup V'_k$ (see Figure~\ref{fig:bisection:2d}).

We next need to choose the facilitating event $\cX$ in Lemma~\ref{lem:two-block}, so that its probability gets close to 1 as the scale $k$ increases. Given $k\ge 0$, let 
\begin{align}
\label{eq:def:Gammak}
\cX_k&{}=\left\{\omega\in\Omega_{V'_k}:\Gamma_k\cap \partial_- V_k'=\varnothing\right\},&
\Gamma_k&{}=\left\{\bx\in V'_k:\bx\to[V'_k]\partial_+ V_k'\right\}.\end{align}
That is, $\cX_k$ is the event that no site on the boundary of $V'_k$ is connected to the other boundary via occupied sites in $V'_k$ (see Figure~\ref{fig:bisection:2d}).

With this definition it is a classical percolation result of Menshikov and Aizenman--Barsky \cites{Menshikov86,Aizenman87} that for some $c>0$ depending only on $q>\qc$, 
\begin{equation}
\label{eq:exp:decay:GOSP}
\varepsilon_k=1-\mu(\cX_k)\le \exp\left(-c\delta_k\right)
\end{equation}
(see \cite{Duminil-Copin16} for a simple proof, \cite{Grimmett99} for more background on percolation and \cite{Hartarsky22GOSP} for more background on oriented percolation). Indeed, taking Observation~\ref{obs:GOSP} into account, $1-\qc$ is the critical parameter of the oriented percolation model on $\bbZ^2$ with edge set $E$. Plugging \eqref{eq:exp:decay:GOSP} into \eqref{eq:East:product}, we easily obtain that the product there is finite.

Our next task is to show bound the second term in \eqref{eq:East:bisection:1} as
\begin{equation}
\label{eq:East:2d:dirichlet}\mu_{V_k}\left(\1_{\cX_k}\var_{\Lambda_k}(f)\right)\le \trel^{V_k}\sum_{\bx\in V_k}\mu_{V_{k}}\left(c_\bx^{\bzero_{\bbZ^2\setminus \Lambda_{k+1}}}\var_\bx(f)\right).\end{equation}
Once \eqref{eq:East:2d:dirichlet} is established, the proof is concluded by the same final adjustments as in Section~\ref{subsubsec:bisection:1d}. In order to prove \eqref{eq:East:2d:dirichlet}, we proceed similarly to the one-dimensional case, but we need to pay attention to the definition of the `rightmost empty site'. For $\gamma\subset V'_k$, let 
\begin{equation}
\label{eq:def:Gammakbar}\overline\gamma=\partial_+V'_k\cup\gamma\cup\left(V'_k\cap\left(\gamma-U'\right)\right).\end{equation}
Then $\partial\Gamma_k=\overline\Gamma_k\setminus\Gamma_k$ (see Figure~\ref{fig:bisection:2d}) will play the role of the rightmost empty site $\xi$ in \eqref{eq:def:xi}. By \eqref{eq:def:Gammak} and \eqref{eq:def:Gammakbar}, $\Gamma_k$ (and $\overline \Gamma_k$) is measurable with respect to $\omega_{\overline\Gamma_k}$ and, if $\cX_k$ occurs, then $\omega_{\partial\Gamma_k}=\bzero$. Finally, observe that, by construction, $\partial\Gamma_k$ is a cut-set separating $\partial_+V'_k$ from $\Lambda_k$. Therefore, recalling \eqref{eq:convexity:var}, this gives
\begin{align*}
\mu_{V_{k}}\left(\1_{\cX_k}\var_{\Lambda_k}(f)\right)&{}\le \sum_{\gamma\subset V'_k\setminus\partial_-V'_k}\mu_{V_{k}}\left(\1_{\Gamma_k=\gamma}\var_{\Lambda_k}(f)\right)\\
&{}\le \sum_{\gamma\subset V'_k\setminus\partial_-V'_k}\mu_{V_{k}}\left(\1_{\Gamma_k=\gamma}\1_{\omega_{\partial\gamma}=\bzero}\var_{V_k\setminus\overline\gamma}(f)\right)\\
&{}\le \sum_{\gamma\subset V'_k\setminus\partial_-V'_k}\trel^{V_k\setminus\overline\gamma}\sum_{\bx\in V_k\setminus\overline\gamma}\mu_{V_{k}}\left(\1_{\Gamma_k=\gamma}c_\bx^{\bzero_{\bbZ^2\setminus \Lambda_{k+1}}}\var_\bx(f)\right),
\end{align*}
where we used $\1_{\omega_{\partial\gamma}=\bzero}c_\bx(\omega)\le c_\bx(\omega)$ for all $\bx\in V_k\setminus\overline{\gamma}$ and $\omega\in\Omega$. Recalling the the monotonicity property \eqref{eq:monotonicity}, we recover \eqref{eq:East:2d:dirichlet} as desired.

Bisection is good for proving $\trel<\infty$. However, if we follow the proof of Theorem~\ref{th:exponential} and take into account Theorem~\ref{th:BP}, we obtain the following very poor bound. 
\begin{corollary}[{Basic FA-2f upper bound}]
\label{cor:FA2f:bound:trivial}
For FA-2f in $d$ dimensions there exists $C>0$ such that for all $q>0$,
\[\trel\le \exp^{\circ 2}\left(Cq^{-1/(d-1)}\right).\]
\end{corollary}
In \cite{Cancrini08} this was improved to $\exp(Cq^{-5})$ using less brutal canonical paths than the ones suggested in the proof of Theorem~\ref{th:exponential}. However, this is still very far from the truth. The next sections examine other techniques for obtaining more accurate bounds.

\section{Long range renormalisation}
\label{sec:LRP}
In order to improve the bound on the relaxation time of FA-2f provided by Corollary~\ref{cor:FA2f:bound:trivial}, we need a new technique by Martinelli and Toninelli \cite{Martinelli19}. The key idea is to transform our KCM, which has short range but unlikely constraints, into a model with long range but very likely constraints. These constraints require the occurrence of a certain \emph{ droplet} far away and of a \emph{good environment}, on which the droplet can move, connecting the droplet location to the origin.  We start by proving a \emph{long range Poincar\'e inequality}, which is the key building block of this technique, in Section~\ref{sec:constrainedpoincare}.  In Section \ref{sec:FA2fbetter} we discuss how to combine it with renormalisation to obtain results about our model of interest.

\subsection{A long range constrained Poincar\'e inequality }
\label{sec:constrainedpoincare} 
Fix a finite probability space $(\mathbb X_0,\nu_0)$ and let 
$(\bbX, \nu)=(\prod_{\bx\in\bbZ^d}\bbX_0,\bigotimes_{\bx\in\bbZ^d}\nu_0)$ be the corresponding product space.
For each $\bx\in \bbZ^d$ let  $\bbN^d_{\bx}:=\bx+(\{0,1,\dots\}^d\setminus\{0\})$. Then fix a finite set $\Delta_0\subset\bbN^d_0$. Let $\cA_0\subset\bbX$ be an event depending only on the restriction $\eta_{\Delta_0}\in\bbX_{\Delta_0}$ of the state $\eta\in\bbX$ to $\Delta_0$. Define $\cA_\bx$ by translating $\cA_0$ by $\bx\in\bbZ^d$. For $\bx\in\bbZ^d$, set
\begin{align*}
r_\bx&{}=\1_{\cA_\bx},&\epsilon&{}=1-\nu_0(\cA_0).
\end{align*}
Then following result gives a taste of the more general \cite{Martinelli19}*{Theorem 2}.
\begin{proposition}[Long range constrained Poincar\'e inequality]
\label{teo:constrainedpoincare}  Assume that \begin{equation}\label{cond}
(1+|\Delta_0|)\epsilon<1/4,\end{equation}
then for any local function $f:\bbX \to \bbR$ it holds that
 \begin{equation}
\label{CP}
\var(f)\le 4 \sum_{\bx\in\mathbb Z^d} \nu\left(r_\bx \var_\bx(f)\right).
  \end{equation}
\end{proposition}

\begin{proof}[Sketch]
The starting point to prove \eqref{CP} is the inequality \cite{Martinelli19}*{Lemma 2.5}
\begin{equation}
\label{AB}
\var(f)\leq\sum_{\bx\in\bbZ^d}
\nu\left(\var_\bx \left(\nu_{\bbN^d_{\bx}}(f) \right)
\right),
\end{equation}
which can be  obtained using the law of total variance and the product form of $\nu$.
Then one can rewrite
 a generic term in the r.h.s.\ of \eqref{AB}  using
 the decomposition
\[
\mu_{\bbN^d_{\bx}}(f)=\mu_{\bbN^d_{\bx}}\left(r_\bx f
\right) + \mu_{\bbN^d_{\bx}}\left(\left(1-r_\bx
\right)f\right),
\]
and repeatedly use Cauchy-Schwarz inequality and convexity of the variance to obtain the final result \eqref{CP}.
\end{proof}

\begin{remark}
\label{rem:pos_tree}
Proposition \ref{teo:constrainedpoincare} also plays a key role for the study of KCM on regular trees. It is actually in this context that it was introduced in \cite{Martinelli13} (also see \cite{Cancrini15} for a refinement).
\end{remark}

In \cite{Martinelli19}, Proposition~\ref{teo:constrainedpoincare} is proved in a more general setting allowing the constraints $r_\bx$ to be the  product of several indicator functions and at the same time transforming \eqref{cond} into a more flexible  condition involving the supports and the probabilities of these events. Using this, Proposition~\ref{3.4} below is deduced in \cite{Martinelli19}*{Proposition 3.4}. Before stating it, we require some more notation.
\begin{definition}[Good path]
\label{def:good:path}
We call a sequence of sites $\gamma=(\bx_i)_{i\ge 0}^k$ with $\bx_i-\bx_{i-1}\in\{\be_j:j\in\{1,\dots,d\}\}$ for any $i\in\{1,\dots,k\}$ an \emph{oriented path of length $k+1$}. Fix a \emph{good} and a \emph{super good} event $\cG,\cS\cG\subset\bbX_0$. Given a configuration $\eta\in\bbX$, we say that the oriented path $\gamma$ is \emph{good}, if $\eta_{\bx_i}\in\cG$ for all $i\in\{0,\dots,k\}$ and $\eta_{\bx_k}\in\cS\cG$. For any $\bx\in\bbZ^d$ and $K\ge 1$, we denote by $\Gamma_\bx^{K}$ the event that for every $i\in\{1,\dots,d\}$ there exists a good oriented path of length at most $K$ starting at $\bx+\be_i$.
\end{definition}
\begin{proposition}[Good path constraint dynamics]
\label{3.4}

There exists $\delta>0$ such that for all $p_1,p_2\in[0,1]$ with $\max(p_2,(1-p_1)(\log(1/p_2))^2)\le \delta$, the following holds. If $\nu_0(\cG)=p_1$ and $\nu_0(\cS\cG)=p_2$, then for any local function $f:\bbX\to\bbR$,
\begin{equation}\label{eq:3.4}
\var (f)\leq 4\sum_{\bx\in\bbZ^d}\nu\left(\1_{\Gamma_\bx^{p_2^{-2}}}\var_{\bx}(f)\right).
\end{equation}
\end{proposition}

\subsection{Combining renormalisation and the long range Poincar\'e inequality}
\label{sec:FA2fbetter}
We next explain how to apply Proposition~\ref{3.4} to FA-2f. As suggested by the fact that we worked with a general state space in Section~\ref{sec:constrainedpoincare}, we intend to use renormalisation. Once again, the argument works in any dimension, but we present it in two-dimensions in order to simplify notation. For any $\bx\in\bbZ^2$, the (renormalised) site $\bx$ will correspond to the rectangle $R_\bx=\ell\cdot\bx+R(\ell,\ell)$ (recall \eqref{eq:def:rectangle}), where $\ell=\lceil3\log(1/q)/q\rceil$. We consider the state space $\bbX_0=\Omega_{R_0}$ and the corresponding states $\eta_\bx=\omega_{R_\bx}$ for $\bx\in\bbZ^2$. 

It remains to choose the good and super good events in Definition~\ref{def:good:path}. A simple choice is to define $\omega\in\cS\cG\subset\Omega_{R_0}$, if 
\begin{equation}
\label{eq:def:SG:naive}
\omega_{R(\ell,\ell)\setminus ((1,1)+R(\ell-2,\ell-2))}=\bzero.\end{equation}
That is, a renormalised site is super good, if the perimeter of the corresponding rectangle is empty. The good event is defined by $\omega\in\cG\subset\Omega_{R_0}$, if
\begin{equation}
\label{eq:def:G}\forall i\in\{0,\dots,\ell-1\},\omega_{(i,0)+R(1,\ell)}\neq\bone,\omega_{(0,i)+R(\ell,1)}\neq\bone.\end{equation}
That is a renormalised site is good, if each of its rows and columns contain at least one empty site. The idea behind \eqref{eq:def:SG:naive} and \eqref{eq:def:G} is that the following statements hold (recall the BP closure $[\cdot]$ from \eqref{eq:def:closure}):
\begin{itemize}
    \item $\cS\cG\subset\{\omega\in\Omega_{R_0}:[\omega]=\bzero_{R_0}\}$;
\item for any $\bx\in\{\be_1,-\be_1,\be_2,-\be_2\}$, if $\omega_{R_0}\in\cS\cG$ and $\omega_{R_{\bx}}\in\cG$, then $[\omega_{R_0\cup R_\bx}]=\bzero$;
\item if $\omega_{R_0}\in\cS\cG$, $\omega_{R_{((-1,0)}}\in\cG$ and $\omega_{R_{(0,-1)}}\in\cG$, then $[\omega_{[-\ell,\ell)^2}]=\bzero$.
\end{itemize}
The reader is invited to verify these deterministic BP claims. They will allow us to transport a super good renormalised site along a good path.

We now have all the ingredients necessary to deal with FA-2f, following \cite{Martinelli19} to prove the following bound greatly improving on the one from Corollary~\ref{cor:FA2f:bound:trivial}.
\begin{theorem}[{FA-2f upper bound up to logarithmic corrections}]
    \label{th:FA2f:MT}
    For FA-2f in $d\ge2$ dimensions there exists a constant $C>0$ such that for any $q$ small enough,
    \[\trel\le \exp\left(\frac{C(\log(1/q))^2}{q^{1/(d-1)}}\right).\]
\end{theorem}
\begin{proof}[Sketch]
The reasoning proceeds in three steps. We first apply Proposition~\ref{3.4} as described above. We then need to bound the generic term in the right hand side of \eqref{eq:3.4}, say for $\bx=0$. The occurrence of the event $\Gamma_0^{p_2^{-2}}$ guarantees good paths of length at most $p_2^{-2}$. 

We can view a good path as a one-dimensional general KCM with FA-1f constraint, as covered by Proposition~\ref{prop:FA1f:East:generalised}. Roughly speaking, this enables us to transform the generic term into a sum of terms of the form \begin{equation}
\label{eq:long-range:short-range}\mu\left(\1_{\omega_{R_{\by}}\in\cS\cG}\var_{R_{\bx}}(f|\cG)\right)
\end{equation}
for $\bx$ and $\by$ neighbours along the good path. The prefactor incurred in this transformation is at most $p_2^{-C}$ for some constant $C>0$. Note that \eqref{eq:long-range:short-range} only features a short range constraint.

The final step is to bound \eqref{eq:long-range:short-range} by the Dirichlet form of FA-2f in $R_{\bx}\cup R_\by$ with $\bone$ boundary condition. This can be done using canonical paths. If done brutally, using the legal paths of Section~\ref{sec:legal:paths}, this would lead to a prefactor of order $q^{-2\ell^2}$. However, proceeding a little more carefully (only creating one or two empty columns/rows and moving them rather than emptying the entire rectangles), one can easily reduce this cost to $q^{-2\ell}$. Putting everything together and computing $p_2\ge q^{4\ell}$, we obtain the desired result (see \cites{Martinelli19,Martinelli19a} for more details).
\end{proof}

The three steps in the last proof can be understood as follows. The first step (applying Proposition~\ref{3.4}) allows us to reduce the study of the infinite volume KCM to one in a large but finite volume ($1/p_2^2$, which is exponential in $1/q$) containing sufficient empty sites to efficiently empty the origin. The second step goes from the `global' scale $p_2^{-2}$ to the much smaller `mesoscopic' scale $\ell$. Finally the third step goes from the mesoscopic to the scale of the lattice. In the proof of Theorem~\ref{th:FA2f:MT}, we used a generalised FA-1f dynamics for the second and third steps, as well as a very simple choice of $\cS\cG$ event. In order to go further, we will need to reconsider all these choices.

\section{Sharp threshold}
\label{sec:sharp}
Our next goal is a much stronger result that Theorem~\ref{th:FA2f:MT} providing a sharp threshold for FA-2f. This is the most precise result for KCM of this type and is due to Hartarsky, Martinelli and Toninelli \cite{Hartarsky23FA}*{Theorem~1.3} (also see that reference for a quantitative bound on the second order correction in two dimensions).
\begin{theorem}[Sharp threshold for FA-2f]
\label{th:FA2f:sharp:threshold}
For FA-2f in $d\ge 2$ dimensions,
    \begin{equation}
    \label{eq:FA2f:sharp:threshold}q\log \tau_0\to d\cdot\lambda(d,2)\end{equation}
    in $\bbP_{\mu_q}$-probability, as $q\to 0$, where $\lambda(d,2)$ is the constant from Theorem~\ref{th:BP} for $j$-neighbour BP. Furthermore, $q\log\bbE_{\mu_q}(\tau_0)\to d\cdot\lambda(d,2)$.
\end{theorem}
In other words, for FA-2f we have $\tau_0=(\tbp)^{d+o(1)}$ as $q\to0$.

Before moving on to the proof, let us mention that it would be good to improve the statement above to a relaxation time result (the lower bound follows from \eqref{eq:tau0} and \eqref{eq:FA2f:sharp:threshold}).
\begin{conjecture}[{FA-2f relaxation time}]
For FA-2f in $d\ge 2$, as $q\to0$, we have
\[q\log \trel\to d\cdot\lambda(d,2).\]
\end{conjecture}

\subsection{Lower bound: combinatorial bottleneck}
\label{sec:FA2f:lower}
We start by discussing the lower bound of \eqref{eq:FA2f:sharp:threshold}. It is a relatively simple consequence of known results in BP and generalises well to other models. Recall that from \eqref{eq:t0t0bp} (also recall \eqref{eq:tau0}) and Theorem~\ref{th:BP}, we have $\liminf_{q\to0}q\log\tau_0\ge \lambda(d,2)$, but we would like to improve this bound. The heuristics behind the improvement is the following combinatorial bottleneck (see \cite{Hartarsky23FA}*{Section~2} for more details).
\begin{proof}[Sketch of the lower bound of Theorem~\ref{th:FA2f:sharp:threshold}]
According to Theorem~\ref{th:BP} and Corollary~\ref{cor:legal:path}, with high probability, the origin cannot be emptied using only empty sites within distance, say, $1/q^3$ from it in the initial configuration. We may therefore consider the first time $\tau$, when the origin can be emptied using only empty sites at distance at most $1/q^3$ from it. That is, setting $\ell=\lceil1/q^3\rceil$ and $\Lambda=[-\ell,\ell]^d$, we define
\[\tau=\inf\left\{t\ge 0:\left[\omega_\Lambda(t)\right]^{\bone_{\bbZ^d\setminus\Lambda}}\right\},\]
where $\omega_\Lambda(t)$ is the restriction of the stationary FA-2f process to $\Lambda$ at time $t$ and we recall \eqref{eq:def:closure} and \eqref{eq:def:BP:finite}. 

The crucial observation is that in $\omega_\Lambda(\tau)$, there is a site $\bx$ at the boundary of $\Lambda$ such that
\begin{equation}
\label{eq:FA2f:bottleneck}\left[\omega_\Lambda(\tau)\right]^{\bone_{\bbZ^d\setminus\Lambda}}_0=1\neq\left[\omega_\Lambda^\bx(\tau)\right]^{\bone_{\bbZ^d\setminus\Lambda}}_0.\end{equation}
That is, at $\tau$ a site at the boundary of $\Lambda$ becomes empty and this is essential (pivotal in the percolation jargon) to being able to empty the origin inside the box $\Lambda$.

Since we are working with the stationary process, we may perform a union bound over the attempted updates in $\Lambda$ as in Section~\ref{subsubsec:bottleneck:deduction}. Thus, it suffices to show that for any $\varepsilon>0$ and $q$ small enough, 
\begin{equation}
\label{eq:FA2f:bottleneck:bound}
\mu_q(\cA)\le \exp\left(\frac{-d\cdot\lambda(d,2)+\varepsilon}{q}\right),
\end{equation}
where $\cA$ is the union over $\bx\in\partial\Lambda$ of the event in \eqref{eq:FA2f:bottleneck}.

Note that \eqref{eq:FA2f:bottleneck:bound} only makes reference to 2-neighbour BP, so we have successfully reduced the problem for FA-2f to its BP counterpart. The bound \eqref{eq:FA2f:bottleneck:bound} is indeed known in BP and is, in fact, the main step in the proof of Theorem~\ref{th:BP} for this update family.
\end{proof}

In rough terms, we harnessed the fact that, in order for the origin to be `locally emptiable', a certain `critical droplet' (similar to the $\cS\cG$ event in Section~\ref{sec:FA2fbetter}) needs to be present at the origin at some point in time. We then plugged this combinatorial bottleneck into the standard bound of Section~\ref{subsubsec:bottleneck:deduction}.

\subsection{Coalescing and branching simple symmetric exclusion process}
\label{subsec:CBSEP}
We require one more ingredient as preparation for the upper bound in Theorem~\ref{th:FA2f:sharp:threshold}. It comes in the form of a model strongly related to the generalised FA-1f of Section~\ref{sec:general:KCM}, but not belonging to the class of KCM we defined in Section~\ref{subsec:markov}. For the sake of simplicity, we introduce only the binary (non-generalised) model and refer to \cite{Hartarsky22CBSEP} for the generalised version.

Let $G=(V,E)$ be the box $V=\{1,\dots,\ell\}^d$ with its usual graph structure, where $\ell$ is a positive integer. We consider the state space $\Omega=\{0,1\}^V$ as usual. We define $\Omega_+=\Omega\setminus \{\bone\}$ to be the event that there exists at least one empty site. Similarly, for any edge $e=\{x,y\}\in E$ we refer to $(\omega_x,\omega_y)\in \{0,1\}^{\{x,y\}}$ as the state of $e$ in $\omega$ and write $E_e=\{\omega\in\Omega:\omega_x\omega_y=0\}$ for the event that $e$ is not occupied (at least one of its vertices is empty). Given $p\in
(0,1)$, let $\pi=\bigotimes_{x\in V}\pi_x$ be the product Bernoulli measure in which each vertex is empty with probability $p$ and let $\mu(\cdot):=\pi(\cdot| \Omega_+)$. Given an edge
$e=\{x,y\}$, we write
$\pi_e:=\pi_x\otimes\pi_y$.

\emph{CBSEP} is a continuous time Markov chain on $\Omega_+$ for which the state of any edge $e\in E$ such that $E_e$ occurs  is resampled with rate one w.r.t.\ $\pi_e(\cdot| E_e)$. Thus, any edge containing exactly one empty site moves the empty site to the other endpoint of the edge (the \emph{SEP move}) with rate $(1-p)/(2-p)$ and creates an extra empty site at the occupied endpoint
(the \emph{branching move})
with rate $p/(2-p)$. Moreover, any edge containing two
empty sites occupies one of the two chosen uniformly (the \emph{coalescing move}) with rate $2(1-p)/(2-p)$. The chain is readily seen to be reversible with respect to $\mu$ and
ergodic on $\Omega_+$, because it can reach the configuration with all sites empty. If $c(\omega,\omega')$ denotes the jump rate from $\omega$ to $\omega'$, the Dirichlet
form $\cD^{\mathrm{CBSEP}}(f)$ of the chain has the expression
\begin{align}
\label{eq:def:cD:CBSEP}
\cD^{\mathrm{CBSEP}}(f)&{}=\frac 12 \sum_{\omega,\omega'\in\Omega_+}\mu(\omega)c(\omega,\omega')\left(f(\omega')-f(\omega)
  \right)^2\\
  &{}=\sum_{e\in E}\mu(\1_{E_e}\var_{\pi_e}(f| E_e)).\nonumber\end{align}
  
Notice that the branching and coalescing moves of CBSEP are exactly the moves allowed in FA-1f (recall Section~\ref{sec:FA1f}). Moreover, the SEP move can be reconstructed by a branching and a coalescing move. This leads to a comparison between the corresponding Dirichlet forms (see e.g.\ \cite{Levin09}*{Section 13.4}).

Although the two models are clearly closely related, we stress that CBSEP has many advantages over \mbox{FA-$1$f}, making its study simpler. Most notably, {CBSEP} is \emph{attractive} in the sense of interacting particle systems, that is the natural stochastic order is preserved by the dynamics (see \cite{Liggett05}*{Section III.2} for background). Furthermore, it is natural to embed in {CBSEP} a continuous time random walk $(W_t)_{t\ge 0}$ on $G$ such that {CBSEP} has an empty site at $W_t$ for all $t\ge 0$. This feature is challenging to reproduce for \mbox{FA-$1$f} \cite{Blondel13}. Finally, it is possible to move an empty site in CBSEP without creating more empty sites, contrary to what is the case in FA-1f (recall Section~\ref{sec:FA1f}).

The main result on CBSEP we will need is the following \cite{Hartarsky23FA}*{Proposition 5.2}.\footnote{We record a mistake in \cite{Hartarsky23FA} leading to the weaker bound stated here. Indeed, in the last but one sentence of the proof of \cite{Hartarsky23FA}*{Proposition~5.2} in Appendix~B, the bounds on the cover time of the random walk and logarithmic Sobolev constant of CBSEP are not correctly imported from \cite{Levin09}*{Chapter 11} and \cite{Hartarsky22CBSEP}*{Corollary 3.2}. This mistake has no impact on the rest of the paper.}
\begin{proposition}[CBSEP relaxation time]
\label{prop:CBSEP}
Assume that $d\ge2$ and consider a sequence of box sizes $\ell_n$ and parameters $p_n$ such that $p_n\ell_n^d\to\infty$ and $p_n\to0$. Then, for some $C>0$ and all $n$ large enough,
\[\trel^{\mathrm{CBSEP}}\le \frac{C\log^{3}(1/p_n)}{p_n}.\]
\end{proposition}
\begin{proof}[Sketch]
The first step of the proof is to renormalise the model by considering boxes of volume approximately $1/p_n$. This brings us to treating CBSEP with parameter $p$ of order 1 on arbitrary volume and CBSEP on a box of volume approximately $1/p_n$. The former relaxation time is uniformly bounded (see \cite{Hartarsky23FA}*{Lemma B.1} and \cite{Hartarsky22CBSEP}*{Theorem~1} for details), while the latter is bounded by $C\log^2(1/p_n)/p_n$ thanks to \cite{Hartarsky22CBSEP}*{Corollary 3.1}.

The proof of the first upper bound follows from the fact that $\trel^{\mathrm{FA-1f}}$ is bounded for $q$ bounded away from 0 (recall Theorem~\ref{th:exponential}) together with a comparison between the Dirichlet forms of CBSEP and FA-1f. The second upper bound may be proved, using canonical paths along the lines of Theorem~\ref{th:FA1f} (also recall Conjecture~\ref{conj:FA1f}, whose upper bound is known from \cite{Cancrini08} for $d=2$, and its analogue for $d\ge3$ from \cite{Shapira20}), see \cite{Hartarsky22CBSEP}*{Proposition~4.6}.

While this concludes the sketch for CBSEP, its generalised version is more subtle to analyse. Indeed, in \cite{Hartarsky22CBSEP}*{Theorem~2} it is shown that one can bound the mixing time (and therefore the relaxation time) of generalised CBSEP, using the mixing time of CBSEP and the cover time of the continuous time simple random walk on the box of interest. The cover time is classically bounded by $\log^2(1/p_n)/p_n$ \cite{Levin09}*{Chapter 11}, but more work is needed to bound the mixing time of CBSEP, as opposed to its relaxation time. The approach of \cite{Hartarsky22CBSEP} is to prove a logarithmic Sobolev inequality, which classically bounds the mixing time. This is achieved, using deep input from \cites{Lee98,Alon20}. An alternative approach to bounding the mixing time of FA-1f can be found in \cites{Pillai17,Pillai19}.
\end{proof}

\subsection{Upper bound}
\label{sec:FA2f:upper}
We are now ready to discuss the proof of the upper bound of Theorem~\ref{th:FA2f:sharp:threshold} in dimension $d=2$, following \cite{Hartarsky23FA}. On the high level, the proof resembles Section~\ref{sec:FA2fbetter}. The first two steps bringing us to a mesoscopic scale are quite general, while reaching the mesoscopic scale from the microscopic one is much more delicate and specific to FA-2f.
\subsubsection{Reduction to the mesoscopic scale via CBSEP}
\label{subsec:FA2f:reduction:to:meso}
Recall from Theorem~\ref{th:BP} that $\lambda(2,2)=\pi^2/18$. Fixing some $\varepsilon>0$ and setting 
\begin{equation}
\label{eq:def:t*}
t_*=\exp\left(\frac{\pi^2
+\varepsilon}{9q}\right),\end{equation}
our goal is to prove that 
\begin{equation}\label{eq:FA2f:goal}
\lim_{q\to 0}\bbP_\mu\left(\tau_0>t_*\right)=0.
\end{equation}
We next use finite speed of propagation (recall \eqref{eq:speed:propagation}) to show that 
it is unlikely that, before time $t_*$, the state of the origin is influenced by the configuration outside the box $B=[-2t_*,2t_*]^2$. Thus, it suffices to prove \eqref{eq:FA2f:goal} for FA-2f with  boundary condition $\bone_{\bbZ^2\setminus B}$.

The next step is a renormalisation to the generalised CBSEP model discussed in Section~\ref{subsec:CBSEP}. The renormalisation itself resembles the one of Section~\ref{sec:FA2fbetter}, but the choice of super good events is much more delicate. We fix a mesoscopic scale $\ell=\lceil q^{-9}\rceil$ (the exponent is fairly arbitrary) and divide $B$ into boxes $Q_\bx=\ell\bx+R(\ell,\ell)$ (recall \eqref{eq:def:rectangle}) for $\bx\in\bar B=[-\lceil 2t_*/\ell\rceil,\lceil 2t_*/\ell\rceil]\cap \bbZ^2$. Each renormalised site $\bx\in\bar B$ will be in one of two states---good and super good. As in Section~\ref{sec:FA2fbetter}, $\bx\in\bar B$ is \emph{good}, if $\omega_{Q_\bx}$ has at least one empty site on each row and column of $Q_\bx$. We denote the corresponding event by $\cG_\bx$.

By a union bound, it holds that
\[1-\mu\left(\bigcap_{\bx\in\bar B}\cG_\bx\right)\le 5t_*^2(1-q)^\ell\le 1/t_*^4,\]
owing to our choice of sufficiently large $\ell$. Since this probability is so small, as in Section~\ref{subsubsec:bottleneck:deduction} (union bound over attempted updates, using stationarity), it is likely that all renormalised sites in $\bar B$ remain good at all times up to $t_*$. 

We will choose the super good event $\cS\cG_\bx$ in such a way that
\begin{equation}
\label{eq:SG:proba}
\mu\left(\cS\cG_\bx\right)\ge \exp\left(-\frac{\pi^2/9+\varepsilon/2}{q}\right).\end{equation}
Comparing this with \eqref{eq:def:t*}, we similarly compute that it is likely that, up to time $t_*$, at all times there is at least one super good renormalised site in $\bar B$. Thus, we therefore assume that
\[\cE=\bigcap_{\bx\in\bar B}\cG_\bx\cap\bigcup_{\bx\in\bar B}\cS\cG_\bx\]
occurs for all $t\le t_*$. The event $\cE$ corresponds to $\Omega_+$ in Section~\ref{subsec:CBSEP} and super good renormalised sites play the role of empty sites for CBSEP.

By a standard result (see \cite{Asselah01}*{Theorem~2}) similar to \eqref{eq:F0F1} (also recall \eqref{eq:def:Trel}), but taking into account that we require not exiting $\cE$, in order to prove \eqref{eq:FA2f:goal}, it suffices to establish a Poincar\'e type inequality of the form
\begin{equation}
\label{eq:FA2f:goal2}
\frac{\cD_{\bone_{\bbZ^2\setminus B}}(f)}{\var(f|\cE)}\ge \exp\left(-\frac{\pi^2/9+2\varepsilon/3}{q}\right)
\end{equation}
for any function $f:\Omega_{B}\to\bbR$ such that $f(\omega)=0$, if $\omega\not\in\cE$. We are now ready to apply the Poincar\'e inequality for generalised CBSEP provided by Proposition~\ref{prop:CBSEP}. Recalling \eqref{eq:SG:proba}, this yields
\begin{equation}
\label{eq:FA2f:CBSEP}\var(f|\cE)\le Cq^{-3}\exp\left(-\frac{\pi^2/9+\varepsilon/2}{q}\right)\sum_{\bx\sim\by}\mu\left(\1_{\cS\cG_{\bx,\by}}\var_{Q_\bx\cup Q_\by}\left(f|\cS\cG_{\bx,\by}\right)|\cE\right),\end{equation}
where $\cS\cG_{\bx,\by}=\cS\cG_\bx\cup\cS\cG_\by$ for neighbours $\bx\sim\by$ in $\bar B$. Note that, if we had used FA-1f (recall Section~\ref{sec:FA1f}) instead of CBSEP at this point, the exponent above would become roughly $2\pi^2/(9q)$ instead of $\pi^2/(9q)$, which is not enough for proving Theorem~\ref{th:FA2f:sharp:threshold}.

Notice that relating $\var_{Q_\bx\cup Q_\by}(f|\cS\cG_{\bx,\by})$ to the terms of the Dirichlet form of FA-2f in $Q_\bx\cup Q_\by$ (recall \eqref{eq:dirichlet}) is essentially equivalent to establishing a Poincar\'e inequality for FA-2f on this mesoscopic volume, conditioned to remain in the event $\cS\cG_{\bx,\by}$. Indeed, some simple but delicate manipulations (see \cite{Hartarsky23FA}*{Claim 5.5}) allow deducing \eqref{eq:FA2f:goal2} from \eqref{eq:FA2f:CBSEP} and the Poincar\'e inequality
\begin{align}
\label{eq:FA2f:meso}
\var_{Q}(f|\cS\cG)&{}\le \gamma(Q)\sum_{\bz\in Q}\mu\left(c_\bz^{\bone_{\bbZ^2\setminus Q}}\var_\bz(f)|\cS\cG\right),\\\nonumber\gamma(Q)&{}\le\exp\left(\frac{\varepsilon}{7q}\right)\end{align}
with $Q=Q_\bx\cup Q_\by$ and $\cS\cG=\cS\cG_{\bx,\by}$, which holds for all $f:\Omega_Q\to\bbR$. The next section is devoted to the proof of \eqref{eq:FA2f:meso}.

\subsubsection{Mesoscopic Poincar\'e inequality: the matryoshka doll technique}
\label{subsubsec:matryoshka}
The inequality \eqref{eq:FA2f:meso}, which we seek to establish, should be interpreted as stating that, once a critical droplet is present at a given location, it is rather easy to completely change the state there. Note that we have not yet specified our choice of super good event beyond requiring \eqref{eq:SG:proba} to hold. Our actual choice and the proof of \eqref{eq:FA2f:meso} go hand in hand and follow a multi-scale renormalisation scheme which we refer to as the \emph{matryoshka doll technique}. For the sake of simplicity of the presentation, we rather present the argument in a way that yields a sharp threshold for slight variants of the model. The additional technical difficulties are dealing with more boundary conditions in \eqref{eq:FA2f:meso}, requiring at least one empty site on every two consecutive lines and implementing a microscopic FA-1f dynamics on the boundary of a droplet. The interested reader can find these features in \cite{Hartarsky23FA}. The model we will treat is the KCM corresponding to the modified two-neighbour update family obtained by removing the lowest two rules in Figure~\ref{fig:FA2f}. We refer to it as modified FA-2f and seek to prove \eqref{eq:FA2f:meso} but with $\cS\cG$ event satisfying \eqref{eq:SG:proba} with $\lambda(2,2)$ replaced by $\lambda'=\pi^2/6$, which is the correct sharp threshold constant for this model (the matching lower bound is proved exactly as described in Section~\ref{sec:FA2f:lower}).

The matryoshka doll technique requires us to keep track of several features of a droplet $\Lambda$ in parallel:
\begin{itemize}
    \item geometry, that is, the shape and size of the droplet;
    \item super good event $\cS\cG(\Lambda)\subset\Omega_\Lambda$, which guarantees that the droplet's state can be resampled efficiently, conditionally on the super good event;
    \item occurrence probability, that is, $\mu_q(\cS\cG(\Lambda))$;
    \item relaxation time, that is, the smallest constant $\gamma(\Lambda)\ge 1$ such that \eqref{eq:FA2f:meso} holds.
\end{itemize}
Given the multi-scale nature of the technique, all of the above are defined or bounded recursively for a sequence of droplets $\Lambda^{(m)}$, starting from a single empty site and reaching the droplet $Q$ from Section~\ref{subsec:FA2f:reduction:to:meso}.

\paragraph{Geometry} In the case of (modified) FA-2f, droplets are simply rectangles. Let $m_0=\lfloor1/\sqrt q\rfloor$,
\begin{align}
\label{eq:def:ellm:Lambdam}
\ell_m&{}=\begin{cases}
m+1&m< m_0,\\
\left\lfloor \frac{e^{(m-m_0)\sqrt q}}{\sqrt q}\right\rfloor &m\ge m_0,
\end{cases}&\Lambda^{(m)}&{}=
    R\left(\ell_{\lceil m/2\rceil},\ell_{\lfloor m/2\rfloor}\right),
\end{align}
for any $m\ge 0$. Thus, $\Lambda^{(m)}$ form a nested sequence of rectangular droplets, each second one being a square. At each step only the width or only the length (depending on parity) is increased. The corresponding length scales increase linearly up to $1/\sqrt q$ and then exponentially, reaching our final size of interest $\ell=1/q^9$ after approximately $\frac{17\log(1/q)}{2\sqrt q}$ steps. The exact choice of scales is not of fundamental importance, but some care is needed. The scale $1/\sqrt q$ is known to be relevant thanks to BP results (recall \eqref{eq:2nBP}), but it is also not crucial for our purposes.

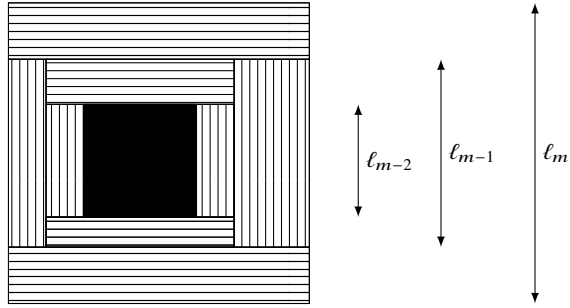
\begin{figure}
  \centering
  \begin{tikzpicture}[>=latex,x=0.5cm,y=0.5cm]
\draw[fill] (10,3.3) rectangle (13,6.3);
\draw[pattern=vertical lines] (9,3.3) rectangle (14,6.3); 
\draw[pattern=horizontal lines] (9,2.5) rectangle (14,3.3); 
\draw[pattern=horizontal lines] (9,2.5) rectangle (14,3.3); 
\draw[pattern=horizontal lines] (9,6.3) rectangle (14,7.5); 
\draw[pattern=horizontal lines] (8,1) rectangle (16,2.5); 
\draw[pattern=horizontal lines] (8,7.5) rectangle (16,9); 
\draw[pattern=vertical lines] (8,2.5) rectangle (9,7.5); 
\draw[pattern=vertical lines] (14,2.5) rectangle (16,7.5); 
\draw [<->] (22,1)--(22,9) node [midway, right] {$\ell_{m}$};
\draw [<->] (17.3,3.3)--(17.3,6.3) node [midway, right]{$\ell_{m-2}$};
\draw [<->] (19.5,2.5)--(19.5,7.5) node [midway, right] {$\ell_{m-1}$};
\draw  (8, 1)-- (8,9);
\draw  (16,1)-- (16,9);
\draw  (8, 1)-- (16,1);
\draw  (8, 9)-- (16,9);
\draw  (8, 7.5)-- (16,7.5);
\draw  (8, 2.5)-- (16,2.5);
\draw  (9, 2.5)-- (9,7.5);
\draw  (14, 2.5)-- (14,7.5);
\draw  (9, 3.3)-- (14,3.3);
\draw  (9, 6.3)-- (14,6.3);
\draw  (10, 3.3)-- (10,6.3);
\draw  (14, 3.3)-- (14,6.3);
\end{tikzpicture}
\caption{An example structure of the super good event $\cS\cG(\Lambda^{(2m)})$. Only the super good translates of $\Lambda^{(n)}$ for $m-n\in\{0,\dots,6\}$ are shown. The hatched regions are required to contain at least one empty site per line in the direction of the hatching.}
\label{fig:FA2f:droplets}
\end{figure}

\paragraph{Super good event}
The definition of the super good event is guided by the idea that lower scale droplets should be allowed to move freely within a larger scale droplet. This is vital for obtaining an efficient relaxation mechanism (that is, a small Poincar\'e constant $\gamma(Q)$ in \eqref{eq:FA2f:meso}). For any $\bx\in\bbZ^2$, this intuition leads us to define $\cS\cG(\{\bx\})=\{\omega_\bx=0\}$ (recall that $\Lambda^{(0)}=\{0\}$ is just a single site) and, for $m\ge 0$,
\begin{multline}
\label{eq:def:SG}
\cS\cG\left(\bx+\Lambda^{(2m+1)}\right)=\bigcup_{s=0}^{\ell_{m+1}-\ell_m}\Big(\cS\cG\left((s,0)+\bx+\Lambda^{(2m)}\right)\\\cap\bigcap_{\substack{t\in\{0,\dots,\ell_{m+1}-1\}\\t\not\in\{s,\dots,s+\ell_m-1\}}}\left\{\omega\in\Omega_{\bx+\Lambda^{(2m+1)}}:\omega_{(t,0)+\bx+R(1,\ell_m)}\neq\bone\right\}\Big)
\end{multline}
and similarly for odd scales. In words, at each scale we require that some translate of the lower scale droplet occurs and that each of remaining rows or columns (depending on parity) contains at least one empty site. The definition is illustrated in Figure~\ref{fig:FA2f:droplets}.

\paragraph{Occurrence probability}
We next turn to lower bounding $\mu(\cS\cG(\Lambda^{(m)})$. The argument goes back to \cite{Aizenman88} and is only about BP, but we include it, as it is informative. Set $f:(0,\infty)\to(0,\infty):x\mapsto-\log(1-e^{-x})$, which is decreasing and convex. Systematically taking $s=0$ in \eqref{eq:def:SG}, we get that for any $m\in[0,17\log(1/q)/\sqrt q]$,
\begin{align*}\mu\left(\cS\cG\left(\Lambda^{(m)}\right)\right)&{}\ge \prod_{n=0}^{\lceil m/2\rceil}\left(1-(1-q)^{\ell_n}\right)^{2(\ell_{n+1}-\ell_n)}\\
&{}\ge \exp\left(-2\sum_{n=0}^{\lceil m/2\rceil}(\ell_{n+1}-\ell_n)f(q\ell_n)\right)\\
&{}=\exp\left(-\frac2q\int_0^\infty f+\frac{o(1)}{q}\right)=\exp\left(-\frac{\pi^2+o(1)}{3q}\right),\end{align*}
where the asymptotic notation is as $q\to 0$. The error in the approximation of the Riemann sum by an integral above is carried out in more detail in \cite{Hartarsky23FA}*{Appendix A}, while the integral can be computed using the series expansion of $\log(1-\cdot)$ and the fact that $\zeta(2)=\pi^2/6$. We have thus concluded the proof of the analogue of \eqref{eq:SG:proba} in the context of modified FA-2f.

\paragraph{Relaxation time}
Proving \eqref{eq:FA2f:meso} is the most challenging part of \cite{Hartarsky23FA}. The proof proceeds by establishing the recursive bound 
\begin{equation}
\label{eq:FA2f:recursion}\gamma\left(\Lambda^{(m)}\right)\le e^{C\log^2(1/q)}\gamma\left(\Lambda^{(m-1)}\right).\end{equation}
for some constant $C>0$. In turn, proving \eqref{eq:FA2f:recursion} is done via a version of the bisection technique, whose base step is unusually delicate. 

For concreteness, we only discuss one parity. Fix some $m$ such that $\ell_{m}\le 9\log(1/q)/q$. Set $\Lambda'=\Lambda^{(2m+1)}=R(\ell_{m+1},\ell_m)$ and $\Lambda=\Lambda^{(2m)}=R(\ell_m,\ell_m)$. The bisection will be used to decrease the width difference $\ell_{m+1}-\ell_m$ down to 1 in roughly $\log(\ell_{m+1}-\ell_m)$ steps. For translates of $R(l,\ell_m)$ with $l\in(\ell_m,\ell_{m+1})$, we extend \eqref{eq:def:SG} by replacing $\ell_{m+1}-\ell_m$ by $l-\ell_m$ and $\Lambda^{(2m+1)}$ by $R(l,\ell_m)$. We consider the rectangles $R^{(k)}=R(\ell_m+d_k,\ell_m)$ with $d_k=\lceil (2/3)^k(\ell_{m+1}-\ell_m)\rceil$ for $k\ge 0$. We then seek to prove a recursive bound of the form
\begin{equation}
\label{eq:FA2f:recursion:step}
\gamma\left(R^{(k)}\right)\le a_k\gamma\left(R^{(k+1)}\right)\end{equation}
with $a_k\le q^{-C}$ for some constant $C>0$.

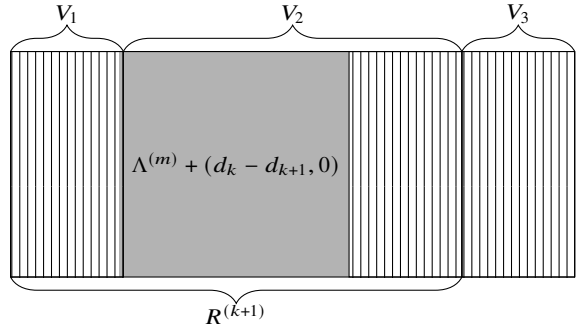
\begin{figure}
\centering
\begin{tikzpicture}[>=latex,x=0.5cm,y=0.5cm]
\fill[opacity=0.3] (-0.1,6) -- (-0.1,0) -- (6,0) -- (6,6) -- cycle;
\draw (-3,0)-- (12,0);
\draw (-3,6)-- (12,6);
\draw (9,0)-- (9,6);
\draw (12,0)-- (12,6);
\draw (9,6)-- (0,6);
\draw (-3,6)-- (-3,0);
\draw[pattern=vertical lines] (-3,0) rectangle (0,6); 
\draw[pattern=vertical lines] (9,0) rectangle (12,6); 
\draw[pattern=vertical lines] (6,0) rectangle (9,6); 
\draw (3,3) node {$\Lambda^{(m)}+(d_k-d_{k+1},0)$};
\draw [decorate,decoration={brace,amplitude=10pt}]   (9,0) -- (-3,0) node [midway,yshift= - 0.5cm] {$R^{(k+1)}$};
\draw [decorate,decoration={brace,amplitude=10pt}] (-3,6) -- (0,6) node [midway,yshift= 0.5cm] {$V_1$};
\draw [decorate,decoration={brace,amplitude=10pt}] (0,6) -- (9,6) node [midway,yshift= 0.5cm] {$V_2$};
\draw [decorate,decoration={brace,amplitude=10pt}] (9,6) -- (12,6) node [midway,yshift= 0.5cm] {$V_3$};
\end{tikzpicture}
\caption{
\label{fig:3block}The partition of $R^{(k)}$ into the rectangles $V_1$, $V_2$, $V_3$.}
\end{figure}

In order to prove \eqref{eq:FA2f:recursion:step}, we use an auxiliary dynamics somewhat similar to the two-block one of Lemma~\ref{lem:two-block}, but using a different mechanism. Contrary to Lemma~\ref{lem:two-block}, which uses an East mechanism, the dynamics behind Lemma~\ref{lem:three-block} below is non-oriented. It has three sites and has two types of constrained updates occurring at rate 1. The first update resamples the first and second sites (blocks) conditionally on some event occurring there before and after the update. The second update is similar for the second and third sites. The first two blocks together correspond to $R^{(k+1)}$, all three blocks form $R^{(k)}$ and the second and third blocks form a translate of $R^{(k+1)}$ (see Figure~\ref{fig:3block}). This leads us to the following lemma adapted from \cite{Hartarsky23FA}*{Proposition 3.5}.

\begin{lemma}[{Three-block dynamics}]
\label{lem:three-block}There exists $C>0$ such that the following holds. Let $(\bbX,\pi)$ be the product of three finite probability spaces $(\bbX_i,\pi_i)_{i=1}^3$. Fix some events $\cA_1\subset\bbX_1$, $\cA_3\subset\bbX_3$, $\cB_{1,2}\subset\bbX_1\times\bbX_2$ and $\cB_{2,3}\subset\bbX_2\times\bbX_3$. Set $\cH=\cB_{1,2}\times\cA_3$ and $\cK=\cA_1\times \cB_{2,3}$. Consider the Dirichlet form
\[\cD(f)=\pi\left(\1_{\cH}\var(f|\cB_{1,2}\times\{\omega_3\})+\1_{\cK}\var(f|\{\omega_1\}\times\cB_{2,3})|\cH\cup\cK\right)\]
defined for any $f:\bbX\to\bbR$. Consider some event $\cC_2\subset\bbX_2$ such that $\cA_1\times\cC_2\times\cA_3\subset(\cH\cap\cK)$. Then
\begin{equation}
\label{eq:3block}
\var_{1,2,3}(f|\cH\cup\cK)\le C\frac{\pi_{1,2}(\cB_{1,2})\pi_{2,3}(\cB_{2,3})}{\pi_1(\cA_1)(\pi_2(\cC_2))^2\pi_3(\cA_3)}\cD(f)\end{equation}
\end{lemma}
\begin{proof}[Sketch]
The mechanism behind the proof is the following. Consider two copies of the chain described above with Dirichlet form $\cD$, coupled by attempting the same updates. Observe that the following sequence of update attempts guarantees that the two chains reach the same state. First attempt to resample sites 1 and 2 to so that $\cA_1\times\cC_2$ occurs. Then, before any other update is attempted, update sites 2 and 3 (after the first update, $\cK$ necessarily occurs, regardless whether the first update attempt was successful) so that $\cC_2\times\cA_3$ occurs. Finally, again before any other update is attempted, update sites 1 and 2. Estimating the rate at which such a sequence of updates is attempted yields the desired result (see \cite{Hartarsky23FA}*{Proposition 3.5} for more detail).
\end{proof}

In order to prove \eqref{eq:FA2f:recursion:step}, using Lemma~\ref{lem:three-block}, we consider the blocks $V_1,V_2,V_3$ depicted in Figure~\ref{fig:3block}. The events $\cA_1$ and $\cA_3$ require each column to contain at least one empty site. The events $\cB_{1,2}$, $\cB_{2,3}$ and $\cC_2$ are the super good events for $V_1\cup V_2$, $V_2\cup V_3$ and $V_2$, which were already defined. In order to show that the fraction in \eqref{eq:3block} is at most $q^{-C}$ for some $C>0$, it remains to observe that 
\begin{equation}
\label{eq:uniform position}
\pi_1(\cA_1)\pi_2(\cC_2)\ge \pi_{1,2}(\cB_{1,2})/d_k.
\end{equation}
The last inequality itself follows from the observation that the event in \eqref{eq:def:SG} for each $s$ has equal probability and the position depicted in Figure~\ref{fig:3block} guarantees the occurrence of $\cH\cap\cK$.

The above proves \eqref{eq:FA2f:recursion:step} as desired. However, we are not done proving \eqref{eq:FA2f:recursion}. Indeed, iterating \eqref{eq:FA2f:recursion:step} yields
\begin{equation}
\label{eq:FA2f:recursion:step:final}\gamma\left(\Lambda^{(2m+1)}\right)\le e^{C\log^2(1/q)}\gamma\left(R\left(\ell_{m}+1,\ell_m\right)\right)\end{equation}
for some $C>0$, the rectangle on the right being the thinnest rectangle $R^{(k)}$ we can obtain, which is one column wider than the desired $\Lambda^{(2m)}$. The reason for this is visible in Figure~\ref{fig:3block}. Indeed, if $d_{k}=1$, we cannot fit a translate of $\Lambda^{(2m)}$ in $V_2$ in such a way that both $V_1$ and $V_3$ remain non-empty. Thus, this base case requires a separate argument contained in \cite{Hartarsky23FA}*{Proposition 3.7, Lemma 4.10}, which we briefly discuss next.

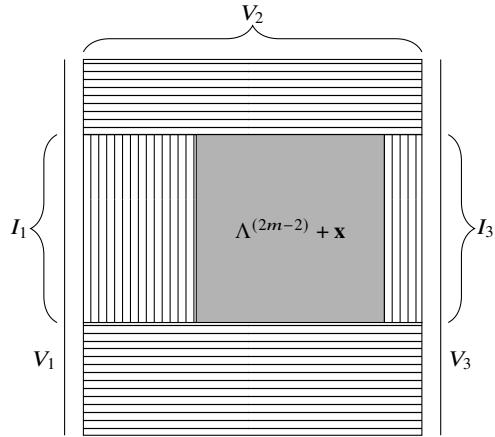
\begin{figure}[b]
\centering
\begin{tikzpicture}[>=latex,x=0.5cm,y=0.5cm]
\fill [opacity=0.3] (3,3) rectangle (8,8);
\draw (5.5,5.5) node {$\Lambda^{(2m-2)}+\bx$};
\draw (-0.5,0)--(-0.5,10);
\draw (9.5,0)--(9.5,10);
\draw [pattern=vertical lines] (0,3) rectangle (3,8);
\draw [pattern=vertical lines] (8,3) rectangle (9,8);
\draw [pattern=horizontal lines] (0,8) rectangle (9,10);
\draw [pattern=horizontal lines] (0,0) rectangle (9,3);
\draw[anchor=east] (-0.5,2) node {$V_1$};
\draw [decorate,decoration={brace,amplitude=10pt}] (0,10.2) -- (9,10.2) node [midway,yshift=0.5cm] {$V_2$};
\draw [decorate,decoration={brace,amplitude=10pt}] (9.7,8) --
(9.7,3) node [midway,xshift=0.5cm] {$I_3$};
\draw [decorate,decoration={brace,amplitude=10pt}] (-0.7,3) -- (-0.7,8) node [midway,xshift=-0.5cm] {$I_1$};
\draw[anchor=west] (9.5,2) node {$V_3$};
\end{tikzpicture}
\caption{
\label{fig:3block:base}
Partition of $R(\ell_m+1,\ell_m)$ into the rectangle $V_2$ and the lines $V_1$ and $V_3$. The internal structure of the contracted super good event for $V_2$ is shown. The empty sites in $V_1$ and $V_3$ need to be in $I_1$ and $I_3$ respectively in order to match it.}
\end{figure}

The proof that 
\begin{equation}
\label{eq:FA2f:recursion:base}
\gamma\left(R(\ell_m+1,\ell_m)\right)\le q^{-C}\gamma\left(\Lambda^{(2m)}\right)\end{equation}
for a constant $C>0$ proceeds along somewhat similar lines to \eqref{eq:FA2f:recursion:step}, but is more subtle. We do use a decomposition as in Figure~\ref{fig:3block} with $V_1$ and $V_3$ consisting of a single column. In order to take into account the fact that the remaining $V_2$ is slightly smaller than $\Lambda^{(2m)}$, we define a ``contracted'' version of $\cS\cG(\Lambda^{(2m)})$ which requires as much of its internal structure as it is possible to fit into $V_2$. That is, we require the super good translate of $\Lambda^{(2m-2)}$ (assuming $m\neq 0$) and lines containing at least one empty site within $V_2$, mimicking Figure~\ref{fig:FA2f:droplets}. This contracted super good event plays the role of $\cC_2$ in an appropriate analogue of Lemma~\ref{lem:three-block} (the events $\cA_1$ and $\cA_3$ now need to fit well with the position of the super good translate of $\Lambda^{(2m-2)}$, see Figure~\ref{fig:3block:base}). Finally, using a double iteration of \eqref{eq:uniform position}, we bound the resulting factor as desired, since the position of the super good translate of $\Lambda^{(2m-2)}$ is somewhat uniform inside $\Lambda^{(2m)}$.

 Putting \eqref{eq:FA2f:recursion:step:final} and \eqref{eq:FA2f:recursion:base} together, we obtain \eqref{eq:FA2f:recursion}. In turn, iterating \eqref{eq:FA2f:recursion} and a trivial bound at scale $m=0$ gives \eqref{eq:FA2f:meso} with $\gamma(Q)\le \exp(o(1/q))$, as $q\to 0$, as desired. This concludes the sketch of Theorem~\ref{th:FA2f:sharp:threshold}.

\section{Conclusion}
In conclusion, let us summarise the new techniques presented in the present chapter, recalling Section~\ref{subsec:tools:BP} and~\ref{subsec:tools:1d}. 

\noindent\textbf{Combinatorial bottlenecks.} Thanks to this method, we were able to prove precise lower bounds. The approach of Section~\ref{sec:FA2f:lower} is quite general. It is based on BP results, but more precise than the BP lower bound we saw in \eqref{eq:t0t0bp}. The intuition behind it is that a critical droplet needs to visit the origin in order to change its state, which does not happen before the time given by the inverse probability of a droplet.

\noindent\textbf{Bisection.} In Section~\ref{sec:bisection:2d} we saw that the bisection technique can also be applied in higher dimensions by alternating two-block steps in different directions.

\noindent\textbf{Long range renormalisation.} This is a new technique in our repertoire. It allows us to focus on a one-dimensional path leading a critical droplet to the origin. This effectively reduces the problem to treating the movement of a critical droplet on its own scale, rather than a much larger one. It has the advantage of yielding upper bounds on $\trel$, but the one-dimensional character of its path does not allow obtaining sharp thresholds.

\noindent\textbf{Matryoska dolls.} This is our state of the art technique for proving sharp upper bounds. It may be viewed as a very adaptable multi-scale renormalisation scheme, whose flexibility will be unleashed in the next chapter. In a sense, bisection in higher dimensions is also an instance of matryoshka dolls. The idea of the method is to recursively prove Poincar\'e inequalities on progressively larger and larger droplets. At each step, we have full freedom in choosing the auxiliary dynamics. Once a sufficiently large scale is reached via the matryoshka dolls, we conclude by a single step renormalisation, whose auxiliary global dynamics may also be adapted to our needs. A major advance of this technique is that it allows us to prove very precise results incorporating tailored relaxation mechanisms, completely bypassing the need to build any explicit canonical paths. While for modified FA-2f  canonical paths reflecting the multi-scale structure of droplets could be envisioned, such approaches very quickly go out of hand in more general settings such as the ones investigated in the next chapter.

\noindent\textbf{Auxiliary dynamics.}
Let us review  the auxiliary dynamics we have plugged into the matryoshka dolls technique so far.  In Section~\ref{sec:bisection:2d}, we always used the simple two-block dynamics of Lemma~\ref{lem:two-block}. In Section~\ref{sec:FA2fbetter}, we employed a one-dimensional FA-1f global dynamics. In Section~\ref{subsec:CBSEP} we developed a new possible global dynamics given by the CBSEP that was applied to FA-2f in Section~\ref{subsec:FA2f:reduction:to:meso}. Finally, in Section~\ref{subsubsec:matryoshka}, we developed a three-block alternative, Lemma~\ref{lem:three-block}, to the two-block Lemma~\ref{lem:two-block}. The three-block dynamics was used repeatedly in one direction and we then switched direction repeatedly, until reaching the desired mesoscopic scale. The key feature of CBSEP and the three-block lemma is that they allow us to move a very unlikely droplet without paying the price of its creation from scratch, but rather just a little bit of internal reshuffling.


	\chapter{Universality}

\label{chap:universality}
  
\abstract{In this chapter, we consider KCM with completely arbitrary update families. The goal of universality is to identify what kinds of behaviour are possible within this vast collection of models and to efficiently classify all update families with respect to their behaviour. We begin by considering one-dimensional models, which feature three rough universality classes, represented by the three nearest neighbour models studied in Chapter~\ref{chap:1d}. We then provide some background on the two-dimensional BP universality theory, before moving on to two-dimensional KCM. We present the complete rough and refined universality theory for KCM in two dimensions and cover the essential elements of their proofs, building on the techniques exhibited in Chapters~\ref{chap:1d} and~\ref{chap:FA2f}.}

\section{KCM universality in one dimension}
Let us warm up by considering one-dimensional KCM. We recall that in Theorem~\ref{th:FA1f} and~\ref{th:East}, we saw that the correct scaling of $\trel$ (or $\tau_0$) for FA-1f is $1/q^3$, while it is $\exp((\log(1/q))^2/(2\log 2))$ for the East KCM. For FA-2f, we have $\qc=1$, as discussed in Section~\ref{sec:FA2f:1d}. The outcome of the universality theory in one dimension is that the only possible behaviours are those of FA-1f, East and FA-2f. 

In order to state the universality result, we need to define stable directions, which will determine the universality class.
\begin{definition}[Stable directions]
Fix a one-dimensional update family $\cU$. We say that the positive (resp.\ negative) direction is \emph{unstable}, if there exists $U\in\cU$ such that $U\subset\{-1,-2,\dots\}$ (resp.\ $\{1,2,\dots\}$) and \emph{stable} otherwise.
\end{definition}
Indeed, it is not hard to check that FA-1f, East and FA-2f have zero, one and two stable directions respectively. The stability of a direction governs whether empty sites can reproduce in that direction.

The following result is the rough universality classification in one dimension due to Mar\^ech\'e, Martinelli, Morris and Toninelli~\cites{Mareche20combi,Mareche20Duarte,Martinelli19a}.

\begin{theorem}[One-dimensional rough universality]
\label{th:univ:1d}
For a one-dimensional KCM with update family $\cU$ we have that
\begin{itemize}
    \item if $\cU$ has two unstable directions, then $\qc=0$ and, for some $C>0$,
    \[\lim_{q\to 0}\bbP_\mu\left(1/C\le \frac{\log\tau_0}{\log(1/q)}\le C\right)=1;\]
    \item if $\cU$ has one unstable direction, then $\qc=0$ and, for some $C>0$,
    \[\lim_{q\to 0}\bbP_\mu\left(1/C\le \frac{\log\tau_0}{\log^2(1/q)}\le C\right)=1;\]
    \item if $\cU$ has no unstable direction, then $\qc=1$.
\end{itemize}
\end{theorem}
\begin{remark}
    The asymptotics in Theorem~\ref{th:univ:1d} also hold for $\trel$ instead of $\tau_0$.
\end{remark}
The proof of this result follows the same lines as what we have already seen, so, rather than repeating it, we only provide the appropriate pointers. The lower bound for update families with two unstable directions follows from the fact that there is typically no empty site at distance much less than $1/q$ from the origin. The corresponding upper bound follows from {\eqref{eq:F0F1},} Proposition~\ref{prop:FA1f:East:generalised} and renormalisation (empty renormalised sites correspond to completely empty intervals of sites with large fixed length). The upper bound for update families with one unstable direction follows from Theorem~\ref{th:general:1d}. The corresponding lower bound is proved as in Section~\ref{subsec:East:lower}, again with a renormalisation in order to apply the combinatorial Proposition~\ref{prop:comb}, considering a renormalised site empty if it contains at least one empty site (not to be confused with the upper bound renormalisation above). Finally, the result for update families with no unstable directions is immediate, since the state of any sufficiently long interval of occupied sites can never be modified.

It would be interesting to know the sharp asymptotics of $\log\tau_0$ for general $\cU$, as in the case of FA-1f and East, but this matter is still open. It is good to note that the corresponding problem for BP is easy (see \cite{Hartarsky22phd}*{Proposition 1.3.4}).

\section{BP universality in two dimensions}
Before we move on to two-dimensional universality for KCM, we require some background on the side of BP. Since our focus is on KCM, we take these BP results for granted and refer the interested reader to \cites{Morris17a} for a detailed survey of the methods involved.

\subsection{Rough universality in BP}
\begin{figure}
\centering
\begin{subfigure}{0.3\textwidth}
\begin{center}
\begin{tikzpicture}[line cap=round,line join=round,x=0.3\textwidth,y=0.3\textwidth]
                \clip (-1.1,-1.1) rectangle (1.1,1.1);
				\draw[very thin](0,0) circle (0.3\textwidth);
			    \draw [very thick,color=red] 
			    	plot[domain=0:pi/2,variable=\t]
					({cos(\t r)},{sin(\t r)});
                \draw (0.5,0.5) node {$\infty$};
		\end{tikzpicture}
  \caption{East (Figure~\ref{fig:East})}
\end{center}
\end{subfigure}
\begin{subfigure}{0.3\textwidth}
\begin{center}
\begin{tikzpicture}[line cap=round,line join=round,x=0.3\textwidth,y=0.3\textwidth]
                \clip (-1.1,-1.1) rectangle (1.1,1.1);
				\draw[very thin](0,0) circle (0.3\textwidth);
		\end{tikzpicture}
\caption{FA-$1$f (Figure~\ref{fig:FA1f})}
\end{center}
\end{subfigure}
 \begin{subfigure}{0.3\textwidth}
\begin{center}
\begin{tikzpicture}[line cap=round,line join=round,x=0.3\textwidth,y=0.3\textwidth]
                \clip (-1.1,-1.1) rectangle (1.1,1.1);
				\draw[very thin](0,0) circle (0.3\textwidth);
			    \fill [color=red] (0,1) circle (2pt) node[anchor=north,black] {$1$};
			    \fill [color=red] (0,-1) circle (2pt) node[anchor=south,black] {$1$};
			    \fill [color=red] (1,0) circle (2pt) node[anchor=east,black] {$1$};
			    \fill [color=red] (-1,0) circle (2pt) node[anchor=west,black] {$1$};
		\end{tikzpicture}
\caption{FA-$2$f (Figure~\ref{fig:FA2f})}
\end{center}
 \end{subfigure}
 \begin{subfigure}{0.3\textwidth}
\begin{center}
\begin{tikzpicture}[line cap=round,line join=round,x=0.3\textwidth,y=0.3\textwidth]
                \clip (-1.1,-1.1) rectangle (1.1,1.1);
				\draw[very thin](0,0) circle (0.3\textwidth);
			    \fill [color=red] (1,0) circle (2pt) node[anchor=east,black] {$1$};
			    \draw [very thick,color=red] 
			    	plot[domain=pi/2:3*pi/2,variable=\t]
					({cos(\t r)},{sin(\t r)});
                \draw (-0.7,0) node {$\infty$};
		\end{tikzpicture}
\caption{Duarte (Figure~\ref{fig:Duarte})}
\end{center}
 \end{subfigure}
\begin{subfigure}{0.3\textwidth}
\begin{center}
\begin{tikzpicture}[line cap=round,line join=round,x=0.3\textwidth,y=0.3\textwidth]
                \clip (-1.1,-1.1) rectangle (1.1,1.1);
				\draw[very thin](0,0) circle (0.3\textwidth);
			    \draw [very thick,color=red] 
			    	plot[domain=3*pi/2:3*pi,variable=\t]
					({cos(\t r)},{sin(\t r)});
                \draw (0.5,0.5) node {$\infty$};
		\end{tikzpicture}
\caption{North-East (Figure~\ref{fig:NE})}
\end{center}
 \end{subfigure}
 \begin{subfigure}{0.3\textwidth}
\begin{center}
\begin{tikzpicture}[line cap=round,line join=round,x=-0.3\textwidth,y=0.3\textwidth]
                \clip (-1.1,-1.1) rectangle (1.1,1.1);
				\draw[very thin](0,0) circle (0.3\textwidth);
			    \draw [very thick,color=red] 
			    	plot[domain=0:pi/4,variable=\t]
					({cos(\t r)},{sin(\t r)});
			    \draw [very thick,color=red] 
			    	plot[domain=pi/2:3*pi/4,variable=\t]
					({cos(\t r)},{sin(\t r)});
			    \draw [very thick,color=red] 
			    	plot[domain=pi:5*pi/4,variable=\t]
					({cos(\t r)},{sin(\t r)});
			    \draw [very thick,color=red] 
			    	plot[domain=3*pi/2:7*pi/4,variable=\t]
					({cos(\t r)},{sin(\t r)});
                \draw (0.65,0.27) node {$\infty$};
                \draw (-0.27,0.65) node {$\infty$};
                \draw (-0.65,-0.27) node {$\infty$};
                \draw (0.27,-0.65) node {$\infty$};
		\end{tikzpicture}
\caption{Spiral (Figure~\ref{fig:spiral})}
\end{center}
 \end{subfigure}
\caption{Stable directions of the two-dimensional update families from Figure~\ref{fig:rules}.\label{fig:stable}}
\end{figure}
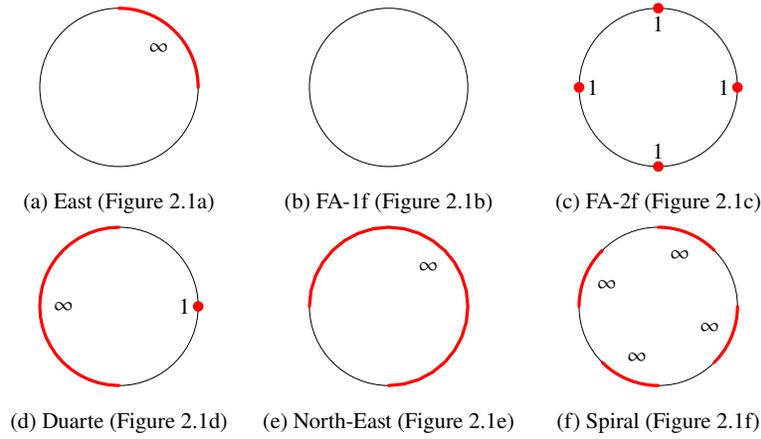

We start by generalising the definition of stable directions. We use $\|\cdot\|$ and $\<\cdot,\cdot\>$ to denote the Euclidean norm and scalar product. Let $S^{1}=\{\bu\in\bbR^2:\|\bu\|=1\}$ be the unit circle. We call the elements of $S^{1}$ \emph{directions}. We consider the open half-plane with outer normal $\bu\in S^{1}$
\begin{equation}\label{eq:openhalf}\bbH_\bu=\left\{\bx\in\bbR^2:\<\bx,\bu\><0\right\}.\end{equation}
\begin{definition}[Stable directions]
Fix an update family $\cU$. We say that $\bu\in S^{1}$ is \emph{unstable}, if there exists $U\in\cU$ such that $U\subset \bbH_\bu$ and \emph{stable} otherwise. A direction is called \emph{strongly stable}, if it belongs to the interior of the set of stable directions. {A stable direction $\bu\in S^1$ is \emph{isolated stable}, if there exists an open interval $I\subset S^1$ such that the only stable direction in $I$ is $\bu$.}
\end{definition}
The stable directions of the models in Figure~\ref{fig:rules} are given in Figure~\ref{fig:stable}. Stable directions allow us to define the rough universality classes in two dimensions.

\begin{definition}[Rough universality partition]
\label{def:2d:rough}
Let $\cC=\{\bbH_\bu\cap S^1:\bu\in S^1\}$ denote the set of open semicircles of $S^1$. An update family $\cU$ is:
\begin{itemize}
\item \emph{supercritical} if there exists $C\in\cC$ containing no stable direction. If additionally
\begin{itemize}
\item there exist two non-opposite stable directions, $\cU$ is \emph{rooted};
\item there do not exist two non-opposite stable directions, $\cU$ is \emph{unrooted}.
\end{itemize}
\item \emph{critical} if every $C\in\cC$ contains a stable direction and there exists $C\in\cC$ containing finitely many stable directions.\\
\item \emph{subcritical} if every $C\in\cC$ contains infinitely many stable directions. If additionally
\begin{itemize}
\item there exists an unstable direction, $\cU$ is \emph{nontrivial};
\item all directions are stable, $\cU$ is \emph{trivial}.
\end{itemize}
\end{itemize}
\end{definition}
Comparing Definition~\ref{def:2d:rough} with the one-dimensional case of Theorem~\ref{th:univ:1d}, we see that the new rough universality classes in two dimensions are the critical and subcritical nontrivial ones. The following rough universality theorem for BP is due to Balister, Bollob\'as, Przykucki, Smith and Uzzell \cites{Balister16,Bollobas15}.

\begin{theorem}[Two-dimensional rough universality for BP]
\label{th:2d:univ:rough}
Let $\cU$ be a two-dimensional update family. If $\cU$ is
\begin{itemize}
\item supercritical, then $\qc=0$ and, for some $C>0$,
\[\lim_{q\to 0}\mu\left(1/C\le \frac{\log\tbp}{\log(1/q)}\le C\right)=1;\]
\item critical, then $\qc=0$ and, for some $C>0$,
\[\lim_{q\to 0}\mu\left(1/C\le \frac{\log\log\tbp}{\log(1/q)}\le C\right)=1;\]
\item subcritical nontrivial, then $0<\qc<1$.
\item subcritical trivial, then $\qc=1$.
\end{itemize}
\end{theorem}

We note that a complete generalisation of Theorem~\ref{th:2d:univ:rough} has already been established in arbitrary dimension by Balister, Bollob\'as, Morris and Smith \cites{Balister22,Balister24,BalisterNaNb}, but its KCM counterpart is still missing. The reader interested in subcritical models is also encouraged to consult \cites{Hartarsky22Toom,Toom80,Gray99,Hartarsky22sharpness,Hartarsky21} for different approaches to bounding $\qc$. While the behaviour of nontrivial subcritical models is very interesting, it is quite challenging and not much is known currently, so we discard them in the sequel. We further refer to \cite{Hartarsky22phd}*{Section~1.5.1} for a detailed account of the history of universality in BP and KCM.

\subsection{Refined universality in BP}
The bounds on $\tbp$ for critical models in Theorem~\ref{th:2d:univ:rough} are quite loose due to the iterated logarithm, particularly compared to results for specific models (recall Theorem~\ref{th:BP}). In two dimensions, it is possible to obtain much more precise asymptotics. In order to state them, we require a refinement of the notion of stable direction (see Figure~\ref{fig:stable} for examples). Recall that $|\cdot|$ is the number of empty sites and $[\cdot]$ is the closure from \eqref{eq:def:closure}.
\begin{definition}[Difficulty]
    \label{def:alpha}
The \emph{difficulty} $\alpha(\bu)$ of $\bu\in S^1$ is 
\begin{itemize}
    \item $0$ if $\bu$ is unstable;
    \item $\infty$ if $\bu$ is stable, but not isolated stable;
    \item $\min\{|Z|:Z\subset\bbZ^2,|[\bzero_{\bbH_\bu\cup Z}\cdot\bone_{\bbZ^2\setminus(\bbH_\bu\cup Z)}]_{\bbZ^2\setminus\bbH_\bu}|=\infty\}$ otherwise.
\end{itemize}
The \emph{difficulty} of $\cU$ is
\[\alpha=\alpha(\cU)=\min_{C\in\cC}\max_{\bu\in C}\alpha(\bu).\]
We say that a direction $\bu\in S^1$ is \emph{hard} if $\alpha(\bu)>\alpha$. We say that $\cU$ is \emph{unbalanced} if there exist two opposite hard directions and \emph{balanced} otherwise.
\end{definition}
In words, the difficulty of an isolated stable direction is the smallest number of empty sites needed to empty an infinite number of sites with the help of an empty half-plane with outer normal the direction. The difficulty of the update family is given by the easiest open semi-circle, a semi-circle being as hard as the hardest direction it contains. Comparing Definitions~\ref{def:2d:rough} and~\ref{def:alpha}, it can be shown \cite{Bollobas23} that an update family is supercritical if $\alpha=0$, subcritical if $\alpha=\infty$ and critical if $0<\alpha<\infty$. Among the critical examples of Figure~\ref{fig:rules}, FA-2f is balanced and Duarte is unbalanced, both having difficulty $1$ (see Figure~\ref{fig:stable}).

With these notions, the refined universality result for BP of Bollob\'as, Duminil-Copin, Morris and Smith \cite{Bollobas23} is the following.
\begin{theorem}[Two-dimensional refined universality for BP]
\label{th:2d:univ:refined:BP}
    Let $\cU$ be a two-dimensional critical update family of difficulty $\alpha$. Then, for some $C>0$,
    \[\lim_{q\to0}\mu\left(1/C\le \frac{q^\alpha\log\tbp}{(\log(1/q))^\gamma}\le C\right)=1,\]
    where $\gamma=0$, if $\cU$ is balanced, and $\gamma=2$, if $\cU$ is unbalanced.
\end{theorem}
Naturally, Theorem~\ref{th:2d:univ:refined:BP} is consistent with Theorem~\ref{th:BP} for 2-neighbour and Duarte BP.

\section{KCM universality in two dimensions}
\subsection{Statement}
With BP universality at hand, we may turn to KCM. For subcritical models, in view of Theorems~\ref{th:ergo} and~\ref{th:exponential} combined with the BP rough universality Theorem~\ref{th:2d:univ:rough}, we do not have anything new to say. We therefore focus on supercritical and critical models in the next result due to Mar\^ech\'e, Martinelli, Morris and Toninelli \cites{Mareche20combi,Mareche20Duarte,Martinelli19a}. Recall Definition~\ref{def:2d:rough}.
\begin{theorem}[Two-dimensional rough universality for KCM]
\label{th:universality:rough:KCM:2d}For any two-dimensional KCM with update family $\cU$ we have that
\begin{itemize}
    \item if $\cU$ is supercritical unrooted, then for some $C>0$,
    \begin{equation}
    \label{eq:rough:supercrit:unrooted}\lim_{q\to0}\bbP_\mu\left(1/C\le \frac{\log \tau_0}{\log(1/q)}\le C\right)=1;\end{equation}
    \item if $\cU$ is supercritical rooted, then for some $C>0$, 
    \begin{equation}
    \label{eq:rough:supercrit:rooted}\lim_{q\to0}\bbP_\mu\left(1/C\le \frac{\log \tau_0}{\log^2(1/q)}\le C\right)=1;\end{equation}
    \item if $\cU$ is critical, then for some $C>0$,
    \begin{equation}
    \label{eq:rough:critical}\lim_{q\to0}\bbP_\mu\left(1/C\le \frac{\log\log \tau_0}{\log(1/q)}\le C\right)=1.\end{equation}
\end{itemize}
The same asymptotics hold for $\trel$ instead of $\tau_0$.    
\end{theorem}
The proof of Theorem~\ref{th:universality:rough:KCM:2d} will be explained in Section~\ref{subsec:rough}, but before that, let us first discuss refined universality. Recalling Definition~\ref{def:alpha}, we only need the following vocabulary in order to define the refined KCM universality classes.
\begin{definition}[Further refined universality types]
\label{def:2d:refined}
A critical two-dimensional update family $\cU$ is \emph{rooted}, if there exist two non-opposite hard directions and \emph{unrooted} otherwise. We say that $\cU$ is \emph{semi-directed}, if there is exactly one hard direction and \emph{isotropic} if there are no hard directions.    
\end{definition}
Notice that balanced unrooted update families are either semi-directed or isotropic. The above notions allow us to state the refined universality result obtained over a series of works of Hartarsky, Mar\^ech\'e, Martinelli, Morris and Toninelli \cites{Hartarsky22univlower,Hartarsky24univupper,Martinelli19a,Hartarsky21a,Hartarsky20} (see \cite{Hartarsky24univupper} for the final step and a detailed discussion), also relying on the BP refined universality Theorem~\ref{th:2d:univ:refined:BP} of \cite{Bollobas23}.
\begin{theorem}[Two-dimensional refined universality for KCM]
\label{th:2d:univ:refined:KCM}
    Let $\cU$ be a two-dimensional critical update family of difficulty $\alpha$. Then for some $C>0$,
    \[\lim_{q\to0}\bbP_\mu\left(1/C\le \frac{q^{\alpha\cdot\beta}\log\tau_0}{(\log(1/q))^\gamma(\log\log(1/q))^\delta}\le C\right)=1,\]
    where the exponents $\beta$, $\gamma$ and $\delta$ are given in the following table depending on whether $\cU$ has finite or infinite number of stable directions; is balanced or unbalanced; is rooted or unrooted; is semi-directed or isotropic, if balanced and unrooted.
    \begin{center}
        \begin{tabular}{c| >{\centering\arraybackslash}m{4cm}| >{\centering\arraybackslash}m{2cm}| >{\centering\arraybackslash}m{2cm}}
         \multirow{2}{*}{$\beta, \gamma, \delta$}&Infinite stable directions
         & \multicolumn{2}{c}{Finite stable directions}\\\cline{2-4}
         & Rooted & Rooted & Unrooted\\
        \hline
       Unbalanced & $2,4,0$ & $1,3,0$ & $1,2,0$\\
        \hline
        \multirow{2}{*}{Balanced} & \multirow{2}{*}{$2,0,0$} & \multirow{2}{*}{$1,1,0$} & S.-dir.~$1,0,1$\\\cline{4-4}
        &&&{Iso.~$1,0,0$}\\
    \end{tabular}\end{center}
    In particular, $\beta=2$, if $\cU$ has a strongly stable direction, and $\beta=1$ otherwise.
\end{theorem}
It is not hard to check that the seven refined universality classes in Theorem~\ref{th:2d:univ:refined:KCM} do exhaust the critical rough universality class of Theorem~\ref{th:universality:rough:KCM:2d}. We emphasise that Theorem~\ref{th:2d:univ:refined:KCM} is currently the sharpest result available for any critical KCM with the exception of FA-2f and slight variations thereof, namely modified FA-2f and the Frob\"ose KCM (recall Section~\ref{subsubsec:matryoshka} and see \cites{Frobose89,Gravner09}). As in Theorem~\ref{th:FA2f:sharp:threshold}, the upper bounds in Theorem~\ref{th:2d:univ:refined:KCM} are not known to hold for $\trel$.

Interestingly, an analogue of Theorem~\ref{th:2d:univ:refined:KCM} is not available in one dimension. It would be good to fill this gap by solving the following problem.
\begin{problem}[{Unrooted scaling in one dimension}]
\label{prob:1d:refined:univ}
    Let $\cU$ be a one-dimensional update family with two unstable directions. Prove existence of and determine $\alpha\in(0,\infty)$ such that \[\lim_{C\to\infty}\liminf_{q\to0}\bbP_\mu\left(1/C\le q^\alpha\tau_0\le C\right)=1.\]
\end{problem}
It should be noted that, while the analogous problem for BP is an exercise \cite{Hartarsky22phd}*{Proposition~1.3.4}, in view of Theorem~\ref{th:FA1f}, Problem~\ref{prob:1d:refined:univ} is not.

\subsection{Rough universality proofs}
\label{subsec:rough}
We next outline the proof of Theorem~\ref{th:universality:rough:KCM:2d}. Firstly, the lower bound of \eqref{eq:rough:supercrit:unrooted} is immediate, since there is typically no empty site at distance much smaller than $q^{-1/2}$ from the origin. The lower bound in \eqref{eq:rough:critical} follows from Theorem~\ref{th:2d:univ:rough} together with \eqref{eq:t0t0bp} and \eqref{eq:tau0}. The upper bounds in \eqref{eq:rough:supercrit:unrooted} and \eqref{eq:rough:supercrit:rooted} are proved using Proposition~\ref{prop:FA1f:East:generalised} and \eqref{eq:F0F1} together with a simple renormalisation, as for Theorem~\ref{th:univ:1d}. Namely, each empty site corresponds to a suitably oriented rectangle whose sites are all empty. The exact shape of this rectangle is chosen based on the proof of BP rough universality, so that, if it is empty, it is able to reproduce an empty copy of itself (see \cites{Bollobas15,Martinelli19a}). The upper bound in \eqref{eq:rough:critical} proceeds like the proof of Theorem~\ref{th:FA2f:MT}, but using a generalised one-dimensional KCM with East instead of FA-1f constraint, still covered by Proposition~\ref{prop:FA1f:East:generalised}.

We are left with the lower bound in \eqref{eq:rough:supercrit:rooted} which is the only one requiring additional ideas with respect to what we have already seen. The overall scheme remains the same---we seek to establish a combinatorial bottleneck akin to the one of Proposition~\ref{prop:comb} and convert it into the desired bound as in Section~\ref{subsubsec:bottleneck:deduction}. However, contrary to the one-dimensional case of Theorem~\ref{th:univ:1d}, the combinatorial bottleneck cannot be deduced from the one-dimensional Proposition~\ref{prop:comb} via renormalisation. We next present a sketch of the proof of the following result of Mar\^ech\'e \cite{Mareche20combi}*{Theorem~4}, which is the crucial ingredient.

\begin{proposition}[Combinatorial bottleneck for rooted models]
    \label{prop:Laure:combi}Let $\cU$ be a two-dimensional update family which is not supercritical unrooted. There exists an integer $C>0$ such that the following holds for any integer $n\ge 1$. Consider the $\cU$-KCM on $\Lambda_n=\{-Cn2^n,\dots,Cn2^n\}^2$ with boundary condition $\bzero_{\bbZ^2\setminus \Lambda}$. Let $V(n)$ be the set of all configurations that the process can reach from $\bone_\Lambda$ via a legal path (recall Definition~\ref{def:legal:path}) in which all configurations contain at most $n$ empty sites. Then $\omega_0=1$ for all $\omega\in V(n)$.
\end{proposition}
    \begin{figure}
        \centering
        \begin{tikzpicture}[x=0.2cm,y=0.2cm]
            \draw (-9,-9) rectangle (9,9);
            \draw (-3,-3) rectangle (3,3);
            \draw [fill=black, fill opacity=0.5,even odd rule] (-5,-5) rectangle (5,5) (-7,-7) rectangle (7,7);
            \draw (0,0) node{$\Lambda_{n-1}$};
            \draw (9,0) node[right]{$\Lambda_{n}$};
        \end{tikzpicture}
        \caption{Illustration of the proof of Proposition~\ref{prop:Laure:combi}. The buffer zone is shaded.}
        \label{fig:Laure}
    \end{figure}
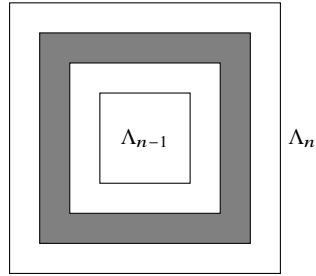
\begin{proof}[Sketch]
    For simplicity, we focus on the two-dimensional East KCM (see Figure~\ref{fig:East}). The proof proceeds by induction on $n$ claiming that for any $\omega\in V(n)$, $|\omega_{\Lambda_{n-1}}|\le n-1$. Indeed, iterating this fact, we obtain that for any $\omega\in V(n)$, $|\omega_{\Lambda_1}|\le 1$, so that the single empty site in $\Lambda_1$ cannot be at distance more than $C$ from the boundary of $\Lambda_1$, so it cannot reach the origin.

    In order to show the claim, using the reversibility of legal paths (recall Definition~\ref{def:legal:path}), we may instead prove that, if we start from $\omega\neq\bone_{\Lambda_n}$ such that $\omega_{\Lambda_n\setminus\Lambda_{n-1}}=\bone$ and never visit configurations with more than $n$ empty sites, then we cannot reach the $\bone_{\Lambda_n}$ configuration. To do this, we prove by a second induction, on the number of steps in the legal path, that the following two conditions remain valid. Firstly, a frame-shaped buffer zone around $\Lambda_{n-1}$ with no empty site remains intact (see Figure~\ref{fig:Laure}). Secondly, there always remains an empty site in the internal region encircled by the buffer, so the dynamics cannot reach $\bone_{\Lambda_n}$.

    We know that so far an empty site remains trapped in the internal region encircled by the buffer, so we only have $n-1$ empty sites available for disrupting the buffer from the outside, which is impossible by the induction hypothesis on $n$. Therefore, it suffices to show that we may not disrupt the buffer from the inside either. By projecting the two-dimensional East model on each axis it is clear that no empty site can enter the right and top parts of the buffer from the inside, and the projections of the top and rightmost empty sites in the region inside the buffer need to remain where they were initially. The left part of the buffer (and similarly for the bottom one) cannot be reached from the inside, because at least one empty site needs to remain as far right as the rightmost initial one was, so we only have $n-1$ empty sites with which to reach the left part of the buffer, which is impossible by induction hypothesis on $n$.
\end{proof}
With the above sketch in mind, we encourage the reader to consult the full proof in \cite{Mareche20combi}, which actually takes only four pages. We also remark that Proposition~\ref{prop:Laure:combi} and its proof generalise immediately to any dimension.

\subsection{Refined universality proofs}

\subsubsection{Lower bounds}
\label{subsec:refined:lower}
\begin{figure}
\centering
\begin{subfigure}{0.49\textwidth}
    \begin{tikzpicture}[x=0.5cm,y=0.5cm]
\draw (-1.5,1.5) circle (0.15);
\draw (-2,0.5) circle (0.15);
\draw (-1.5,1) circle (0.15);

\draw (1,1) circle (0.15);
\draw (1.5,0.5) circle (0.15);

\draw (-2,-2) circle (0.15);
\draw (-1.5,-2.5) circle (0.15);
\draw (-2.5,-2) circle (0.15);
\draw (-1.5,-3) circle (0.15);

\draw (1.5,-2) circle (0.15);
\draw (1,-2.5) circle (0.15);
\draw (2,-2) circle (0.15);
\draw[step=0.5,gray,very thin](-2.5,-3)grid(-0.5,-1);
\draw[step=0.5,gray,very thin](0,-3)grid(2,-1);
 \draw[step=0.5,gray,very thin](-2.5,-0.5)grid(-0.5,1.5);
\draw[step=0.5,gray,very thin](0,-0.5)grid(2,1.5);
\draw (-1.5,-2) node[cross=2pt,rotate=0] {};
\draw (1,-2) node[cross=2pt,rotate=0] {};
\draw (-1.5,0.5) node[cross=2pt,rotate=0] {};
\draw (1,0.5) node[cross=2pt,rotate=0] {};
		\end{tikzpicture}\quad
		\begin{tikzpicture}[x=1cm,y=1cm]
				\clip (-1.2,-1.2) rectangle (1.2,1.2);
				\draw[very thin](0,0) circle (1);
			    \fill [color=red] (0,1) circle (2pt) node[anchor=north,black] {$2$};
			    \fill [color=red] (1,0) circle (2pt) node[anchor=east,black] {$2$};
				\fill [color=red] (-1,0) circle (2pt) node[anchor=west,black] {$1$};
				\fill [color=red] (0,-1) circle (2pt) node[anchor=south,black] {$1$};
		\end{tikzpicture}
    \caption{Balanced model of Section~\ref{subsec:refined:lower}.\label{fig:log1}}
\end{subfigure}
\begin{subfigure}{0.49\textwidth}

    \begin{tikzpicture}[x=0.5cm,y=0.5cm]
\draw (-1.5,1.5) circle (0.15);
\draw (-2,0.5) circle (0.15);
\draw (-2.5,0.5) circle (0.15);
\draw (-1.5,1) circle (0.15);

\draw (1,1) circle (0.15);
\draw (1.5,0.5) circle (0.15);
\draw (2,0.5) circle (0.15);

\draw (-2,-2) circle (0.15);
\draw (-1.5,-2.5) circle (0.15);
\draw (-2.5,-2) circle (0.15);
\draw (-1.5,-3) circle (0.15);

\draw (1.5,-2) circle (0.15);
\draw (1,-2.5) circle (0.15);
\draw (2,-2) circle (0.15);
\draw[step=0.5,gray,very thin](-2.5,-3)grid(-0.5,-1);
\draw[step=0.5,gray,very thin](0,-3)grid(2,-1);
 \draw[step=0.5,gray,very thin](-2.5,-0.5)grid(-0.5,1.5);
\draw[step=0.5,gray,very thin](0,-0.5)grid(2,1.5);
\draw (-1.5,-2) node[cross=2pt,rotate=0] {};
\draw (1,-2) node[cross=2pt,rotate=0] {};
\draw (-1.5,0.5) node[cross=2pt,rotate=0] {};
\draw (1,0.5) node[cross=2pt,rotate=0] {};
		\end{tikzpicture}\quad
		\begin{tikzpicture}[x=1cm,y=1cm]
				\clip (-1.2,-1.2) rectangle (1.2,1.2);
				\draw[very thin](0,0) circle (1);
			    \fill [color=red] (0,1) circle (2pt) node[anchor=north,black] {$2$};
			    \fill [color=red] (1,0) circle (2pt) node[anchor=east,black] {$2$};
				\fill [color=red] (-1,0) circle (2pt) node[anchor=west,black] {$1$};
				\fill [color=red] (0,-1) circle (2pt) node[anchor=south,black] {$2$};
		\end{tikzpicture}
    \caption{Unbalanced model of Section~\ref{subsec:refined:upper}.\label{fig:log3}}
\end{subfigure}
\caption{The update rules, stable directions and difficulties of two example critical rooted update families with difficulty $\alpha=1$.}
\end{figure}
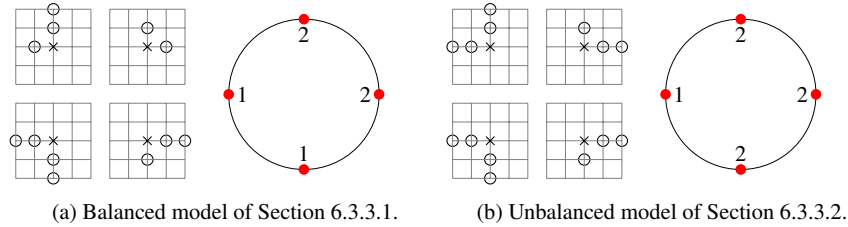

We start with the lower bounds in Theorem~\ref{th:2d:univ:refined:KCM}, since they are closely related to Proposition~\ref{prop:Laure:combi}. We follow \cite{Hartarsky22univlower} and refer to that work for the formal proof. Somewhat surprisingly, all seven refined universality classes are governed by the same combinatorial bottleneck, but on different length scales and for different reasons. For the sake of concreteness, we focus on the model from Figure~\ref{fig:log1}, whose difficulty is $\alpha=1$.

Morally speaking, in this model the smallest mobile entity (`droplet') is an empty square of size roughly $1/q$, similarly to FA-2f (recall Chapter~\ref{chap:FA2f}). Indeed, typically, on its left and bottom sides, one can find an empty site, which allows it to empty the column of sites on its left and the row of sites below it. However, it is essentially impossible for the droplet to grow up or right, as this requires two consecutive empty sites and those are typically only available at distance $1/q^2$ from the droplet. We will only work in the box $\Lambda=\{-1/q^{7/4},\dots,1/q^{7/4}\}^2$, so such couples of empty sites are not available for most columns and rows. Thus, we expect droplets to essentially follow the dynamics of the two-dimensional East KCM (see Figure~\ref{fig:East}).

On a very high level, we proceed in the same way as in the proof of Proposition~\ref{prop:Laure:combi}. However, there are several glaring problems in making the above reasoning rigorous. Firstly, much like in Figure~\ref{fig:FA2f:droplets}, droplets can be more complex than empty squares of size $1/q$. Thus, one needs to identify an event which says whether a droplet is present and this event should be deterministically necessary for empty sites to  spread. Moreover, the event should have probability of the correct order $\exp(-1/q)$, as suggested by the BP Theorem~\ref{th:2d:univ:refined:BP}. It turns out that the notion of `spanning' introduced in \cite{Bollobas23} for the proof of the lower bound of Theorem~\ref{th:2d:univ:refined:BP} for unbalanced models, following \cite{Cerf99}, is flexible enough for our purposes. Roughly speaking, a droplet (rectangle) is spanned if the empty sites present inside it are sufficient to empty a connected set touching all its sides (see \cite{Morris17a}*{Section~8} for an overview of the BP result and its proof). We call a droplet critical if it has size roughly $1/q$. 

If our goal were to only find a single critical droplet, we could proceed as in Section~\ref{sec:FA2f:lower}. Indeed, for isotropic and unbalanced unrooted models the proof of Section~\ref{sec:FA2f:lower} or directly \eqref{eq:t0t0bp} combined with Theorem~\ref{th:BP} suffices. However, for other refined universality classes, including the one of the model of Figure~\ref{fig:log1}, we need to work much harder. In order to obtain the exponent $\gamma=1>0$ in Theorem~\ref{th:2d:univ:refined:KCM}, we need many droplets. Unfortunately, given a configuration, spanned critical droplets may overlap, so, in order to obtain good bounds on the probability of the configuration, we need to consider disjointly occurring ones (droplets occur disjointly, if they admit disjoint witness sets of empty sites). We may then define the number of spanned critical droplets as the maximal number of disjointly occurring ones. Considering the KCM on $\Lambda$ with $\bzero_{\bbZ^2\setminus\Lambda}$ boundary condition, starting from a typical initial configuration (which in particular does not contain more than three empty sites close to each other), our aim is to prove that, before we can empty the origin, we need to visit a configuration with at least order $\log(1/q)$ disjointly occurring critical droplets.

If droplets did follow a two-dimensional East dynamics exactly (up to renormalisation), this would follow from Proposition~\ref{prop:Laure:combi}. But this is not the case. Indeed, by changing their internal structure, droplets may move a bit without creating another droplet, as we saw in Section~\ref{sec:FA2f:upper}. Worse yet, they are not really forbidden to move right or up, but simply are not likely to be able to do so wherever they want: it depends on the dynamic environment.

In order to handle these problems, we need the crucial notion of crossing. Consider a vertical strip $S$ of width $1/q^{3/2}$ of our domain $\Lambda$. Roughly speaking, we say that $S$ has a crossing if the following two events occur. Firstly, the empty sites in $S$ together with the entire half-plane to the left of $S$ are enough to infect a path from left to right in $S$. Secondly, $S$ does not contain a spanned critical droplet. Notice that these two events have opposite monotonicity in the configuration. Employing BP tools, it can be shown (see \cite{Hartarsky22univlower}*{Appendix B}) that the probability of a crossing decays exponentially with the width of $S$ at our scales of interest. In particular, the probability under $\mu$ that such a strip $S$ is crossed is of order $\exp(-1/q^{3/2})$. Roughly speaking, the proof proceeds by splitting $S$ into smaller strips which are either crossed by a single spanned droplet of subcritical size, or contain a pair of adjacent empty sites. Having appropriate bounds on the probability of spanned subcritical droplets as a function of their size, one may prove the desired bound on crossing by a union bound over the partition of $S$ into the smaller strips. We note that in the case of infinite number of stable directions, bounding the probability of crossings is quite different \cite{Hartarsky22univlower}*{Appendix~B.2}, but remains possible.

Having established such a bound on the probability of crossings, we may incorporate them into the combinatorial bottleneck---we are satisfied if we visit either a configuration with a crossed strip of width $1/q^{3/2}$, or with $\log(1/q)$ disjointly spanned critical droplets. The lack of crossings allows us to exclude the possibility of a droplet reaching the right side of the vertical strip $S$ without help from the right of $S$, since the KCM dynamics can never infect more than what bootstrap percolation can (recall Section~\ref{sec:legal:paths}).

With these additional inputs, the proof scheme of Proposition~\ref{prop:Laure:combi} can be carried out to give the lower bounds in Theorem~\ref{th:2d:univ:refined:KCM}. The only difference between refined universality classes comes in the choice of length scales, in the bounds on the probability of spanning droplets and crossing strips and their proofs.

\subsubsection{Upper bounds}
\label{subsec:refined:upper}
We next turn to the upper bounds in Theorem~\ref{th:2d:univ:refined:KCM}. Contrary to the lower bounds, the proofs of upper bounds are highly dependent on the refined universality class. However, there are two classes, for which all the elements of the proof have already been discussed. The weakest upper bound corresponding to $\beta=2,\gamma=4,\delta=0$ in Theorem~\ref{th:2d:univ:refined:KCM}, which applies to all models, but is only sharp for unbalanced critical families with infinite number of stable directions, in fact follows from the proof of the rough Theorem~\ref{th:universality:rough:KCM:2d} mentioned in Section~\ref{subsec:rough} (see \cite{Martinelli19a}). At the other extreme, the upper bound for isotropic models ($\beta=1,\gamma=0,\delta=0$) is proved similarly to Theorem~\ref{th:FA2f:sharp:threshold}, using the matryoshka doll technique (recall Section~\ref{sec:FA2f:upper}), up to some technical modifications (see \cite{Hartarsky24univupper}*{Section~5}).

Our next goal is to outline the proof of the upper bound of Theorem~\ref{th:2d:univ:refined:KCM} with $\beta=1,\gamma=3,\delta=0$, which applies to any critical KCM with finite number of stable directions, but is only sharp for unbalanced rooted families with finite number of stable directions. This will clarify the relevance of the absence of any strongly stable directions, which governs the value of the most important exponent $\beta$. For the sake of simplicity, we focus on the model depicted in Figure~\ref{fig:log3}. 

Let us start with some heuristic considerations before explaining how they can be turned into the proof originating from \cite{Hartarsky21a}. Since we are interested in an upper bound, we may choose the notion of droplet in a simple way. Namely, a droplet $D$ is (a translate of) an empty square frame of size $C\log(1/q)/q$ and thickness $2$, where $C$ is a suitably large constant (see Figure~\ref{fig:mechanism:1}). Then typically (under $\mu$), there is an empty site in the column to the left of $D$ allowing us to empty $D-(1,0)$. However, it is unlikely to find a pair of adjacent empty sites on any of the other sides of $D$. We conclude that it is easy for $D$ to advance only to the left. An efficient way to perform this leftward motion is given by the legal path for the East KCM from Figure~\ref{fig:East_path}, where each empty site represents an empty translate of $D$.

\begin{figure}[b]
  \begin{subfigure}{\linewidth}
  \centering
\begin{tikzpicture}[x=-0.25cm,y=0.25cm,scale=0.5]
\clip(-15,-5) rectangle (55,11);
\fill[color=gray] (10,0) rectangle (40,10);
  \fill (0,0) -- (10,0) -- (10,10) -- (0,10) -- cycle;
\fill[line width=0pt,color=white,fill=white] (2,2) rectangle (8,8);
\fill[fill=black] (40,10) -- (30,10) -- (30,0) -- (40,0) -- cycle;
\fill[line width=0pt,color=gray,fill=gray] (32,2) rectangle (38,8);
\draw (39.5,10.5) circle (0.3);
\draw (38.5,10.5) circle (0.3);
\draw [thick,->] (12,5)--(27,5);
\draw [<->] (-2,0) -- (-2,10) node[midway, right]{$\frac{C\log(1/q)}{q}$};
\draw [<->] (0,-1) -- (40,-1) node[midway, below]{$1/q^2$};
\end{tikzpicture}
\end{subfigure}
\begin{subfigure}{\linewidth}
\centering
\begin{tikzpicture}[x=-0.25cm,y=0.25cm,scale=0.5]
  \clip(-12,-0.5) rectangle (52,11);
\fill[color=gray] (10,0) rectangle (40,11);
\fill (0,0) rectangle (10,11);
\fill[line width=0pt,color=white,fill=white] (2,3) rectangle (8,8);
\fill[fill=black] (30,0) rectangle (40,11);
\fill[line width=0pt,color=white,fill=gray] (32,3) rectangle (38,8);
\draw [thick,<-] (12,5)--(27,5);
\end{tikzpicture}
\end{subfigure}
 \caption{
\label{fig:mechanism:1}
The mechanism for the droplet to grow up.}
\end{figure}
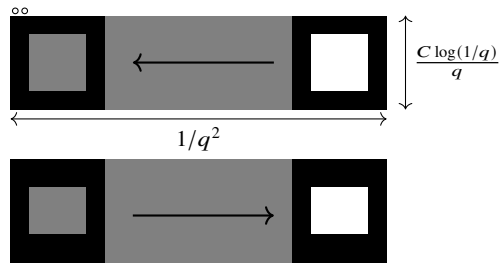
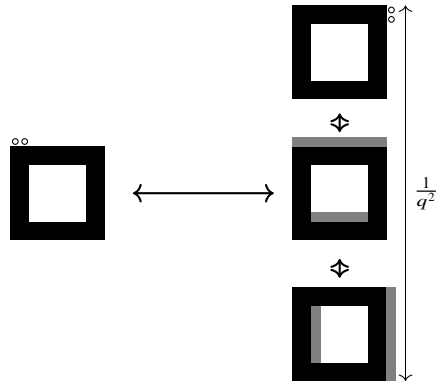
\begin{figure}[t]
 \centering
\begin{tikzpicture}[x=0.25cm,y=0.25cm,scale=0.5]
  \clip(-55,0) rectangle (15,40);
  \begin{scope}[shift={(-10,15)}]
  \fill (0,0) -- (10,0) -- (10,10) -- (0,10) -- cycle;
\fill[line width=0pt,color=white,fill=white] (2,2) rectangle (8,8);
\fill[fill=black] (-20,0) rectangle (-30,10);
\fill[line width=0pt,color=white,fill=white] (-22,2) rectangle (-28,8);
\draw (-29.5,10.5) circle (0.3);
\draw (-28.5,10.5) circle (0.3);
\draw [thick,<->] (-2,5)--(-17,5);

\fill[color=gray] (0,10) rectangle (10,11);
\fill[color=gray] (2,2) rectangle (8,3);
 \end{scope}
  
  \begin{scope}[rotate=90]
  \fill (0,0) -- (10,0) -- (10,10) -- (0,10) -- cycle;
\fill[line width=0pt,color=white,fill=white] (2,2) rectangle (8,8);
\fill[fill=black] (40,10) -- (30,10) -- (30,0) -- (40,0) -- cycle;
\fill[line width=0pt,color=white,fill=white] (32,2) rectangle (38,8);
\draw (39.5,-0.5) circle (0.3);
\draw (38.5,-0.5) circle (0.3);

\fill[color=gray] (0,0) rectangle (10,-1);
\fill[color=gray] (8,8) rectangle (2,7);
\draw [thick, <->] (11,5)--(13,5);
\draw [thick, <->] (26.5,5)--(28.5,5);
\end{scope}
\draw [<->] (2,0) -- (2,40) node[midway, right]{$\frac{1}{q^2}$};
\end{tikzpicture}
\caption{The mechanism for the droplet to grow to the right by
  making a long excursion up, each of whose steps is a long excursion to the left, 
  as in Figure~\ref{fig:mechanism:1}.}
\label{fig:mechanism:2}
\end{figure}

The key idea is that it suffices to perform this East-like motion for a distance of order $1/q^2$, in order to find a pair of adjacent empty sites on the row above the droplet. Once the droplet reaches them, it is able to move one step up. It is then possible to revert the East path to bring the droplet to the original position, but shifted one lattice step up. This procedure effectively yields a step in the hard up direction. We may then iterate this idea, moving upwards in an East-like manner, where each step up is, in fact, a long East-like path to the left and back. This way, we eventually reach a pair of adjacent empty sites allowing the droplet to move to the right, etc. See Figures~\ref{fig:mechanism:1} and~\ref{fig:mechanism:2} for an illustration of the mechanism. Based on the above heuristics, we expect droplets to be able to move freely in all directions by creating only about $\log(1/q)$ additional droplets at a time. Hence, we expect the time necessary for droplets to move to be of order $\rho^{-\log(1/q)}$, where $\rho\approx q^{8C\log(1/q)/q}$ is the probability of a single droplet under $\mu$. This gives the desired $\exp(8C\log^3(1/q)/q)$ time scale.

\begin{figure}[b]
 \centering
\begin{tikzpicture}[line cap=round,line join=round,x=-0.15cm,y=0.15cm]
\draw [thick] (0,0) rectangle (30,5) node[midway]{$B$};
\draw [thick] (0,5.3) rectangle (30,30) node[midway]{$R_0$};
\draw [thick] (-0.3,0) rectangle (-15.3,30) node[midway]{$R_1^-$};
\draw [thick] (30.3,0) rectangle (45.3,30) node[midway]{$R_1^+$};
\fill (0,0) rectangle (5,5);
\fill[line width=0pt,color=white,fill=white] (1,1) rectangle(4,4);
\draw (2.5,2.5) node{$\mathring D$};
\draw [<->] (0,29)--(30,29) node[midway,below] {$C \frac{\log(1/q)}{q^2}$};
\draw [<->] (-16,0)--(-16,30) node[midway,right] {$C \frac{\log(1/q)}{q^2}$};
\draw [<->] (-1,0)--(-1,5) node[midway,right] {$C \frac{\log(1/q)}{q}$};
\draw (2,4)--(5,7) node[above left]{$D$};
\end{tikzpicture}
\caption{Geometry of the matryoshka dolls in Section~\ref{subsec:refined:upper} for the update family depicted in Figure~\ref{fig:log3}.}
\label{fig:mechanism:3}
\end{figure}
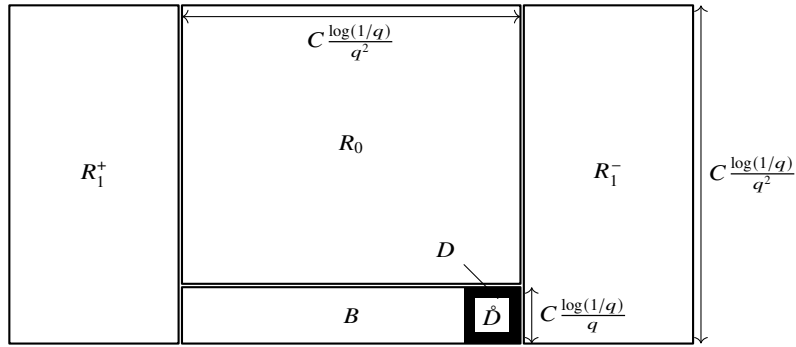

In order to turn the above into a proof, we use the matryoshka doll technique from Section~\ref{subsubsec:matryoshka}. The geometry of the consecutive regions is given in Figure~\ref{fig:mechanism:3}. The super good event $\cS\cG(\Lambda)$ for $\Lambda=D\cup \mathring D\cup B\cup R_0\cup R_1^+\cup R_1^-$ requires that:
\begin{itemize}
    \item the droplet $D$ is empty;
    \item each column of the \emph{base} $B$ contains an empty site
    \item each row of the rectangle $R_0$ contains a pair of adjacent empty sites;
    \item each column of the rectangles $R_1^+$ and $R_1^-$ contains a pair of adjacent empty sites.
\end{itemize}
Thanks to our choice of geometry, it is not hard to check that the latter three events are very likely under $\mu$, so $\mu(\cS\cG(\Lambda))\approx \rho=q^{|D|}\approx \exp(8C(\log(1/q))^2/q)$. Recalling Section~\ref{sec:finite_vol}, we then seek to prove the Poincar\'e inequality
\begin{align}
\label{eq:snail:Poincare}\mu_{\Lambda}(f|\cS\cG(\Lambda))&{}\le \gamma(\Lambda)\cD_{\bone_{\bbZ^2\setminus\Lambda}}(f),&
\gamma(\Lambda)&{}\le e^{(C\log(1/q))^3/q}\end{align}
for any local function $f:\Omega\to\bbR$. Once \eqref{eq:snail:Poincare} is proved, concluding the proof of Theorem~\ref{th:2d:univ:refined:KCM} for the model under consideration can be done along the lines of Section~\ref{subsec:FA2f:reduction:to:meso}.

The proof of \eqref{eq:snail:Poincare} proceeds in a roughly similar way to \eqref{eq:FA2f:meso}, by proving Poincar\'e inequalities successively for $D$, $D\cup\mathring D$, $D\cup\mathring D\cup B$,\dots, $\Lambda$. The one for $D$ is trivial, since $\cS\cG(D)=\{\bzero_D\}\times\Omega_{\bbZ^2\setminus D}$ is a single configuration if restricted to $D$. The inequality for $D\cup \mathring D$ is proved by dividing $\mathring D$ into vertical strips of width 2 and using Proposition~\ref{prop:FA1f:East:generalised} for generalised FA-1f in one dimension. This yields the relaxation time bound
\[\gamma\left(D\cup\mathring D\right)\le q^{C^2\log(1/q)/q},\]
where $\gamma(D\cup\mathring D)$ is defined as in \eqref{eq:snail:Poincare}.

Proving that $\gamma(D\cup \mathring D\cup B)\le \exp(e^{C^2(\log(1/q))^3/q})$ (and similarly for adding the remaining rectangles one by one) is done along the lines of the proof of \eqref{eq:FA2f:recursion}. Namely, we use bisection to reduce the length of the base until it reaches 1. However, the factor $a_k$ in the analogue of \eqref{eq:FA2f:recursion:step} is only bounded by $\rho$, because we use the original two-block dynamics of Lemma~\ref{lem:two-block} rather than the non-oriented three-block variant, Lemma~\ref{lem:three-block}. This reflects the fact that our heuristics is based on East-like motion rather than CBSEP-like.

Thus, the only remaining ingredient is dealing with adding a single column to the left of $D\cup\mathring D$. To do this, we observe that, viewing the empty droplet $D$ as a boundary condition, the problem reduces to dealing with FA-1f on a one-dimensional segment, given that it contains at least one empty site. This was already done in Theorem~\ref{th:general:1d}. This completes the sketch of the proof of Theorem~\ref{th:2d:univ:refined:KCM} for the model in Figure~\ref{fig:log3}. Let us note that, the fact that Theorem~\ref{th:general:1d} covers inhomogeneous one-dimensional KCM on their ergodic components is crucial for dealing with more general update families at this point.

For the remaining unbalanced refined universality class ($\beta=1,\gamma=2,\delta=0$), the proof of \cite{Hartarsky24univupper} is still quite similar, using a CBSEP instead of East-type dynamics in the above argument (Lemma~\ref{lem:three-block} versus Lemma~\ref{lem:two-block}). However, balanced models, particularly with finite number of stable directions, are much more delicate. The reason is that one needs to treat scales below the critical one as well, taking into account the non-trivial internal structure of critical droplets which has a multi-scale form as in Figure~\ref{fig:FA2f:droplets}. However, the present situation is more complex than in Section~\ref{sec:FA2f:upper}, because some directions are hard, so one needs to use East-like dynamics on some scales and CBSEP-like on others, also carefully choosing directions in which to grow, depending on the scale. When coupled with the not necessarily rectangular geometry required for general models, as well as bounding conditional probabilities of droplets as in Section~\ref{subsubsec:matryoshka}, but in the absence of symmetry, the proof of \cite{Hartarsky24univupper} becomes quite involved. We direct the reader to \cite{Hartarsky24univupper}*{Section~2} for a detailed description of all mechanisms underlying the proof.

\section{Conclusion}
No new tools were encountered in the present chapter, as compared to previous ones. Instead, we saw how to combine and generalise several of the techniques we were already familiar with, in order to obtain extremely general and precise results. This should clearly showcase the robustness of these methods. For lower bounds we still relied on combinatorial bottlenecks generalising the one for the East model and also incorporating BP ideas. Rough upper bounds used long range renormalisation, while refined ones were proved via the matryoshka dolls technique. One of the takeaways from universality is that a thorough understanding of lower-dimensional models (in our case, one-dimensional) together with respecting the natural geometry and directional preferences of the model can allow one to understand higher dimensional models. The universality viewpoint not only gives a unified framework for understanding the landscape of KCM theory, but also historically supplied the motivation and playground for developing many of the tools presented in the previous chapters.

	\chapter{Out of equilibrium}

\label{chap:out}
  
\abstract{In this chapter, we study KCM with initial state not distributed according to the stationary measure. We start by presenting detailed results on the East model, illustrating the kind of results one would like: exponential convergence to equilibrium after a temperature quench, mixing time cutoff in finite volume, etc. We then treat KCM in full generality at the cost of weakening the results. We conclude with open problems, such as bringing the two aspects above together, and some additional out-of-equilibrium settings. This chapter is independent from Chapters \ref{chap:1d}-\ref{chap:universality}.\\}

In the last three chapters, we discussed properties of the stationary KCM, that is, starting with initial distribution given by the invariant non-trivial product measure $\mu=\mu_q$. While this is the first setting to explore in order to understand these models, it is not the only one of interest. Indeed, both from the mathematics and the physics perspective, it is relevant to study KCM with other initial conditions, that is, out of equilibrium. 
An initial basic question is to determine under which conditions on $\cU$, $q$ and the initial configuration $\omega
$, the law of the corresponding infinite volume KCM $(\eta(t))_{t\ge0}$ converges to the equilibrium measure $\mu$, as time goes to  infinity. 
It is natural to expect that, in the ergodic regime, $q>\qc$, provided $\omega
$ has ``enough empty sites'' (at the very least, one should have $[\omega
]=\bzero$ in view of Section~\ref{sec:legal:paths}) it should hold that for any local function $f$,
\begin{equation}
\lim_{t\to\infty}|\bbE_{\omega
}(f(\eta(t)))-\mu(f)|=0\label{hope},\end{equation}
Similarly, if  the process is initialised according to a distribution $\nu$, one could expect that if $\nu$ has a sufficiently high density of empty sites and $q>\qc$, we should have
\begin{equation}\lim_{t\to\infty}|\bbE_{\nu}(f(\eta(t)))-\mu(f)|=0\label{hope_nu}.\end{equation}
A natural choice is $\nu=\mu_{q_0}$ with $q_0\neq q$. In the physics jargon this is known as a \emph{temperature quench}:
 abruptly changing the temperature from one value to another (recall \eqref{eq:temperature}).\footnote{
 The extreme
 case, $q=1$, is well understood. Indeed, the process becomes a continuous-time version of BP, and essentially behaves like BP.
 }

Unfortunately, robust tools to prove \eqref{hope} and \eqref{hope_nu} are  not yet available. Indeed, with the notable exception of the East model (see Section~\ref{sec:East-out}), results are far from being satisfactory and are limited to a restrictive regime of high $q$ (see Section \ref{sec:high-temp}).
The first and foremost reason for this is the fact that, even though the constraints \eqref{eq:def:cx} are monotone, KCM are not attractive (that is, the product partial order is not preserved by the semi-group of the process, see \cite{Liggett05}*{Sections II.2 and III.2} for background). This is due to the fact that the presence of more empty sites may make certain constraints satisfied and therefore allow certain empty sites to become occupied. Consequently, many of the powerful techniques  (e.g.\ censoring or coupling arguments) which have been developed for the study of other Glauber dynamics (e.g.\ the contact process and stochastic Ising model, see \cites{Liggett05,Liggett99,Martinelli99}), fail for KCM.
To make matters worse, the usual Holley--Stroock strategy \cite{Holley89} to prove convergence to equilibrium does not apply. Indeed, this approach uses the finiteness of the logarithmic Sobolev constant, which implies hypercontractivity of the semigroup. However, due to the presence of constraints, the logarithmic Sobolev constant is infinite for KCM (see Corollary \ref{cor:LS})
and the technique fails.

Another natural
question is how the mixing time (recall \eqref{eq:def:tmix}) scales with the volume
$\Lambda^{(n)}=\{0,\dots,n\}^d$
, as $n\to\infty$ with $\cU$ and $q$ fixed.
For $q>\qc$, we expect linear scaling: there exists
$C>0$ such that, for any 
$\delta\in(0,1)$ and $n$ large enough depending on $\delta$,
\begin{equation}\label{hope_tmix}
t_{\rm{mix}}^{(n)}(\delta):=t_{\rm{mix}}^{\Lambda^{(n)}}(\delta)\leq Cn.
\end{equation}
While lower bounds linear in $n$ are easy to obtain (see Proposition \ref{prop:tmix:easy}), proving \eqref{hope_tmix} is more challenging and will be addressed in the next sections. Once linear upper and lower bounds are established, it is natural to seek a finer \emph{cutoff} result.
Namely, we expect that there exists $v$ independent of $\delta$ such that
\begin{align}
\label{hope_cutoff}
t_{\rm{mix}}^{(n)}(\delta)&{}= v n +\epsilon_{n,\delta},&\lim_{n\to\infty}\epsilon_{n,\delta}/n&{}=0.\end{align}
Proving this cutoff is essentially equivalent to establishing a limit shape result like those known e.g.\ for the contact process (see \cite{Durrett80}). This question was raised by Kordzakhia and Lalley \cite{Kordzakhia06} for the North-East KCM and is readily supported by simulations (see \cite{Kordzakhia06}*{Figure~1}).

\section{Oriented KCM}
Before turning to the East model, let us gather a few useful facts which hold more generally for oriented models. For $\vec u\in\bbR^d\setminus \{0\}$, we consider the open half-space $\bbH_{\vec u}=\{\vec x\in\bbR^d:\<\vec x,\vec u\><0\}$, as in \eqref{eq:openhalf}.

\begin{definition}[Oriented KCM]
\label{def:oriented}
Fix a dimension $d\ge 1$ and a $d$-dimensional update family $\cU$. We say that the $\cU$ is oriented (and that the corresponding KCM is \emph{oriented}) if 
 there exists $\vec v\in \bbR^d\setminus\{0\}$ such that $\bbH_{\vec v}\supset\bigcup_{U\in\cU}U$.
\end{definition}
\begin{remark}[{Oriented examples}]
Among the models introduced in Section \ref{sec:rules} only the East model and the North-East model are oriented (in every dimension). In both cases
a possible choice is $\vec v=-\sum_{i=1}^d \vec e_i$.
\end{remark}
The propositions below state two very handy properties shared by all oriented KCMs: 
\begin{itemize}
    \item dependence propagates only in one direction among well-chosen hyperplanes;
    \item conditionally on a given site having been legally updated, its occupation variable has its equilibrium distribution.
\end{itemize}

\begin{proposition}[{Oriented dependence}]
\label{prop:handy1}
Fix an oriented update family $\mathcal U$, a site $\vec x\in\mathbb Z^d$ and let $\bbH_{\vec v}\supset\bigcup_{U\in\cU}U$.
Recalling the graphical construction (see Section \ref{subsec:markov}), the restriction $(\eta_{\vec x+\bbH_v}(t))_{t\ge 0}$ of the KCM process is independent of the initial condition, clock rings and coin tosses in $\bbZ^d\setminus(\vec x+\bbH_{\vec v})$.
\end{proposition}
In the case of the one-dimensional East model,
Proposition \ref{prop:handy1} states that a site $x\in\bbZ$ is only influenced by the restriction of the process to its right.
\begin{corollary}[{Exact equilibrium}]
\label{cor:handy2} Fix an oriented update family $\mathcal U$, an initial configuration $\omega$, and a site $\vec x\in\mathbb Z^d$.  Let $\mathcal E_{\vec x}(t)$ for $t\geq 0$ be the event that there has already been a legal update (see Section \ref{subsec:markov}) at $\vec x$ by time $t$. Then
\[\mathbb P_{\omega}(\eta_t(x)=1|\mathcal E_x(t))=q.\]
\end{corollary}
The proofs of Proposition~\ref{prop:handy1} and Corollary~\ref{cor:handy2} are left as an exercise to the reader. 

The following result was proved by Chleboun and Martinelli \cite{Chleboun13}, using an idea similar to Proposition~\ref{prop:handy1}.
\begin{theorem}[Quasi-linear mixing for oriented KCM]
\label{th:CM}
Fix a dimension $d\ge 1$ and a $d$-dimensional oriented update family $\cU$.
For any $q>\qc$ there exists $C>0$ such that 
\[t^{(n)}_{\mathrm{mix}}(\delta)\le C n\log n\]
for any $\delta\in(0,1)$ and $n$ large enough depending on $\delta$.
\end{theorem}
Though Theorem~\ref{th:CM} is expected to be sub-optimal (see \eqref{hope_tmix} and Conjecture \ref{conj:tmix} below), we recall the main ideas of its simple and instructive proof. 
\begin{proof}[Sketch]
Let $\Lambda=\{1,\dots,n\}^d$. We proceed iteratively by decomposing $\Lambda$ into its sections 
\[\left\{\bx\in\Lambda:\<\bx,\bv\>=\lambda\right\}\] 
for $\lambda\in\bbR$ and $\bv\in\bbR^d\setminus\{0\}$ as in the definition of an oriented update family. Choosing the direction $\bv\in\bbZ^d$, we obtain order $n$ such sections of cardinality at most order $n^{d-1}$. 

Observe that sites on the first hyperplane (corresponding to the smallest value of $\lambda$ above) are unconstrained thanks to the boundary condition, so a classical coupon collector argument (see e.g.\ \cite{Levin09}*{Section~5.3.3}) shows that its mixing time is of order at most $\log n$. The idea is to show that after a time of order $i\log n$, the distribution of the state of the first $i$ hyperplanes is very close to the stationary one. Then the inductive step is performed using the fact that $\trel<\infty$ whenever $q>\qct$, by Theorem~\ref{th:exponential} and $\qc=\qct$ for oriented update families by \cite{Hartarsky22sharpness}*{Corollary 1.8} (recall Conjecture~\ref{conj:qc:qct}, which is established in this case). Indeed, using $\trel<\infty$ along the lines of the proof of Theorem~\ref{th:exponential}, we show that it is likely that, at time of order $\log n$, each site of the last hyperplane has received a legal update, assuming that the initial marginal on the first $i-1$ sections is exactly (or close to) the equilibrium measure $\mu_q$. Of course, the fact that the equilibrium marginal is preserved over time is true thanks to the oriented nature of the constraint.
\end{proof}

\section{East model}
\label{sec:East-out}
\subsection{Results}
In addition to orientation, the East model has other helpful features enabling the proof of rather detailed results outlined next. Let $\bbN=\{0,1,\dots\}$ and $\Delta_{\vec x}^d=\vec x-(\bbN^d\setminus\{0\})$.

\begin{theorem}
[Exponential convergence in all dimensions]\label{theo:East_out} Consider the East model on $\mathbb{Z}^d$
with $q\in(0,1]$. Fix $\vec x\in\mathbb Z^d$ and an initial configuration $\omega \in \Omega$ such that $\omega_{\vec x+\bbN^d}\neq\bone_{\vec x+\bbN^d}$. Then there exists $m=m(\omega,q)>0$ and for each local function $f$ with support contained in $\Delta_{\vec x}^d$ there exists $C=C(f,\omega,q)$ 
 such that for $t>0$ sufficiently large it holds that
\[
  \left| \mathbb{E}_{\omega}(f(\eta(t))) - \mu(f) \right|\leq C e^{-m t}.
  \]
\end{theorem}
\begin{remark}\label{rem:East}Theorem \ref{theo:East_out}, which was proved in one dimension by Cancrini, Schonmann, Martinelli and Toninelli \cite{Cancrini10} and extended by Mar\^ech\'e  \cite{Mareche19a} to any dimension, implies  relaxation to equilibrium in the sense of \eqref{hope} with exponential decay
for the minimal possible initial condition, namely as soon as $\omega$ has at least one  empty site in $\vec x+\bbN^d
$ for any $\vec x\in\mathbb Z^d$.
Using Theorem \ref{theo:East_out}, one can also easily prove convergence to equilibrium in the sense of \eqref{hope_nu} with  $\nu= \mu_{q_0}$ for any $q_0\in(0,1]$ (see \cite{Cancrini10}*{Theorem 4.3} and  \cite{Mareche19a}*{Theorem 2.2 and Remark 2.1}). Furthermore, this convergence occurs exponentially fast, namely 
there exists $m=m(q_0,q)>0$ and, for $f$ local, there exists $C=C(f,q_0,q)$ such that
\[
  \left| \mathbb{E}_{\mu_{q_0}}(f(\eta(t))) - \mu(f) \right|\leq C e^{-m t}.
  \]
In dimension one, the dominant term of the time scale of convergence to equilibrium when $q\downarrow 0$ has also been established. It holds that $\trel \leq m^{-1} \leq \trel \log (1/q)=e^{(\log (1/q))^2/(2\log 2)} \log (1/q)$. The lower bound follows
by using an argument similar to the one in Section \ref{subsubsec:bottleneck:deduction}, the upper bound may be found in \cite{Faggionato13}*{Theorem 3.5}. We also refer to \cite{Blondel13b}*{Proposition~4.3} for non-local functions.
\end{remark}

\begin{theorem}[Linear mixing time in all dimensions]\label{th:East:precut} For the East model in any dimension   \eqref{hope_tmix} holds. Furthermore, the result also holds when instead of the completely empty boundary condition we fix any  ergodic boundary condition (see Definition \ref{def:ergoBC}).
\end{theorem}

We will provide a proof of Theorem \ref{theo:East_out} in Section \ref{sec:East-out1}
for $d=1$ and refer the reader to \cite{Mareche19a} for $d>1$. Concerning the proof of Theorem \ref{th:East:precut}, which is based on similar  ingredients,  we refer the reader to \cite{Chleboun15}. 

Let us now turn to finer results which have been proved only in the one dimensional case.
Consider the  East model on $\mathbb Z$ with initial condition $\omega\in\Omega$ such that $\omega_{-\bbN}=\bone_{-\bbN\setminus\{0\}}\cdot\bzero_{\{0\}}$. Let $X_t=\min\{x\in\bbZ:\eta_x(t)=0\}$ be the position of the leftmost empty site of $\eta(t)$, which we call the \emph{front}. Notice that $X_t$ makes only nearest neighbour jumps. Indeed, an empty site can appear or disappear only when its left neighbour is empty.
Since, once a site is legally updated its occupation variable is forever set to equilibrium (see Corollary~\ref{cor:handy2}), one
could imagine that to the right of the front
the distribution is $\mu$. If this were true, the front would move  as a biased random walk: negative increments would  occur at rate $q$ (the constraint is always satisfied on the site $X_{t}-1$), and positive increments  at rate $q(1-q)$ (the occupation variable  at the position of the front can be updated to occupied only if $X_t+1$  is empty).
This would yield a speed $q^2$ to the left. 
While it is not true that the configuration seen from the front has distribution $\mu$, the following result due to Ganguly, Lubetzky and Martinelli \cite{Ganguly15}, confirms that
(as for a biased random walk) the front moves at a negative speed with normal fluctuations and that its concentrated passage times imply cutoff (see \eqref{hope_cutoff}) with a window $\sqrt n$.  
\begin{theorem}[East in one dimension: CLT for the front] \label{theo:CLT} Let $d=1$ and $\omega\in\Omega$ be such that $\omega_{-\bbN}=\bone_{-\bbN\setminus\{0\}}\cdot\bzero_{\{0\}}$. For any $q\in(0,1)$, there exist constants $\sigma=\sigma(q)>0$, $v=v(q)<0$ and $C=C(q)>0$
such that, the front $X_t$ of the one-dimensional East process with initial condition $\omega$ satisfies
\begin{align*}
\bbP\left(\lim_{t\to\infty}\frac{X_t}{t}=v\right)&{}=1,&
\left|\mathbb E(X_t)-vt\right|&{}\le C,\\
\lim_{t\to\infty}\frac{\Var(X_t)}{t}&{}=\sigma^2,&
\lim_{t\to\infty}\frac{X_t-vt}{\sigma \sqrt t}
&{}=Z
\end{align*}
in distribution for a standard normal random variable $Z$.
\end{theorem}

\begin{corollary}[Cutoff in one dimension]
\label{th:East:cutoff}
For any $q\in(0,1)$, there exists $C=C(q)>0$ such that the one-dimensional East model satisfies \eqref{hope_cutoff}  with \[|\epsilon_{n,\delta}|\le C\phi^{-1}(1-\delta)\sqrt n,\]
where $\phi$ is the cumulative distribution function of the standard 
normal law.
\end{corollary}

We refer the reader to \cite{Ganguly15} for the proof of Theorem \ref{theo:CLT} based on: 
\begin{itemize}
    \item ergodicity of the process seen from the front (see Theorem \ref{theo:ergoEast} below); 
    \item showing that, after an initial burn-in time, the front increments behave like a stationary sequence of weakly dependent random variables, and applying an ingenious Stein's method argument of Bolthausen \cite{Bolthausen82} to derive the central limit theorem.
\end{itemize}
Corollary~\ref{th:East:cutoff} follows almost immediately. We further refer the reader to \cite{Couzinie22} for a first step in the direction of proving \eqref{hope_cutoff} for East in higher dimensions.

\begin{theorem}[Ergodicity  of the process seen from the front]\label{theo:ergoEast}
    Let $d=1$ and $\omega\in\Omega$ be such that $\omega_{-\bbN}=\bone_{-\bbN\setminus\{0\}}\cdot \bzero_{\{0\}}$. Let $(\eta(t))_{t\ge 0}$ be the one-dimensional East KCM with initial condition $\omega$. For all $t\ge 0$ and $x\in\bbZ$, let $\tilde\eta_x(t)=\eta_{x+X_t}(t)$, defining the process seen from the front. 
There exists a unique measure $\tilde\mu=\tilde\mu(q)$ such that
$\tilde\eta(t)\to\tilde\mu$ in distribution as $t\to\infty$. Moreover, there exist $C=C(q)$ and $m=m(q)>0$, such that for any $x\in\bbZ$
\[\|\tilde\mu-\mu\|_{[x,\infty)}\le C e^{-mx},\]
where, for $\Lambda\subset\mathbb Z$,  $\|\tilde\mu-\mu\|_{\Lambda}$ denotes the total variation distance between the marginals of $\tilde\mu$ and $\mu$ on $\Lambda$.
\end{theorem}
\begin{remark}[{Front velocity}]\label{rem:velo}
The velocity of the front $v$ in Theorem \ref{theo:CLT} can be expressed in terms of the invariant measure as
$v=-q+(1-q)\tilde\mu(\omega_1=0)$.
\end{remark}
The main ingredients for the proof of Theorem~\ref{theo:ergoEast} due to Blondel \cite{Blondel13b}, in addition to the techniques needed for Theorem~\ref{theo:East_out}, are: 
\begin{itemize}
    \item coupling the processes seen from the front, starting with two different initial conditions, in order to prove that their laws converge to the same limit (see \cite{Blondel13b}*{Theorem 4.7});
    \item using the fact that the front moves at most linearly (by the finite speed argument of Proposition~\ref{prop:tmix:easy}) and leaves empty sites behind to prove that far from the front the process is almost at equilibrium (see \cite{Blondel13b}*{Theorem 4.7}). Here the distinguished zero of Definition~\ref{def:distinguished} below plays a key role.
\end{itemize}

\begin{problem}[{Front measure}]
    Theorems \ref{theo:CLT}-\ref{theo:ergoEast} leave various questions unanswered.  For instance, can one quantify  the correlations between adjacent occupation variables in the invariant measure $\tilde\mu
    $? Can one determine the asymptotics of the velocity $v$ of Theorem~\ref{theo:CLT} and Remark~\ref{rem:velo} in the $q\to 0$ limit? 
\end{problem}

\subsection{Exponential decay to equilibrium for East in one dimension}
\label{sec:East-out1}
Let us start by introducing a key notion, due to  Aldous and Diaconis \cite{Aldous02}.

\begin{definition}[Distinguished zero] \label{def:distinguished}
Fix $x\in\mathbb Z$ and an initial configuration $\omega \in \Omega$ with
$\omega_x=0$. Call the site $x$ \emph{distinguished}. We set
$\xi_0=x$. The position $\xi_s\in \mathbb{Z}$ of the \emph{distinguished zero} at time $s >0$ is defined according to the following iterative
rule. For all times $s$ strictly smaller than the first
legal update $t_1$ at $x$, we set $\xi_s=x$, while $\xi_{t_1}=x+1$. Then we wait for the first legal update $t_{2}$ at $x+1$ after $t_1$, at which point we set $\xi_{t_{2}}=x+2$ and so on.
\end{definition}

Note that, almost surely, the trajectory $(\xi_s)_{s\ge 0}$ is right-continuous, piece-wise
constant, increasing by 1 at each jump and not exploding in finite time. 
Also note that, by definition of the legal updates, necessarily the state $\eta_{\xi_s}(s)$ of the East process at the position of the distinguished zero is 0 for any $s\ge 0$, hence the name.
\begin{proposition}[{Properties of the distinguished zero}]
\label{prop:handyEast}
Fix $t>0$. In the setting of Definition \ref{def:distinguished}, conditioning on the knowledge of the trajectory $(\xi_s)_{s \leq t}$ and denoting by $0 < t_1 < \dots< t_{n-1}<t$ the times of the distinguished zero's jumps (with
$t_0=0$, $t_n=t$), the following holds for all $i\in\{0,\dots,n-1\}$
\begin{itemize}\item   for $s\in[t_i,t_{i+1})$,  $\xi_s=x+i$ and the restriction of the process to $\{x,\dots,x+i-1\}$ 
follows an East dynamics with zero boundary condition;
\item at time $t_{i+1}$, the configuration at site $x+i$ is updated according to a Bernoulli$(1-q)$ variable.
\end{itemize}
\end{proposition}
The proof, which is derived using  Definition \ref{def:distinguished}, Proposition \ref{prop:handy1} and Corollary \ref{cor:handy2}, is left to the reader. 
\begin{remark}[{Conditional graphical construction}]
Recall the graphical construction of Section \ref{subsec:markov}. Proposition \ref{prop:handyEast} implies that, given  $(\xi_{t'})_{t'\le t}$,  for any $(s, y)$ satisfying
$s\leq t$ and $y\leq\xi_s$, the variable $\eta_y(s)$ is  determined by the following ``conditional graphical construction''. In the time interval $[0,t_1)$, the occupation variables  in the interval
$\{y,\dots, x-1\}$ evolve as an East process with fixed empty site at $x$, using the clock rings $(t_{z,k})_{z\in\{y,\dots, x-1\}}$ and coin tosses $(s_{z,k})_{z\in\{y,\dots, x-1\}}$. In the time interval $[t_1,t_2)$ the same happens with empty boundary condition at $x+1$ instead of $x$ and so on. Note that, at $t_1$, the clock at $x$ rings and $\eta_x(t_1)$ takes the value of the corresponding coin toss.
\label{rem:conditional}
\end{remark}
Furthermore, Proposition \ref{prop:handyEast} yields the following.
\begin{corollary}[{Equilibrium zone}]
\label{cor:handyEast}
Fix $a\leq x\in\bbZ$ and two measures $\psi^-$ and $\psi^+$ on $\Omega_{a-1-\bbN}$ and $\Omega_{x+1+\bbN}$ respectively. Let
\[
 \psi=
    \begin{cases} \psi^-\otimes \mu_{\{a,\dots,x-1\}}\otimes\delta_{0}\otimes\psi^+  & \text{ if $a<x$}\\
     \psi^-\otimes\delta_{0}\otimes\psi^+  & \text{if $a=x$}
     \end{cases}
\]
Let $(\eta(t))_{t\ge0}$ be the East process with initial distribution $\psi$. Fix $t\ge 0$. Then the conditional distribution of $\eta_{\{a,\dots,\xi_t-1\}}(t)$, given the distinguished zero trajectory $(\xi_s)_{s\le t}$, is the equilibrium one, $\mu_{\{a,\dots,\xi_t-1\}}$.
\end{corollary}
The following proof, due to Cancrini, Martinelli, Schonmann and Toninelli \cite{Cancrini10},
uses as key ingredients  Proposition \ref{prop:handyEast}, Corollary \ref{cor:handyEast} and  the $\mathrm{L}^2$ convergence to equilibrium guaranteed by Theorem~\ref{th:exponential}\ref{item:expo:7} and \eqref{eq:expo}. 
\begin{proof}[Theorem \ref{theo:East_out} for $d=1$]
Assume for simplicity that $\mu(f)=0$ , $x=0$, $\omega_0=0$ and let the support of $f$ be contained in $\{a,\dots,a'\}\Subset-1-\bbN$.
Let $b\leq 0$ be the position of the first empty site in $\omega$ to the right of $a'$.  Make $b$ distinguished and denote by $\xi_s$ its position at time $s$.
Given the trajectory $(\xi_s)_{s \leq t}$, let $0 < t_1 < t_2 < \dots< t_{n-1}<t$ be the times when the distinguished zero jumps,
and set $t_0=0$, $t_n=t$. For $i\in\{0,\dots,n-1\}$, let \begin{align*}
{\Xi}
_i&{}=(\xi_s)_{s\in [t_i,t]}&V_i&{}=\{a,\dots,b+i-1\}.\end{align*}
We claim that
\begin{align} 
\nonumber\left|\bbE_\omega\left(f(\eta(t))\right)\right| & {}\leq
\bbE_\omega \left( \left|\bbE_\omega(f(\eta(t)) | \Xi_0) \right| \right)  \\
 \nonumber &{}\leq(\min(q,1-q))^{-(b-a)}\bbE_\omega \left( \int \mathrm{d}\mu_{V_0}(\omega') \left|\bbE_{{\omega'}}
(f(\eta
'(t)) | \Xi
_0) \right| \right)\\
\nonumber&{}
\leq
(\min(q,1-q))^{-(b-a)} \bbE_\omega
\left(\sqrt{\int \mathrm{d}\mu_{V_0}(\omega')(\bbE_{\omega'}(f(\eta'(t)) |\Xi
_0))^2}\right)
\\
\label{eq:Eastbound}&{}
\leq (\min(q,1-q))^{-(b-a)} \bbE_\omega \left( \sqrt{\var_{V_0} \left(g_t^{(0)}\right)}\right),
\end{align}
where we let $\eta
'(t)$  be the configuration obtained following the conditional graphical construction of Remark \ref{rem:conditional} with $\eta'_{V_0}(0)=\omega'\in\Omega_{V_0}$
and
\[g_t^{(0)}(\omega')=\bbE_{\omega'}(f(\eta
'(t)) | \Xi
_0
).\]
In order 
to obtain  \eqref{eq:Eastbound}, we used 
Remark \ref{rem:conditional} together with the fact that for any $\omega\in\{0,1\}^{\mathbb Z}$ it holds $\mu_{V_0}(\omega|_{V_0})\geq (\min(q,1-q))^{-(b-a)}$ to obtain the second inequality,  Cauchy-Schwarz for the third inequality and, for the last inequality, we used
Corollary \ref{cor:handyEast} and the hypothesis that the support of $f$ is contained in  $[a,b)$, which yield
\begin{equation}
\int \mathrm{d}\mu_{V_0
}(\omega') \bbE_{\omega'}(f(\eta
'(t)) | (\xi_s)_{s \leq t})=
\mu_{V_0
} (f)
=0.\label{berluingalera}
\end{equation}

Let $P_s^{(i)}$ for $ s\in [t_i,t_{i+1})$ be the
Markov semigroup associated to the East process in the interval
$V_i$ with a fixed empty site at the right boundary $b+i$. Then, using Propositions 
\ref{prop:handy1} and \ref{prop:handyEast}, we get 
\begin{align}
\label{fru}
g_t^{(0)}({\omega'})& = \sum_{\sigma \in \{0,1\}^{V_0}}\sum_{\sigma' \in \{0,1\}} 
 P_{t_1}^{(0)}(\omega',\sigma)\mu_b(\sigma')g^{(1)}_{t-t_1}(\sigma\cdot\sigma'),
\end{align}
where, for $s\geq 0$, we define $g_s^{(1)}:\Omega_{V_1}\to\mathbb R$ by 
\begin{equation}\label{g0}g_s^{(1)}(\omega'')=\bbE_{\omega''}(f(\eta
''(s)) | \Xi_1
)
\end{equation}
where $\eta''(s)$    for $s\leq t$ denotes the configuration
in the interval $\{a,\dots,\xi_t-1\}$ obtained starting at time $t_1$ from the configuration $\omega''\in\Omega_{V_1}$ and evolving according to the conditional graphical construction described
in Remark \ref{rem:conditional} applied to the time interval $(t_1, t_1+s]$.
Therefore, using \eqref{eq:expo} we get
\begin{align}\nonumber
 \var_{V_0} \left(g_t^{(0)}\right)  &{}\leq e^{-2 t_1 /T_{\rm{rel}}^{V_0}} \var_{V_0}
\left(
\sum_{\sigma' \in \{0,1\}} \mu_b(\sigma') \mu_{{V_1}}\left(g_{t-t_1}^{(1)}\right)\right)\\
&{}\leq
e^{-2 t_1/T_{\rm{rel}}} \var_{V_1} 
 \left(g_{t-t_1}^{(1)}\right) \label{fru2}
\end{align}
where, as usual, we denote by $T_{\rm{rel}}^{V_0}$ (resp.\ $T_{\rm{rel}}$) the relaxation time of the East model on $V_0$ with empty boundary condition (resp.\ $\mathbb Z$).  In order to obtain  \eqref{fru2} we use convexity of the variance and \eqref{eq:finite:volume:chain}. 
We can now proceed analogously 
to get
\begin{equation}
\var_{V_1} \left(g^{(1)}_{t-t_1}\right)  \leq e^{-2(t_2-t_1)/{T_{\rm{rel}}}}
\var_{V_2}\left(g^{(2)}_{t-t_1-t_2}\right)
\end{equation}
and, by  induction, we get
\[\var_{V_0} \left(g_t^{(0)}\right) \leq e^{-2  t/T_{\rm{rel}}} \var_{\{a,\dots,\xi_t-1\}}(f).\]
Plugging this bound into \eqref{eq:Eastbound} and recalling that $T_{\rm{rel}}>0$ for any $q\in{(0,1]}$ for the East model, we finally get
\[|\bbE_\omega(f(\eta(t)))| \leq c {e}^{-t/\trel}  \bbE_\omega \left( \sqrt{\var_{\{a,\dots,\xi_t-1\}} (f)}
\right) \leq ce^{-t/\trel}  
\|f\|_\infty \]
for some $c=c(q,a,b)>0$.
\end{proof}

While the distinguished zero does not generalise to oriented KCM, a version of it was used in 
\cite{Cancrini10}*{Section 4} for an oriented KCM on a tree.

\section{High vacancy density regime}
\label{sec:high-temp}
In this section we focus on results on KCM with arbitrary update family with $q$ close enough to 1. While the same results should hold whenever $q>\qc$, the current techniques do not allow proving this. 
\subsection{Results}
The next two theorems, proved by Hartarsky and F.~Toninelli \cite{Hartarsky24KCM-CP}, establish \eqref{hope_nu} with an exponential convergence when $\nu$
is a product measure with a  vacancy density $q_0>q_c$ and a linear upper bound for the mixing time \eqref{hope_tmix}.
Both results apply to all models but are restricted to a high vacancy density regime for $q$ and
their proof is delegated to Section~\ref{subsec:out-of-eq:high:proof}.
Before stating them  we need to introduce the notion of trivial subcritical models.
\begin{definition}[{Trivial subcritical}]
  We say that an update family $\cU$ is not \emph{trivial subcritical}, if there exists $U\in \cU$ and a direction $\vec v\in \bbR^d$ such that $\<\vec u,\vec v\><0$ for all $\vec u\in U$. 
\end{definition}
In \cite{Balister24}*{Theorem~7.1} it is proved that $q_c=1$ iff $\cU$ is trivial subcritical. Therefore, excluding trivial subcritical models from the following two results is necessary.

\begin{theorem}[Exponential convergence to equilibrium after a temperature quench]
\label{th:convergence:out-of-eq}
Fix a dimension $d\ge 1$ and an  update family $\cU$ which is not trivial subcritical. For any $\alpha>0$ there exist $\varepsilon=\varepsilon(\alpha)>0$ and $c=c(\alpha)>0$ such that the following holds for the $\cU$-KCM $\eta$. For any $q_0\in[\qct+\alpha,1]$, $q\in[1-\varepsilon,1]$ and local function $f:\Omega\to \bbR$,
\[\left|\bbE_{\mu_{q_0}}(f(\eta(t)))-\mu_q(f)\right|\le \frac{\|f\|_\infty\cdot|\mathrm{supp}f|}{ce^{ct}},\]
where $\mathrm{supp}f$ is the set of sites on whose state the value of $f$ depends.
\end{theorem}
 
\begin{theorem}[Linear mixing at high vacancy density]
\label{th:mixing}
Fix a dimension $d\ge 1$ and a $d$-dimensional update family $\cU$ which is not trivial subcritical. Then there exist $\varepsilon>0$ and $C>0$ such that for any $q\in[1-\varepsilon,1]$, inequality \eqref{hope_tmix} holds 
 for any $\delta\in(0,1)$ and $n$ large enough depending on $\delta$.
\end{theorem}
It should be noted that Theorem~\ref{th:mixing} also applies to domains of non-hypercubic shape, but is stated as is for the sake of simplicity. We refer the reader to \cite{Hartarsky24KCM-CP}*{Section~4} for an account of previous works in this direction, in particular \cites{Blondel13,Mountford19,Mareche20}. Moving on to more precise results, it only remains to report the following one due to Ertul \cite{Ertul22} (also recall Corollary~\ref{th:East:cutoff} for the one dimensional East model).
\begin{theorem}[Cutoff for one-dimensional FA-1f]
\label{th:FA1f:cutoff}
Consider the  FA-1f model in one dimension.
There exists an explicit $\varepsilon>0$ such that for any $q\in(1-\varepsilon,1]$, there exists $v,\alpha,\beta>0$ such that \eqref{hope_cutoff} holds with 
$$-\alpha\sqrt n\leq \epsilon_{n,\delta}\leq \beta\sqrt n.$$
\end{theorem}
Furthermore, the number $v$ in Theorem~\ref{th:FA1f:cutoff}  corresponds to twice the speed at which the rightmost empty site, known as \emph{front}, moves in FA-1f on $\{1,2,\dots\}$ with boundary condition $\bzero_{\{0,-1,\dots\}}$ (compare with Remark \ref{rem:velo} for East). Indeed, as for Corollary~\ref{th:East:cutoff}, the proof of Theorem~\ref{th:FA1f:cutoff} is based on the identification of the front speed thanks to the convergence of the process seen from the front \cite{Blondel19} (as in Theorem~\ref{theo:ergoEast}). Then, considering each interval of occupied sites, one may show that it shrinks at the front speed at both ends, leading to a double speed (see \cite{Ertul22} for more details). Nonetheless, deducing Theorem~\ref{th:FA1f:cutoff} from the analogue of Theorem~\ref{theo:ergoEast} is not immediate, as in the case of Corollary~\ref{th:East:cutoff}.

\subsection{Proofs via cooperative contact processes}
\label{subsec:out-of-eq:high:proof}
We now turn to outlining the proofs of Theorems~\ref{th:convergence:out-of-eq} and~\ref{th:mixing}. In order to simplify the presentation (see \cite{Hartarsky24KCM-CP} for the full details\footnote{Beware that in \cite{Hartarsky24KCM-CP} the roles of 0 and 1 states are exchanged.}), we focus on Theorem~\ref{th:mixing} in the case of FA-2f in two dimensions (see Figure~\ref{fig:FA2f}). The only additional ingredient needed in the general case and in Theorem~\ref{th:convergence:out-of-eq} are one-scale renormalisations along the lines of Section~\ref{sec:exponential:decay}.

The proof proceeds in several steps involving a number of interacting particle systems other than KCM. In Section~\ref{subsubsec:CP}, we reduce the mixing time problem for FA-2f to a more complex problem for a simpler model. Namely, studying the space-time connected components of occupied sites in a cooperative contact process known as the sexual contact process. In Section~\ref{subsubsec:LPP}, we use a comparison with last passage percolation to replace the arbitrary initial condition by the fully empty one. In Section~\ref{subsubsec:BPwD} we discretise time to transform the sexual contact process into a North-East BP with death. In Section~\ref{subsubsec:Toom}, we introduce Toom cycles to show that, at $q$ close enough to 1, with empty initial condition, BP with death mostly has empty sites. Finally, in Section~\ref{subsubsec:chains}, we show that long chains of such Toom cycles are also unlikely.

\subsubsection{Sexual contact process}
\label{subsubsec:CP}
Let us start by introducing the \emph{sexual contact process} (SCP) of \cite{Durrett86}. It is a continuous time Markov process on $\{0,1\}^{\bbZ^2}$ (or in finite volume $\Lambda$ with $\bzero_{\bbZ^2\setminus\Lambda}$ 
boundary condition as in Section~\ref{sec:finite_vol}) which, using a graphical representation similar to the one for KCM 
from Section~\ref{subsec:markov}, is defined as follows. Each site $\bx\in\bbZ^2$ waits an independent exponentially distributed time with mean one before attempting to update. At that time, if both $\bx+\be_1$ and $\bx+\be_2$ are in state 0, the state of $\bx$ becomes $0$ with probability $q$ and becomes $1$ with probability $1-q$. Otherwise the state of $\bx$ remains unchanged with probability $q$ and becomes $1$ with the remaining probability $1-q$. SCP has two notable advantages over FA-2f and one disadvantage. Namely, SCP is attractive and oriented, but its upper invariant measure is not explicit when it is not the trivial $\delta_\bone$.

The above formulation of the process suggests a canonical coupling with \mbox{FA-2f}, using the same clock rings and Bernoulli variables with parameter $q$, the same domain $\Lambda=\{1,\dots,n\}^2$ and boundary condition $\bzero_{\bbZ^2\setminus\Lambda}$. It is an exercise to check the following.
\begin{lemma}[Sexual contact process comparison]
\label{lem:SCP:comparison}
Under the canonical coupling of FA-2f $(\eta(t))_{t\ge 0}$ and SCP $(\zeta(t))_{t\ge 0}$, if $\zeta(0)=\bone_\Lambda$, then $\eta_\bx(t)\le \zeta_\bx(t)$ for all $\bx\in\Lambda$ and $t\ge 0$.
\end{lemma}
However, we require a finer relation from \cite{Hartarsky24KCM-CP}*{Section~7} (also see \cite{Liggett99}*{Section~III.1}, \cite{Busic13}*{Section~4.2} and \cite{Gottschau18}*{Section~1.3} for similar ideas), since the two processes do not share the invariant measure $\mu_q$. Consider a set of \emph{orange} sites $O_t\subset \Lambda$ defined as follows for $t\ge 0$. At time $0$, we set $O_0=\Lambda$. At each clock ring $t\ge 0$ at site $\bx\in\Lambda$, we obtain $O_t$ from $O_{t-}$ by:
\begin{itemize}
    \item removing $\bx$, if $\zeta_\bx(t)=0$;
    \item adding $\bx$, if $\zeta_\bx(t)=1$ and there is an orange site around $\bx$, that is, $O_{t-}\cap \{\bx,\bx+\be_1,\bx+e_2,\bx-\be_1,\bx-\be_2\}\neq \varnothing$;
    \item changing nothing otherwise.
\end{itemize}
The purpose of orange sites is to ensure that FA-2f processes with different initial conditions are coupled outside orange sites.
\begin{lemma}[Orange set coupling]
\label{lem:orange}
Under the canonical coupling of SCP $(\zeta(t))_{t\ge 0}$ with initial condition $\bone_\Lambda$ and two FA-2f processes $(\eta(t))_{t\ge 0}$ and $(\eta'(t))_{t\ge 0}$ with different initial conditions, for all $t\ge 0$, we have 
\begin{equation}
\label{eq:Ot}\left\{\bx\in\Lambda:\eta_\bx(t)\neq\eta'_\bx(t)\right\}\subset O_t.\end{equation}
\end{lemma}
\begin{proof}
We proceed by induction on the number of clock rings in $\Lambda$. Removing $\bx$, if $\zeta_\bx(t)=0$ is justified by Lemma~\ref{lem:SCP:comparison}. Setting $O_t=O_{t-}\cup\{\bx\}$ cannot violate \eqref{eq:Ot}, at time $t$, if the induction hypothesis is verified at $t-$. Finally, assume that there is no orange site around $\bx$. Then by induction hypothesis $\eta(t-)$ and $\eta'(t-)$ are equal around $\bx$ and, therefore, $\eta_\bx(t)=\eta'_\bx(t)$.
\end{proof}

Owing to Lemma~\ref{lem:orange} and the standard results on coupling and mixing times \cite{Levin09}*{Corollary~5.5}, proving Theorem~\ref{th:mixing} is reduced to showing that with high probability the set of orange sites is empty at time $Cn$ for $C$ large enough. We say that two space-time points $(t,\bx)$ and $(t',\bx')$ are connected in $X\subset[0,\infty)\times\Lambda$, if there exists a sequence of segments $[s_i,t_i]\times\{\bx_i\}\subset X$ indexed by $i\in\{0,\dots,N\}$ such that $[s_i,t_i]\cap[s_{i-1},t_{i-1}]\neq \varnothing$ and $\bx_{i-1}$ and $\bx_{i}$ are neighbours in $\Lambda$, such that $(t,\bx)\in[s_0,t_0]\times\{\bx_0\}$ and $(t',\bx')\in[s_N,t_N]\times\{\bx_N\}$. By the definition of orange sites, it is clear that $\bigcup_{t\ge 0}O_t$ is contained in the connected component $\cC$ of $\{(t,\bx)\in[0,\infty)\times\Lambda:\zeta_\bx(t)=1\}$ containing $\{0\}\times\Lambda$. Thus, it suffices to show that, for $q\ge 1-\varepsilon$, 
\begin{equation}
\label{eq:connected:comp}
\bbP\left(\cC\subset[0,Cn]\times\Lambda\right)\ge 1-\delta
\end{equation}
for some constant $C$ independent of $\delta>0$, provided $n$ is large enough. In order to simplify the exposition, we fix $q=1-\varepsilon<1$ in the rest of the section with $\varepsilon$ to be chosen small enough in Section~\ref{subsubsec:BPwD}.

\subsubsection{Last passage percolation}
\label{subsubsec:LPP}
Following \cite{Hartarsky24KCM-CP}*{Section~11}, our next goal is to replace the worst initial condition, $\bone_\Lambda$, of SCP $\zeta$ of Section~\ref{subsubsec:CP} by the best one, $\bzero_\Lambda$. For this we use a comparison with a version of last passage percolation, which will play a similar role to the orange set above. Define the \emph{last passage} set $L_t\subset\Lambda$ for $t\ge 0$ as follows. Set $L_0=\Lambda$ and at each clock ring $t\ge 0$ at site $\bx\in\Lambda$, we obtain $L_t$ from $L_{t-}$ by removing $\bx$, if the following conditions are both satisfied:
\begin{itemize}
    \item at time $t$ the Bernoulli variable with parameter $q=1-\varepsilon$ 
    takes the value 
    0, that is, SCP changes the state of $\bx$ to 1 regardless of the current configuration;
    \item $\{\bx+\be_1,\bx+\be_2\}\cap L_{t-}=\varnothing$.
\end{itemize}
Otherwise, we set $L_t=L_{t-}$.
The proof of the following observation, similar to Lemma~\ref{lem:orange}, is left as an exercise for the reader.
\begin{lemma}[Last passage percolation coupling]
\label{lem:LPP:coupling}
Under the above coupling, we have
\[\left\{\bx\in\Lambda:\zeta_\bx(t)\neq\zeta'_\bx(t)\right\}\subset L_t\]
for any $t\ge 0$, where $\zeta$ and $\zeta'$ are SCP with initial condition $\bone_\Lambda$ and $\bzero_\Lambda$ respectively.
\end{lemma}
We next invoke a classical fact about last passage percolation. We note that much more precise results are available, but not more useful for our purposes.
\begin{lemma}[Linear last passage time]
\label{lem:LPP}
For some absolute constant $C'>0$ we have
\[\bbP\left(L_{C'n/\varepsilon}\neq\varnothing\right)\le \delta\]
for any $n$ large enough depending on $\delta>0$.
\end{lemma}
\begin{proof}[Sketch]
A proof of Lemma~\ref{lem:LPP} generalising to arbitrary dimension was given by Greenberg, Pascoe and Randall \cite{Greenberg09}. The idea is to introduce an exponential metric on the set of possible values of the last passage set $L_t$:
\[d(A,B)= \sum_{\bx\in A\Delta B}e^{\gamma\<\bx,\bone\>},\]
where $\gamma>0$ is a suitably large constant depending only on dimension and $\Delta$ denotes the symmetric difference of sets. One can verify that the expected distance of two canonically coupled last passage percolations starting from neighbouring configurations contracts. Applying the path coupling method (see \cite{Levin09}*{Chapter 14}), this yields that the hitting time of $\varnothing$ is logarithmic in the diameter of the space, which is exponential in $n$.
\end{proof}

Combining Lemmas~\ref{lem:LPP:coupling} and~\ref{lem:LPP} and recalling \eqref{eq:connected:comp}, it is now sufficient to prove that, with probability $1-\delta$, the connected components of $\{(t,\bx)\in[0,\infty)\times\Lambda:\zeta'_\bx(t)=1\}$ intersecting $\{C'n/\varepsilon\}\times\Lambda$ are contained in $[0,Cn]\times\Lambda$. Taking $C>C'/\varepsilon$ and performing a union bound, it suffices to show that for any space-time point $(t,\bx)$, its connected component $\cC_{t,\bx}$ satisfies
\begin{equation}
\label{eq:diameter}\bbP(\diam(\cC_{t,\bx})\ge k)\le e^{-ck}\end{equation}
for some $c>0$ independent of $n$ and all $k$ large enough. Since SCP is attractive, $\zeta'$ has initial condition $\bzero_\Lambda$ and boundary condition $\bzero_{\bbZ^2\setminus\Lambda}$, it suffices to prove \eqref{eq:diameter} for the infinite volume SCP $\zeta''$ with initial condition $\bzero$.
Indeed, by induction on the number of clock rings in $\Lambda$ in the graphical construction, one can show that for every $t\ge 0$ and $\vec x\in\bbZ^2$, we have $(\bzero_{\bbZ^2\setminus\Lambda}\cdot \zeta'(t))_{\vec x}\le \zeta''_{\vec x}(t)$, which implies that the connected component $\cC_{t,\vec x}$ for $\zeta''$ contains the one for $\zeta'$.

\subsubsection{BP with death}
\label{subsubsec:BPwD}
The next step of the proof is a discretisation in time \cite{Hartarsky24KCM-CP}*{Section~8}. Namely, we fix $T$ large enough but such that $T\varepsilon$ is small. This way, in each time interval of length $T$, it is likely that SCP $\zeta''$ with initial condition $\bzero$ attempts to change the state of a given site $\bx$ to 0 (and succeeds, if $\bx+\be_1$ and $\bx+\be_2$ are in state 0), but never attempts to change the state of $\bx$, $\bx+\be_1$ or $\bx+\be_2$ to 1. We declare each space-time point $(m,\bx)\in\{0,1,\dots\}\times\bbZ^2$ \emph{good}, if the above event occurs for site $\bx$ in the time interval $[mT,(m+1)T)$. Note that almost surely, no clocks ring at times $(mT)_{m=0}^\infty$, so these events are 1-dependent and have probability $1-\varepsilon'$ for $\varepsilon'>0$ that can be chosen arbitrarily small, if $\varepsilon$ is small enough. Using the Liggett-Schonmann-Stacey theorem \cite{Liggett97}, we can replace this 1-dependent field of indicators of good sites by an independent one with high marginals (at the cost of changing $\varepsilon'$).

With the good space-time points at hand, we define the discrete time version $\xi$ of SCP, which we refer to as North-East BP with death parameter $\varepsilon'$. Set $\xi(0)=\bzero$. For $m\ge 1$ and $\bx\in\bbZ^2$, define $\xi_\bx(m)=0$, if $(m-1,\bx)$ is good and $\xi_{\bx}(m-1)=0$, or if $(m-1,\bx)$ is good and $\xi_{\bx+\be_1}(m-1)=\xi_{\bx+\be_2}(m-1)=0$. Otherwise, set $\xi_\bx(m)=1$. The name is justified by the fact that, in the absence of non-good space-time points, this process is exactly BP with the North-East update family (see Figure~\ref{fig:NE}). It is not hard to check that, if $\xi_\bx(m)=\xi_\bx(m+1)=0$, then $\zeta''_\bx(t)=0$ for all $t\in[mT,(m+1)T)$. In view of this, we consider the set 
\begin{equation}
\label{eq:m'x'}
X=\left\{(m,\bx)\in\{0,1,\dots\}\times\bbZ^2:\xi_{\bx}(m)=1\right\}.\end{equation}
equipped with all edges of the form 
$((m,\bx),(m',\bx'))$ with 
\begin{equation}
\label{eq:connected:chain}
(m'-m,\bx'-\bx)\in\{-1,0,1\}\times\{\bzero,\be_1,\be_2,-\be_1,-\be_2\}.\end{equation}
Hence, it suffices to prove \eqref{eq:diameter} for the connected component $\cC'_{m,\bx}$ of any space-time point $(m,\bx)$ in $X$.

\subsubsection{Toom cycles}
\label{subsubsec:Toom}
\begin{figure}
    \centering
    \includegraphics[width=\linewidth]{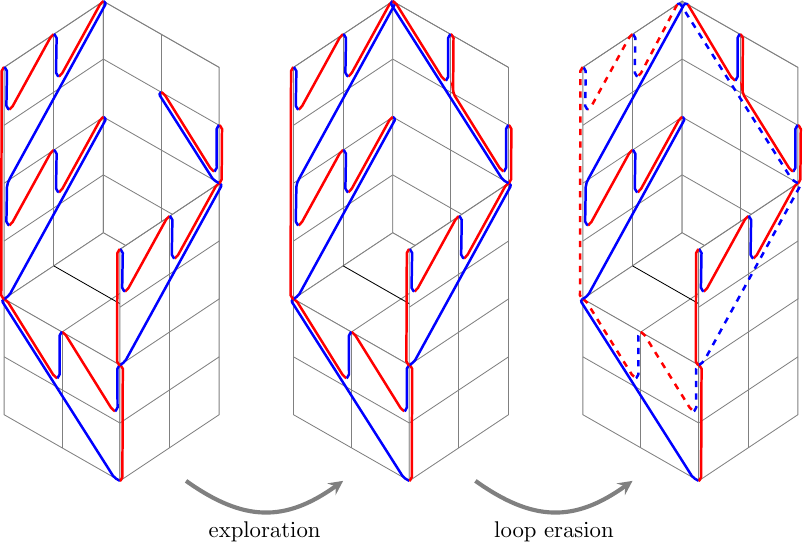}
    \caption{Illustration of the algorithm used to construct a Toom cycle rooted at the bottom space-time point $(m,\vec x)$. Time increases downwards.}
    \label{fig:Toom}
\end{figure}

Before treating connected components $\cC'_{m,\bx}$ in $X$ of \eqref{eq:m'x'}, it is useful to first show that $\bbP(\xi_\bx(m)=1)$ is small for any space-time point $(m,\bx)$. This is an instance of a classical result of Toom \cite{Toom80} for stability of cellular automata subjected to random noise. However, since our cellular automaton is BP with update family consisting of a single rule, it is possible to use a simpler argument of Swart, Szab\'o and Toninelli \cite{Swart22}*{Section~3.5} presented below. We present the construction of a \emph{Toom cycle} in the form of an algorithm illustrated in Figure~\ref{fig:Toom}.

Fix a realisation of the good space-time points such that $\xi_\bx(m)=1$. For each $(m',\bx')$ such that $\xi_{\bx'}(m')=1$ and $(m'-1,\bx')$ is good, fix some $\be(m',\bx')\in\{\be_1,\be_2\}$ such that $\xi_{\bx'+\be(m',\bx')}(m'-1)=1$, which necessarily exists by construction. We construct a sequence of space-time points starting at $(m,\bx)$. Elements of the sequence such that both the previous and next element have larger time coordinate are called \emph{sinks}. Initially the sequence is the single point $\cT_0=(m,\bx)$. If $\be(m,\bx)$ is not defined, we stop and output $\cT_0$. Otherwise, we \emph{explore}: replace the point $(m,\bx)$ by the sequence $\cT_1=(m,\bx),(m-1,\bx),(m,\bx),(m-1,\bx+\be(m,\bx)),(m,\bx)$. For $l\ge 1$, given $\cT_l$, define $\cT_{l+1}$ as follows. Among the sinks $(m',\bx')$ of $\cT_l$ such that $\be(m',\bx')$ is defined, we find the first one that maximises $m'$. If no such element exists, we stop and output $\cT_l$. Otherwise, explore: replace the selected space-time point $(m',\bx')$ in $\cT_l$ by $(m',\bx'),(m'-1,\bx'),(m',\bx'),(m'-1,\bx'+\be(m',\bx')),(m',\bx')$. Denote the result of this exploration operation by $\cT_l'$. If the points $(m'-1,\bx')$ and $(m'-1,\bx'+\be(m',\bx'))$ are not already present in $\cT_l$, we set $\cT_{l+1}=\cT_l'$. If $(m'-1,\bx')$ appears in $\cT_l$, we remove from $\cT_l'$ all vertices of the corresponding sub-cycle except one of the two occurrences of $(m'-1,\bx')$ to obtain $\cT_l''$. Otherwise, set $\cT''_l=\cT'_l$. Finally, we do the same \emph{loop erasion} operation on $\cT''_l$, if $(m'-1,\bx'+\be(m',\bx'))$ appears twice in this sequence. The final result defines $\cT_{l+1}$. The output of this algorithm is called the \emph{Toom cycle rooted at} $(m,\bx)$. A number of combinatorial properties of this object are needed. All of their proofs are fairly simple, but fiddly, so we refer to \cite{Swart22} for the details. Firstly, the Toom cycle is well defined, contains its root and its increments belong to 
\[\{(-1,0,0),(1,0,0),(-1,1,0),(1,-1,0),(-1,0,1),(1,0,-1)\}.\]
Consequently, the number of Toom cycles with given root of length $l$ is at most $6^l$. Secondly, for each sink $(m',\bx')$, the space-time point $(m'-1,\bx')$ is not good and each sink appears only once in the cycle. 

Observe that there are three types of vertices in a Toom cycle: the ones with two neighbours with smaller time coordinate, those with two neighbours with larger time coordinate and the others. We call them \emph{sources}, \emph{sinks} and \emph{internal vertices}. The number $n_*$ of sources and the number of sinks are the same by double counting. One can prove that internal vertices $(m_i,\bx_i)$ such that their neighbours are not sinks or sources satisfy 
\[((m_{i+1}-m_i),\<\bone,\bx_{i+1}-\bx_i\>)=((m_i-m_{i-1}),\<\bone,\bx_{i}-\bx_{i-1}\>)\in\{(-1,0),(1,-1)\}.\]
When the above quantity is $(1,-1)$ (resp.\ $(-1,0)$), we call the internal vertex \emph{blue} (resp.\ \emph{red}) and denote the number of such vertices $n_b$ (resp.\ $n_r$). Examining the time increments, it is clear that $|n_b-n_r|\le 4n_*$. On the other hand, by examining the projected space increments $\<\bone,\bx_{i+1}-\bx_i\>$, we see that $n_b\le 6n_*$. Combining these facts, we obtain that the total length of the Toom cycle is at most $6n_*+n_b+n_r\le 22n_*$.

As noted above, sinks are not good and distinct, so the probability that all sinks of a Toom cycle are bad is at most $(\varepsilon')^{n_*}$. We are then able to use a union bound
\[\bbP\left(\xi_\bx(m)=1\right)\le\varepsilon'+\sum_{l\ge 1} 6^{2l}(\varepsilon')^{2l/22}<(\varepsilon')^{1/12},\]
for $\varepsilon'$ small enough, taking into account that cycles have even length.

\subsubsection{Chains of Toom cycles}
\label{subsubsec:chains}
In Section~\ref{subsubsec:Toom}, we showed that for any space-time point $(m,\bx)$, the probability $\bbP(\xi_\bx(m)=1)$ that it is in state 1 in North-East BP with death parameter $\varepsilon'$ and initial condition $\bzero$ is small. However, following Section~\ref{subsubsec:BPwD}, we need to prove that the probability of a connected component of such points is exponentially low in the diameter of the component. We already have one exponential bound, namely for a given space-time point, the probability that its Toom cycle has length more than $l\ge 1$ is at most $6^l(\varepsilon')^{l/22}$. It therefore remains to treat connected components of Toom cycles. Unfortunately, if Toom cycles for different points intersect at a sink, we lose independence. In order to deal with this issue, we introduce the following notion of chain (of Toom cycles), following \cite{Hartarsky24KCM-CP}*{Section~9.3}.
\begin{definition}[Chain]
    A \emph{chain} is a finite sequence of vertex-disjoint Toom cycles $T_i$ rooted at space-time points $(m_i,\bx_i)$ of lengths $l_i$ such that 
    \begin{equation}
    \label{eq:chain:distance}
    d((m_i,\bx_i),(m_{i+1},\bx_{i+1}))\le 7(l_i+l_{i+1}).\end{equation}
    The length of a chain is $\sum_i l_i$.
\end{definition}
The next result \cite{Hartarsky24KCM-CP}*{Lemma~9.15} shows that we can extract a chain from a large connected component of occupied sites.
\begin{lemma}[Existence of chains]
\label{lem:chain}
There exists a constant $C>0$ (depending only on the dimension, which we fixed equal to 2) such that the following holds. If the connected component $\cC_{m,\bx}'$ of $(m,\bx)$ in $X$ from \eqref{eq:m'x'} has diameter at least $k$, then there exists a chain of length $l\ge k/C$ contained in the ball $B_{(m,\bx)}(Cl)$ of radius $Cl$ centered at $(m,\bx)$.
\end{lemma}
\begin{proof}[Sketch]
    The construction of the chain is obtained algorithmically by pruning the sequence $(T_i)_{i=1}^N$ of Toom cycles rooted at the space-time points of a path from $(m,\bx)$ to the farthest point in its connected component $\cC_{m,\bx}'$. Initially the Toom cycles in this sequence may intersect. We start by discarding the Toom cycles rooted at space-time points at distance at most $6l_1$ from $(m_1,\bx_1)=(m,\bx)$ except $T_1$. We inspect the first remaining Toom cycle $T_i$ with $i>1$ (if any). If $T_1\cap T_i\neq \varnothing$, we also remove $T_1$ from our sequence, in which case we observe that necessarily $l_i> 3l_1$, so
    \[\left\{\left(m_j,\bx_j\right):j<i\right\}\subset B_{(m_1,\bx_1)}(6l_1)\subset B_{(m_i,\bx_i)}(6l_i).\]
    If, on the contrary, $T_1\cap T_i=\varnothing$, we keep both $T_1$ and $T_i$ for the moment and observe that $d((m_1,\bx_1),(m_i,\bx_i))\le 6l_1+\sqrt{2}$ by construction, since the roots initially form a connected path in the sense of \eqref{eq:connected:chain}. 
    
    In subsequent steps, we proceed similarly. Namely, we first remove the Toom cycles rooted at space-time points $(m_j\bx_j)$, with $j>i,$ at distance at most $6l_i$ from $(m_i,\bx_i)$. Then, for the first remaining Toom cycle $T_j$ with $j>i$, we remove those among the remaining ones with smaller index which intersect $T_j$. After the algorithm terminates, denoting the set of remaining indices of Toom cycles by $I$, one can show that $T_i\cap T_j=\varnothing$ for distinct $i,j\in I$ and
    \[\bigcup_{i\in I}B_{(m_i,\bx_i)}(6l_i)\]
    is a connected set containing the initial path. This allows us to extract the desired chain. See \cite{Hartarsky24KCM-CP} for more details.
\end{proof}
With Lemma~\ref{lem:chain} at hand, the proof of Theorem~\ref{th:mixing} is nearly complete. We use a union bound over the possible chains. Recall that Toom cycles in a chain are disjoint (so independent) and their probability is at most $(\varepsilon')^{l/22}$ for a chain length $l$. It therefore remains to show that the number of possible chains of length $l$ is bounded by $C^l$ for some constant $C$ independent of $\varepsilon'$. Indeed, the number of chains modulo the choice of the roots of the Toom cycles is at most $6^l$ as discussed above. Finally, the number of choices for the roots can be bounded using \eqref{eq:chain:distance} and the fact that $\bbZ^2$ does not have super-exponential growth.

Putting the above together, we obtain an exponential bound on $\bbP(\diam(\cC'_{m,\bx})\ge k)$, concluding the proof of Theorem~\ref{th:mixing} for FA-2f in two dimensions.

\section{Other out-of-equilibrium results}
\label{sec:out-of-eq:other}
We next briefly mention a few more works tackling KCM out of equilibrium from angles 
 different from those discussed above.

\subsection{The biased annihilating branching process}

The biased annihilating branching process (BABP) is an interacting particle system closely related to FA-1f, to which several works have been devoted in the 1990s \cites{Sudbury95, Sudbury97,Neuhauser93}. Here sites are updated to $0$ (resp.\ $1$) with rate proportional to  the number of neighbouring empty sites. More precisely BABP has generator \eqref{eq:generator} with $c_{\vec x}$ replaced by $\tilde c_{\vec x}$, where
\[\tilde c_{\vec x}(\omega):=\sum_{i=1}^d\sum_{\varepsilon\in\{-1,1\}}(1-\omega_{\vec x+\varepsilon\vec e_i}).\]
Like FA-1f, BABP 
is reversible w.r.t.\ $\mu$ and not attractive. Furthermore, for all $\omega\in\Omega$,
\[\tilde c_{\vec x}(\omega)=0\Leftrightarrow c_{\vec x}(\omega)=0\]
with $c_{\vec x}$ the FA-1 constraint function (recall \eqref{eq:def:cx} and Section~\ref{sec:rules}). Despite these similarities, BABP has some special features not shared by FA-1f, which make it more tractable. In particular, it enjoys a self-duality property (see \cite{Sudbury95}*{Eq.~34}) and quasi-duality with another model known as double flip process (see \cite{Sudbury97}*{Section 6}). Thanks to these features, which correspond to some special algebraic properties satisfied by the generator, Sudbury \cite{Sudbury97} established that for BABP in any dimension convergence to equilibrium (i.e.\ \eqref{hope_nu}) holds when $\nu=\mu_{q_0}$, provided $q_0\neq 0$. 
Unfortunately, this result does not seem to provide any insight or tool for proving \eqref{hope_nu} at all vacancy densities for FA-1f.

\subsection{FA-1f at low density}
\label{sec:low}The FA-1f model has been particularly investigated, since it is the most accessible model beyond the East one. It has the notable advantage of being non-cooperative in the sense that a single empty site can move around and allow the configuration to be resampled. Moreover, as we saw in Section~\ref{sec:FA1f}, empty sites move more or less like random walks, albeit possibly coalescing when they meet and occasionally branching. If one works in a suitably chosen finite volume with $q$ small and only follows the process over a relatively short amount of time, it can be possible to track these random walks and deal with delicate collision events. More precisely, consider a `critical' volume, that is, a torus (or any other bounded degree finite connected graph) $\Lambda=(\bbZ/n\bbZ)^d$ of cardinal $n^d$ of order $1/q$. Since the graph is connected, by Corollary \ref{cor:legal:path}, the FA-1f model is ergodic on $\Omega\setminus\{\bone\}$, so we exclude the $\bone$ configuration in the sequel.Typically, under $\mu_q$, $\Lambda$ only contains a bounded number of empty sites and they are well separated. If we follow the process up to a time horizon $1/q^2$, say, we do not expect to see more than three empty sites near each other. Since moving an empty site requires creating another one, this means that only binary collisions or branchings occur and these can be analysed.

Furthermore, the above reasoning can also work out of equilibrium, provided that one can control the length of the initial period of time it takes for the number of empty sites to drop down to order 1. This delicate strategy was employed by Pillai and Smith \cites{Pillai17,Pillai19} yielding that the mixing time in the above setting is of order $q^{-2}$, up to logarithmic corrections, for $d\ge 2$. As already discussed in the proof of Proposition~\ref{prop:CBSEP}, one can prove the same result, but also for other graphs by establishing a finite-volume logarithmic Sobolev inequality \cite{Hartarsky22CBSEP} (recall \eqref{eq:logsob}, Corollary~\ref{cor:LS}). In fact, such inequalities already proved useful in \cite{Blondel13} for studying FA-1f at high vacancy density on infinite graphs (e.g.\ $\bbZ^d$) with initial condition such that there are empty sites at bounded distance from all vertices of the graph. However, in the absence of attractiveness, it is usually difficult to relate finite and infinite volume results.

\subsection{Large deviations in trajectory space}
A radically different viewpoint consists in studying trajectories of KCM rather than their state at a given time. Given, e.g., a finite box $\Lambda=\{1,\dots,n\}^d$ with suitable boundary condition $\sigma\in\Omega_{\bbZ^d\setminus\Lambda}$ and a time $t$, the \emph{activity} $\cA(t)$ is the total number of times any site changed its state up to time $t$. One is interested in the large deviation properties of $\cA(t)$ as $t\to\infty$. That is, for $a\in[0,\infty)$, one expects that
\[f(a)=\lim_{\varepsilon\to0}\lim_{n\to\infty}\lim_{t\to\infty}\frac{-\log\bbP_{\mu}(|\cA(t)-a|<\varepsilon)}{n^dt}\]
exists and would like to study the properties of the large deviation rate $f$. The constraint of the KCM impacts the function $f$ in that $f(a)=0$ for all $a<\lim_{n\to\infty}\lim_{\to\infty}\bbE_{\mu}(\cA(t))/(n^dt)$: low activity is easy to achieve. This can be seen as a constraint-induced dynamical phase transition. We refer the reader to \cites{Garrahan09,Garrahan07,Bodineau12,Bodineau12a,Jack06}, where these and finer matters are investigated for FA-1f and East in one dimension.

\subsection{Aging for the one-dimensional East}
Remark \ref{rem:East}
guarantees that, for $q\ll 1$, relaxation to equilibrium after a density quench (namely starting from $\nu=\mu_{q_0}$ with a fixed $q_0\in(0,1]$) occurs at an exponential rate on a time scale of order $T_{\rm rel}\sim e^{c|\log q|^2}$ with $c=(2\log 2)^{-1}$ (see Theorem \ref{th:East}).
In this section we will discuss a peculiar behaviour that occurs  at intermediate times $1\ll t \ll T_{\rm rel}$. 
 
 Fix $\epsilon, q\in(0,1)$, and, for $n\ge 1$, set 
\begin{align}
t_0&{}=1,&t_0^-&{}=0,&t_0^+&{}=\left(\frac{1}{q}\right)^{\epsilon},\nonumber\\
t_n&{}=\left(\frac{1}{q}\right)^n,&t_n ^-&{}=t_n ^{1-\epsilon},&t_n
^+&{}=t_n ^{1+\epsilon}.\label{deftn}
\end{align}
The time intervals $[t_n^-, t_n^+]$ and $[t_n^+,t_{n+1}^-]$ are
called \emph{$n$\textsuperscript{th} active} and \emph{$n$\textsuperscript{th} stalling}
periods respectively. Note  that in the limit $q\to 0$ there is a sharp separation of time scales, $t_n/t_{n'}\to 0$ for $n<n'$ and, if $\epsilon\ll 1$, the $n$\textsuperscript{th} stalling period is much longer than the $n$\textsuperscript{th} active period.

Suppose that $q\ll 1$ and initialise the process from a Bernoulli distribution at density $q_0>q$. Most of the non-equilibrium
evolution will try to remove the excess of empty sites present initially and will thus be dominated by the coalescence of domains corresponding to the intervals separating two consecutive empty sites.
This process must  occur in a 
 cooperative way because, in order to
remove an empty site, another empty site must be created to its right. Furthermore, recalling that creating an empty site at distance $\ell$ requires creating $\sim\log_2 \ell$ additional empty sites (see Proposition \ref{prop:comb}) the  following heuristic picture emerges 
\begin{enumerate}
\item  During the $n$\textsuperscript{th} active period, only the empty sites with another empty site to their right at
  distance less than $2^n$ can be removed; \item at the end of the $n$\textsuperscript{th} active
  period, no empty sites at distance less than  $2^n+1$ are present any more, and there is no empty site that was not present at the beginning of the period;
  \item during the
  $n$\textsuperscript{th} stalling period, 
  nothing  happens: none of the empty sites present at the beginning is destroyed and no new empty site is present at the end.
\end{enumerate}
The above heuristics, set out  by Evans and Sollich in two  physics papers \cites{Sollich99,Sollich03}, was turned into rigorous results by Faggionato, Martinelli, Roberto and Toninelli  in \cite{Faggionato12}.
This implies  that the vacancy density displays a peculiar staircase behavior (see \cite{Faggionato12}*{Theorem 2.5(i)}) and
 the 
two-time autocorrelation function  depends in a non trivial way
on the two times
and not just on their difference (see \cite{Faggionato12}*{Theorem 2.5(i)}). This  
phenomenon is known as \emph{aging}.
In the same paper, sharp results for the statistics of the interval  between two consecutive zeros in the different periods were established (see\cite{Faggionato12}*{Theorem 2.6}). 
The key ingredients for these results are:
\begin{itemize}
    \item proving that  the non equilibrium dynamics of the East model starting from a renewal process is well approximated when $q\downarrow 0$ by
a  hierarchical coalescence process whose rates depend on large deviation probabilities
of the East model;
\item universality results for the scaling limit of this coalescence process (see \cite{Faggionato12a}). 
\end{itemize}

\begin{problem}
It could be interesting to investigate whether  a staircase behaviour for local functions and aging for two-time functions hold also for East in higher dimensions or for other KCM featuring a sharp separation of time scales, in particular for supercritical rooted models where logarithmic energy barrier also occur (see Proposition \ref{prop:Laure:combi}).
\end{problem}

\section{Basic open problems}
\label{sec:open-pbs}

The results discussed in the previous sections leave several basic questions open.

\subsection{Ergodic regime}
Fix a dimension $d\geq 1$, a $d$-dimensional update family $\cU$ and  $q>q_c(\cU)$. The following conjectures  can be viewed  either as an extension to the whole set of $\cU$ of the results of Section~\ref{sec:East-out} for the East model  or as extension to the whole ergodic regime of the high temperature results of Section~\ref{sec:high-temp}.
  
 \begin{conjecture}[Convergence to equilibrium after a temperature quench]
Let $\nu$ be the Bernoulli measure $\mu_{q_0}$, $q_0>\qc$. Then,
\begin{enumerate}
\item for any local function $f$, \eqref{hope_nu} holds;
\item
 furthermore, the convergence occurs exponentially fast.
\end{enumerate}
 \end{conjecture}

\begin{conjecture}[Linear pre-cutoff]
There exists $C$ such that\eqref{hope_tmix}
holds  for any $\delta\in (0,1)$ and $n$ large enough.\label{conj:tmix}
\end{conjecture}

\begin{conjecture}[Cutoff]
There exists \sl{${v}$} such that \eqref{hope_cutoff} holds.
\end{conjecture}

\subsection{Beyond the ergodic regime}
Fix a dimension $d\geq 1$, a $d$-dimensional update family $\cU$ and  $q<q_c(\cU)$.

\begin{problem}[Quench from high to low temperature]
Let the initial distribution $\nu$ be $\mu_{q_0}$ with $q_0>\qc$. Since a.s.\ the initial state has closure $\bzero$, this will be true at any later time (see Lemma \ref{lem:closure}).
However, it is natural to expect a behaviour dominated by the growth of occupied ``clusters'' that can be unblocked only from their boundary. Is this intuition correct and how fast do these clusters grow?  E.g.\ how does $\tbp$ of the current configuration scale with time? Another natural question is whether, locally, the measure converges to the equilibrium density, namely whether $\lim_{t\to\infty}\mathbb E_{\mu_{q_0}}(\eta_0(t))=1-q$. In principle, this is still compatible with the fact that the closure is $\bzero$ at any time. See however \cite{Cancrini10}*{Section 4}  for an example of KCM defined on trees (see Section \ref{sec:tree}) for which such a result is ruled out.
\end{problem}

The inverse regime, $q_0<\qc<q$, is often uninteresting, because the set of sites which can be updated may partition into an infinite collection of independent finite Markov chains. However, in some cases, an infinite chain remains. This leads to a setting similar to inhomogeneous KCM, discussed in Section~\ref{sec:general:KCM} and \ref{sec:inhomogeneous} below.

\section{Techniques}
Finally, let us review the new techniques encountered in this chapter. In fact, the only one that was not new (recall Section~\ref{subsec:tools:BP}) is renormalisation used in the proof of the general case of Theorems~\ref{th:convergence:out-of-eq} and~\ref{th:mixing}.

\noindent\textbf{Distinguished empty site.} Definition~\ref{def:distinguished} provided this essential tool for the analysis of the one dimensional East model. It  can be seen as a  boundary condition moving to the right and leaving equilibrium to its left.

\noindent\textbf{Couplings.} In Sections~\ref{subsubsec:CP}-\ref{subsubsec:BPwD}, we saw useful couplings between KCM and various attractive processes: cooperative contact processes, BP with death and last passage percolation, all based on the graphical construction of the models. The drawback of such comparisons is that, in order for them to be useful, the attractive process may need to be supercritical, leading to an artificial limitation to the high vacancy density regime.

\noindent\textbf{Toom cycles.} This combinatorial tool for proving stability of BP with death, discussed in Section~\ref{subsubsec:Toom}, has broader applications. Namely, suitable extensions can be used for perturbations of arbitrary cellular automata and, in some cases, interacting particle systems \cite{Swart22}.

\noindent\textbf{Chains.} In order to control connected components rather than individual sites, in Section~\ref{subsubsec:chains}, we relied on chains. This is a rather general percolation method not specific to Toom cycles.

\noindent\textbf{Cutoff strategy.} In both Corollary~\ref{th:East:cutoff} and Theorem~\ref{th:FA1f:cutoff} the same strategy was used to prove a linear time cutoff with square root window. While the entire programme has only been implemented for specific models in one dimension, it is, in principle, feasible for any KCM. The first step is proving a (possibly stretched) exponential convergence to equilibrium as in Theorem~\ref{th:convergence:out-of-eq}. Then one shows a positive speed result as in Theorem~\ref{th:mixing}. The third step is showing ergodicity for the process seen from the front and a law of large numbers for the front position. Finally, a central limit theorem for the front can be sought in order to control the cutoff window.


	\chapter{Related settings and models}

\label{chap:other}
  
\abstract{In this final chapter, we conclude by surveying a few additional settings not covered by the models defined in Chapter~\ref{chap:models}, but strongly related. Indeed, so far we have allowed general update families and, in Chapter~\ref{chap:out}, non-equilibrium initial conditions. However, we have restricted our attention to the equilibrium measure $\mu$ being product, constraints being identical at all sites, dynamics changing the state of a single site at a time, the underlying graph being a $d$-dimensional lattice. Each of these hypotheses may be revoked and leads to interesting models and questions, many of which have not yet been explored.}

\section{KCM on other graphs}
\label{sec:tree}
KCM can also be defined on graphs different from $\mathbb Z^d$, including arbitrary graphs \cites{Cancrini09,Hartarsky22CBSEP}, trees \cites{Cancrini15,Martinelli13,Cancrini10}, hyperbolic lattices \cite{Sausset10} and many more such as Bienaym\'e--Galton--Watson trees or various models of random graphs waiting to be explored. See e.g.\ \cite{Aldous05} for a possible application of FA-1f to information storage in a sensor network. The most studied case is FA-$j$f on oriented or non-oriented regular trees of degree $k+1$. In the non-oriented version, the constraint at site $\bx$ is satisfied if at least $j$ of the neighbours of $\bx$ in the tree are empty. In the oriented version, at least $j$ empty sites should be among the $k$ children. As for KCM on $\mathbb Z^d$, the  ergodicity thresholds for the KCM and for the corresponding BP dynamics coincide (see Theorem \ref{th:ergo}). Thanks to the tree structure, it is not difficult to write recursive equations for the critical thresholds \cite{Balogh06a}, yielding that,
both for oriented and non-oriented models, we have $\qc=1$ for $j> k$, $\qc=0$ for $j=1$  and  $0<\qc<1$ for $j\in\{2,\dots,k\}$. In  \cite{Martinelli13}, Martinelli and Toninelli prove that $\trel<\infty$ for all models in the whole ergodic regime, $q>\qc$.  In \cite{Cancrini15}, the same authors together with Cancrini and Roberto, analyse the scaling in the critical regime in the case $j=k$ for the oriented
model and prove a power law divergence of $\trel$  as $q\downarrow \qc$. An analogous scaling is also conjectured in all cases $j\in\{2,\dots,k-1\}$, but has yet to be proven. A fundamental difference between the case $j=k$ and $j\in[2,k-1]$, at the root of the difficulties to handle the missing cases, is that the BP transition is continuous for $j=k$ (it essentially correspond to a standard percolation transition), while for all $j\in\{2,\dots,k-1\}$  it is discontinuous, namely $\mu_{\qc}(\tbp=\infty)>0$. Using a strategy similar to the proof of Theorem \ref{theo:East_out}, for the oriented
model with $k=j=2$, exponential convergence to equilibrium when $q>\qc$ starting from an initial distribution $\mu_{q'}$ with $q'>\qc$ is proven in \cite{Cancrini10} (see Theorem 4.3 therein together with \cite{Martinelli13}*{Theorem 2}). The proof can be readily extended to all oriented models for $j=1$ or $j=k$. The result is conjectured to hold also for the remaining cases $j\in[2,k]$.

\section{Inhomogeneous KCM}
\label{sec:inhomogeneous}
We briefly considered an inhomogeneous setting in Section~\ref{sec:general:KCM}, where the update family defining the constraint is allowed to depend on the site. While our treatment in one dimension is very general, one can also consider such inhomogeneous models in higher dimensions, possibly choosing the update families at random. For instance, one may consider sites $\bx\in\bbZ^d$ having FA-$j_\bx$f constraint with $j_\bx$ chosen i.i.d.\ at random according to some distribution on $\{0,\dots,2d\}$. A few such models are studied in \cites{Shapira20a,Shapira19polluted}, but most remain unexplored.

\section{KCM with interactions}
\label{sec:interaction}
Another natural modification of KCM is to introduce static interactions between occupied sites. This may be achieved by updating each site w.r.t.\ a measure depending on the current state of other sites. For instance, one could consider the $\cU$-constrained Glauber dynamics for the Ising model with inverse temperature $\beta$ with generator \eqref{eq:generator}, where, instead of $\mu_\bx(\omega_\bx=0)=q$, we set 
\[\mu_\bx(\omega_\bx=0)=\frac{1}{1+\exp(\beta\sum_{\by\sim\bx}(2\omega_\by-1))},\]
the sum running over nearest neighbours of $\bx$ in $\bbZ^d$. In fact, this was already considered in \cite{Fredrickson84} together with an external magnetic field. More generally Gibbs measures were considered in \cite{Cancrini09}*{Section~5}. While the initial motivation behind KCM is to investigate the extent to which glassy phenomenology can be explained by purely dynamical means, in reality interactions are certainly present. It is therefore interesting, but probably challenging, to study such models.

\section{Plaquette models}
\label{sec:plaquettes}
Instead of kinetic constraints, \emph{plaquette models} have static interactions as in Section~\ref{sec:interaction}, but of multi-body type. They were introduced to show that kinetic constraints can emerge from static interactions at low temperatures {\cites{Garrahan02,Newman99,Turner15}}. An example is the \emph{square plaquette model} on $\bbZ^2$. For any $\bx\in\bbZ^2$, the plaquette of $\bx$ in this model is the square $P_\bx=\{\bx,\bx+\be_1,\bx+\be_2,\bx+\be_1+\be_2\}$. The Gibbs weights are defined by the Hamiltonian \[-\sum_{\bx\in\bbZ^2}\prod_{\by\in P_\bx}(2\omega_\by-1)\]
and one then considers an unconstrained single site Glauber dynamics. The behaviour of the square plaquette model turns out to be similar to that of FA-1f, while a similar triangular plaquette model with plaquettes of the form $P_\bx=\{\bx,\bx+\be_1,\bx+\be_1+\be_2\}$ is conjectured to behave like the East KCM. Work on these models can be found in \cites{Chleboun17,Chleboun20,Chleboun21}.

\section{Conservative models}
\label{sec:conservative}
The physical motivation behind KCM (recall Chapter~\ref{chap:intro}) views sites of $\bbZ^d$ as mesoscopic volumes whose particle density may be lower or higher, as reflected by the state of the site. However, if we take a microscopic perspective, it is more natural to consider constrained models in which the number of particles is conserved. The first and most classical such models are the Kob--Andersen ones \cite{Kob93}. In KA-$j$f, one may exchange the states of any two neighbouring sites, $\vec x$ and $\vec y$, provided they both have at least $j-1$ empty neighbours in $\bbZ^d\setminus\{\vec x,\vec y\}$. The case $j=1$ coincides with the simple symmetric exclusion process (SSEP) and we will disregard it in the following. 
One can similarly define conservative versions of $\cU$-KCM, obtaining the class of \emph{kinetically constrained lattice gases} (KCLG) (see \cite{Cancrini10a}*{Section 2} for a formal definition of this class). It is immediate to verity that for any $q\in[0,1]$, $\mu_q$ is a reversible measure for these dynamics. As for the non conservative models, we define $\qc$ by \eqref{eq:def:qc}. 
In \cite{Cancrini10a}*{Proposition 2.16}, a  conservative analogue of Theorem~\ref{th:ergo}
is established stating that for any KCLG the ergodicity threshold  coincides with the \emph{excheangeability threshold} defined as the minimal value of $q$ above which, for $\mu_q$-almost every configuration and for any couple of
sites, there exists a legal path exchanging their occupation variables.
The first mathematical result on KCLG is due to Biroli, Fisher and Toninelli \cites{Toninelli04, Toninelli05} who proved  that $\qc=0$ for any  KA-$j$f with $j\in\{2,\dots,d\}$.\footnote{The cases $j>d$ 
trivially lead to $\qc=1$.} The key ingredients of this proof are: 
\begin{itemize}
    \item the  construction of a set of configurations, the so called \emph{frameable configurations}, which can be connected by a legal path to a configuration with well-chosen boundary state (see \cite{Martinelli20}*{Definition 3.3} and \cite{Toninelli05}) for the dynamics on finite volume with occupied boundary condition;
    \item the construction of legal paths connecting any two frameable configurations with equal number of empty sites;
    \item the fact that the $\mu_q$-probability that a configuration is frameable goes sufficiently fast to one as the volume increases (see \cite{Martinelli20}*{Proposition 3.26}). More precisely
    , there exist $c_+,c_->0$ dependent on $d,j$ such that, setting
 $$\Xi_{\pm}(q,j,d)=\exp^{\circ (j-1)}\left(\frac{c_{\pm}}{q^{\frac{1}{d-j+1}}}\right),$$
 when $q\to 0$ and $L\to\infty$ faster (resp.\ slower) than  $\Xi_+(q,j,d)$ (resp.\  $\Xi_-(q,j,d)$), the probability of being frameable goes to one (resp.\ to zero).
\end{itemize}
Combining these results with canonical paths, renormalization arguments and tools borrowed from oriented percolation, Martinelli, Shapira and Toninelli \cite{Martinelli20}*{Theorem 1} prove  that the spectral gap of the KA-$j$f models on $\{1,\dots,L\}^d$ with unconstrained sources at the boundary in any dimension $d$ and for any $j\in [2,d]$ scales diffusively (as for SSEP), namely as $L^{-2}$, and with a density pre-factor of the form  $\Xi_{\pm}(q,j,d)^{-1}$ (while for SSEP there is no density pre-factor).

In \cite{Blondel18}, Blondel and Toninelli consider the behavior of a  tagged particle, and prove (following the ideas sketched in \cite{Toninelli04}), that for all $d\geq 2$ and for any $j\in [2,d]$,     diffusive behaviour holds at any density $q\in
(0,1)$. More precisely, 
if  we distribute the configuration  according to $\mu_q$, condition on the presence of a particle initially at the origin, tag it and denote by $\vec X_t$ its position at time $t$, for some matrix $D(q)$ such that $\vec e_i \cdot D(q) \vec e_i>0$ for all $i\in\{1,\dots,d\}$, it holds that
\begin{equation}
\label{eq:diff}
\lim_{\epsilon\to 0} \epsilon \vec X_{\epsilon^{-2}t}=\sqrt{2D(q)}B_t,
\end{equation}
where $B_t$ is a standard $d$-dimensional Brownian motion and the convergence holds in the sense of weak convergence of path measures.
This contradicts conjectures based on numerical simulations \cites{Kob93,Kurchan97} claiming the occurrence of a diffusive/non-diffusive transition. Furthermore Ertul and Shapira prove (see \cite{Ertul21}*{Theorem 2.3}) upper and lower bound for $D(q)$ of the form $\Xi_{\pm}(q,j,d)^{-1}$ (modulo a logarithmic correction  in the case $d=2$). The fast shrinking to zero  explains why  it was incorrectly inferred from numerical simulation (\cite{Kob93}) that for $d=3$ and $j=3$ a diffusive/non-diffusive transition would occur.

\begin{problem}
\label{prob:conservative}
Two other natural issues in the conservative setting are  \begin{enumerate}
    \item \label{item:prob1}
    determining the evolution of macroscopic density profiles, namely establishing the hydrodynamic limit, and the fluctuations around these profiles;
    \item \label{item:prob2}establishing relaxation at equilibrium in infinite volume.
\end{enumerate}
Concerning \ref{item:prob1}, a natural candidate for the hydrodynamic limit is a parabolic equation of porous media that degenerates when the density approaches one. As for fluctuations, it is reasonable to expect they would be Edward--Wilkinson Gaussian fluctuations as for SSEP \cite{DeMasi84}.
Establishing these results in the presence of constraints is particularly challenging (see \cite{Goncalves09} where this is achieved for a different KCLG). Concerning \ref{item:prob2}, a first result \cite{Cancrini10} shows that there exists $C(q)>0$ such that for any local function $f$  it holds that, for all $t>0$,
\[\var_{\mu_q} (P_t f) \leq \frac{C(q) \|f\|_{\infty}}{t}.\]
We expect the correct behavior to be of the form $t^{-d/2}$, as for SSEP. 
\end{problem}

Other kinetically constrained lattice gases can be found in \cites{Bertini04,Goncalves09,Bonorino20,Goncalves11,dePaula21,Jackle91a,Ertel88,Nagahata12,Shapira23}. A model that is currently very actively considered is known as the \emph{facilitated exclusion process}. In this one-dimensional model, a particle is allowed to jump to a neighbouring empty site, provided its other neighbour is occupied. We direct the reader to \cites{Blondel21,Erignoux24,Erignoux24a,Ayre24,DaCunha24,Erignoux23,Goldstein21,Blondel20,Massoulie24} for work on this topic.

\section{Tracer diffusion}
\label{sec:tracer}
In the previous section we described the tagged particle behavior for KCLG. Though the non-conservative dynamics of KCM is not diffusive, one can define the following similar problem. 
Consider  a stationary KCM evolving from a configuration distributed according to $\mu_q$  and inject at time zero a particle (the tracer) at the origin. The tracer moves like a modified random walk attempting to jump at rate one to a site chosen uniformly at random among its nearest neighbours, with the jump being allowed if and only if both the sites occupied by the walker before and after the move are empty (see \cite{Blondel14} for a precise definition).
Note that the KCM constitutes a dynamical random 
environment in which the tracer evolves, and is not influenced by the motion of the tracer.
Blondel proved (see \cite{Blondel14}*{Proposition 3.1 and 3.2}) that if the underlying KCM has a positive spectral gap, the tracer has a diffusive behavior with a non-degenerate diffusion matrix, namely  \eqref{eq:diff} holds. 
Furthermore for the FA-1f model in any dimension $d\geq 1$, it holds that $\vec e_i \cdot D(q)\vec e_i\sim q^2$ for 
$q\downarrow 0$ and $i\in\{1,\dots d\}$. Instead, for the East model \cite{Blondel14}*{Theorem 3.5} proves that $D$ scales as the spectral gap (modulo power law corrections in $q$). This corrects the conjecture that had been put forward by physicists (\cites{Jung04,Jung05}) affirming that for the East model $D$ would scale as $\trel^{-\xi}$ 
with $\xi<1$.

Asymmetric tracers on stationary KCM have also been the object of investigation. See for example \cite{Avena16} for results on a tracer on the one-dimensional East model and with positive (resp.\ negative) drift  when  on its current position the occupation variable of the East model is occupied (resp.\ empty).

\section{Upper triangular matrix walk}
We conclude with a further context in which KCM arise naturally beyond the study of glassy dynamics and interacting particle systems. Let $G_n$ be the group of $n\times n$ upper triangular matrices with entries in the two-element field $\bbF_2$ and ones on the diagonal. The following Markov chain was considered in \cite{Coppersmith00}. At each step, with probability $1/2$ nothing happens and, for each $i\in\{1,\dots,n-1\}$, with probability $1/(2n-2)$, we add row $i+1$ to row $i$. This corresponds to performing a lazy random walk on $G_n$ with generator set $(I+E_{i,i+1})_{i=1}^{n-1}$. If we restrict our attention to column $j$ of the matrix, this Markov chain becomes exactly the East process with parameter $q=1/2$ on the segment $\Lambda_j=\{1,\dots, j-1\}$ with boundary condition $\bzero_{\bbZ\setminus\Lambda_j}$. We direct the reader to \cites{Peres13,Ganguly19} for works on this random walk. 


\cleardoublepage
\backmatter
\phantomsection
\addcontentsline{toc}{chapter}{References}
\bibliography{Bib}

\end{document}